\documentclass[11pt]{article}
\usepackage{comment}
\usepackage[normalem]{ulem}
\usepackage{amsfonts}
\usepackage{amssymb,amsmath}
\usepackage{amsthm}
\usepackage{latexsym}
\usepackage{graphics}
\usepackage{epic}
\usepackage{epsfig}
\usepackage{psfrag}
\usepackage{subfigure}
\usepackage{bbm}
\usepackage{tikz}
\usepackage{dsfont}
\usepackage{enumitem}
\usepackage{stmaryrd}
\usepackage{xcolor}
\usepackage{hyperref} 
\usepackage{soul} 
\usepackage{subcaption}

\usepackage[margin=3cm]{geometry}

\numberwithin{equation}{section}

\def\E{\mathbb{E}}
\def\P{\mathbb{P}}
\def\R{\mathbb{R}}

\def\1{\mathbbm{1}}
\def\N{\mathbb{N}}

\def\eps{\varepsilon}

\def\dimE{q}

\newtheorem{thm}{Theorem}[section]
\newtheorem{lem}[thm]{Lemma}

\newtheorem{dfn}{Definition}
\newtheorem{prop}[thm]{Proposition}

\newtheorem{hyp}{Assumption}

\theoremstyle{remark}
\newtheorem{rem}{Remark}

\usepackage{geometry}
\usepackage[textwidth=2.5cm]{todonotes}
\presetkeys{todonotes}{size=\scriptsize}{}

\begin{document}

\title{Niche construction as an emerging phenomenon from fast ecological and slow evolutionary-environmental timescales in individual-based models} 
\author{Nicolas Champagnat$^1$, Coralie Fritsch$^1$, Pablo A. Marquet$^{2,3,4,5}$, \\Cristóbal Quiñinao$^2$, Leonardo Videla$^6$, Nicol\'as Zalduendo$^6$}

\footnotetext[1]{Universit\'e de Lorraine, CNRS, Inria, IECL, F-54000 Nancy, France.}
\footnotetext[2]{Facultad de Ciencias Biol\'ogicas, Pontificia Universidad Cat\'olica de Chile, Avenida Libertador Bernardo O’Higgins 340, 8331150, Santiago, Chile.}
\footnotetext[3]{The Santa Fe Institute, Santa Fe NM 87501, USA.}
\footnotetext[4]{Centro de Modelamiento Matem\'atico (CMM), Universidad de Chile, International Research Laboratory 2807, CNRS, CP 8370456 Santiago, Chile.}
\footnotetext[5]{Wissenschaftskollege zu Berlin, Wallotstrasse 19, 14193, Berlin, Germany.}
\footnotetext[6]{Universidad de Santiago de Chile, Facultad de Ciencia, Departamento de Matem\'atica y Ciencia de la Computaci\'on, Las Sophoras 173, Estaci\'on Central, Santiago, Chile.\\
E-mail: nicolas.champagnat@inria.fr, coralie.fritsch@inria.fr, pmarquet@uc.cl, cquininao@uc.cl, leonardo.videla@usach.cl, nicolas.zalduendo@usach.cl}

\maketitle

\begin{abstract}
  We consider an individual-based model of niche construction built upon birth-death processes of $d$ interacting biological entities (e.g., species, subpopulations or mutant strains in the same species) immersed
  in a slowly-varying environment which is influenced by the population state. In this setting, extinction and/or
  emergence of entities with negligible size or rarity, on long time scales can be observed. We prove that the joint dynamics of the environment and
  the species sizes on a logarithmic scale undergo a piecewise deterministic Markov process in a large time regime, which can be
  approximated by an explicit dynamical system in the limit of large populations. This convergence result holds true under the
  assumption that no ``bad event'' occurs during the population dynamics, corresponding to situations where several species become
  dominant or go extinct simultaneously. We also prove that no such bad event occurs for Baire-almost all initial conditions of the
  population sizes and the environmental resources. We apply this framework to study the long term coexistence of two specialist
  species consuming two resources, with a ``joint'' niche construction where each species constructs the niche of the other while
  depleting its resources. We also study an example of immune escape in cancer, where the environmental variable is associated with the
  state of the immune system and the species are associated with three types of tumor cells.
\end{abstract}


\pagebreak
\tableofcontents
\section{Introduction}
\label{sec:intro}

Niche construction is a process whereby organisms actively modify their surrounding environment~\cite{laland1996,odling-smee2003}. Modifications of the environment can occur through metabolism, engineering activities that create habitat for other species (e.g.\ burrows or nests) or habitat choice, and can change abiotic or biotic features of the environment. Modification of the ecosystem usually generates ecological and evolutionary feedback on the species, so the term of niche construction generally also refers 
to the evolutionary process that generates selection through the feedback produced by environmental changes~\cite{lewontin1983}. 

The evolutionary dynamics of niche construction results from the combination of three mechanisms: ecology (i.e.\ population dynamics), environmental change, and evolution. Each of these mechanisms can take place on different timescales~\cite{nguyen2021}. For example, in terrestrial ecosystems, plants can alter various soil properties (e.g.\ fertility, temperature, humidity) and thus influence nutrient cycling~\cite{knops2002}. This environmental change can occur on the ecological timescale for instance by uplift previously inaccessible soil resources using deep roots~\cite{chapin1997}, or it can occur on longer timescales, for instance by modifying nutrient mineralization, either through their litter quality~\cite{pastor1984} or by creating favorable abiotic conditions for decomposers under their canopy~\cite{cross1999}. In some cases, the environmental change can take place on short timescales, but it may impact the community for much longer timescales through ecological legacies. For instance, a termite mound may grow as the termite colony grow along the ecological timescale, but the changes of the mound could then affect local environment for millennia, creating new sources of selection that act on very long niche construction timescales~\cite{martin2018}. In other cases, the timescale of environmental change can be very long, 
such as for crustacean species whose shells accumulate beneath the ocean, gradually forming hard substrata which influences the occupation of other species~\cite{kidwell1983}. Other examples are given by changes of cultural practice in human communities which led to adaptation, such as the co-evolution of genes for adult lactase persistence and milk use, or plant cultivation and animal farming, which favored the spread of diseases such as malaria and sickle-cell anemia, leading to the selection of alleles that confer resistance to these diseases~\cite{durham1991,laland2010}. Although it is known in these cases that niche construction led to faster adaptation, it is unknown whether the changes of cultural practices occurred fast or progressively on the evolutionary timescale.

Most niche construction models assume that environmental changes occur on the same timescale as ecology, as in the seminal work of  references~\cite{laland1996,laland1999}, and usually incorporate evolutionary effects on a longer timescale, typically considering the environmental state as a parameter (see e.g.\ the survey~\cite{laland2010} on models of niche construction shaped by human cultural practices) or assuming a quasi-equilibrium (or adaptive dynamics) setting~\cite{kylafis2008}. New models were recently proposed to deal with more complex timescales interaction, considering either that environmental change occurs on a very slow evolutionary timescale~\cite{tanaka2020} (the timescale of adaptive dynamics~\cite{metz1996,geritz1998,champagnat2022}), or that evolution is fast~\cite{nguyen2021}. 

In this work, we propose a general methodology for niche construction modeling in cases where (part of) the environmental variables evolve on the evolutionary timescale, which is assumed to be slow compared to the ecological timescale, but not much slower (as opposed to~\cite{tanaka2020}). Our goal is to provide a macroscopic model able to approximate on the evolutionary timescale birth-death dynamics while accounting for extinction or emergence of negligible subpopulations on the same timescale. We assume that the number of different biological entities, which can be interpreted as different species in interaction and/or different subpopulations or mutant populations in the same species, is fixed and finite. The order of magnitude of the total population size will be parametrized by a scaling parameter $K$, related to the so-called carrying capacity of the ecological community under consideration. Assuming that the environmental variables evolve according to a dynamical system depending on the state of the population on the timescale $O(\log K)$ (long compared to the ecological timescale $O(1)$ of birth and death events), our main result gives the full characterization of the asymptotic behavior of the population dynamics when $K\to+\infty$ on the long timescale $O(\log K)$. In particular, we will characterize the dominant subpopulation (of size of the order of $K$) at each time $t$ (on the long timescale, i.e.\ after normalizing time by $1/\log K$) and approximate the size of each of the other subpopulations as $K^{\beta(t)}$ for some deterministic exponent $\beta:\R_+\to[0,1]$, where $\beta(t)=0$ means that the subpopulation is extinct and $\beta(t)\in(0,1)$ means that the subpopulation is non-extinct but negligible compared to the dominant subpopulation. This will allow us to characterize the extinction time of a subpopulation as the first time when its exponent $\beta(t)$ hits $0$, and its emergence time (i.e.\ the time at which it becomes dominant) as the first time where $\beta(t)=1$. The joint dynamics of the environmental variables and the exponents $\beta(t)$ of each subpopulation will be given by a time-discontinuous dynamical system where discontinuities occur at each time of change of dominant population, and each exponent evolves according to the fitness vector field (defined as the growth rate of the corresponding non-dominant subpopulation) determined by the dominant subpopulation and the environmental state. Our population scaling is similar to the one of the references~\cite{durrett-mayberry-2011,bovier2019crossing,coquille2021stochastic,champagnat2021stochastic,esser2024,brouard2024,costa2026} which studied evolutionary dynamics without environmental variables. In~\cite{esser2025,esser2025b}, a similar scaling was applied to an individual-based model with a periodically changing environment which does not depend on the population state and evolves on a fast timescale. The novelty of our work is the inclusion of an environmental dynamics which evolves jointly with the population on the evolutionary timescale. 

The limit dynamics in our main result is only obtained in cases where no bad event occurs, which would prevent us from characterizing the limit with our method. Such events can occur either when two or more negligible or rare subpopulations become dominant simultaneously, when one subpopulation becomes extinct exactly when another one becomes dominant, or if the fitness of a negligible subpopulation vanishes exactly at the same time at which it becomes dominant or extinct. The first two cases were already known in former works~\cite{coquille2021stochastic,champagnat2021stochastic}, the third one is new and is due to the fact that fitness depends on the changing environment. None of the previous works on similar scalings provided a characterization of the initial conditions and the model parameters which lead to such a bad event. Here, we provide the first result of this kind, showing that for generic model parameters (i.e.\ ``typical'' model parameters in a sense that will be made precise later) with two entities or subpopulations and all initial conditions in a Baire co-meager set, no bad event occurs and the convergence to the limit dynamics holds true for all positive time. We also extend these results to the case of three or more entities or subpopulations.

The joint dynamics of the environmental variables and the exponents $\beta(t)$ of each subpopulation describes a niche construction process where adaptation corresponds to the emergence of new subpopulations. The fact that the environment changes on the evolutionary timescale induces a change of the fitness landscape responsible for this adaptation, and considering not too slow evolutionary timescales (compared to adaptive dynamics studies of niche construction~\cite{kylafis2008,tanaka2020}) induces a progressive change in the size of the rare subpopulations which leads either to its extinction or dominance  on the same timescale. In this sense, our method allows to design macroscopic models which encompass the standard feature of many examples of niche construction, where adaptation and environmental changes both should occur on a long timescale. To illustrate this, we study in detail two simple examples of niche construction. The first one is a model of indirect mutualism, where 
two specialist species consume two resources, inspired by~\cite{macarthur1972,levine1976}, with a ``joint'' niche construction where each species constructs the niche of the other while depleting its resources. We show that, depending on the initial condition, coexistence in a (quasi-)periodic dynamics, or extinction of both species, one on a short timescale and the second one on a long timescale, may occur. The second example deals with immune escape in cancer inspired by~\cite{tanaka2020} and shows that, depending on the parameters, the immune system can maintain the tumor in a controlled state, or induce a periodic chase of emerging cancer variants. 

The paper is organized as follows: in Section~\ref{sec:BDP}, we define our individual-based model and state our main results: the
convergence of the population process towards a deterministic niche-construction model with extinction and emergence of
subpopulation, and the fact that this convergence holds true for generic models and initial conditions. We then study in
Section~\ref{sec:examples} our two examples. Preliminary results useful for the proof of our main results are given in
Section~\ref{sec:prelim}, and the proof of the convergence of the environment and the population exponents is given in
Section~\ref{sec:proof-main}. Finally, the proof that the limit dynamics does not encounter bad event for generic initial conditions is given
in Section~\ref{sec:induction-never-stops}. 

\section{Niche construction model and main results}
\label{sec:BDP}
In this section, we first describe our model and the main assumptions in Subsection~\ref{sec:model}, then we state our main result in Subsection~\ref{sec:main} and we explain the strategy of proof in Subsection~\ref{sec:main-2}. 

\subsection{A birth-death model of niche construction}
\label{sec:model}

We consider $d$ populations, that we call species hereafter, but which can be either seen as different species or groups of individuals with different traits or different phenotypes/genotypes in one or several species. For all $i\in\{1,2,\ldots,d\}$, we denote by $N^K_i(t)\in\mathbb{N}$ the number of
individuals of species $i$ at time $t\geq 0$, and we denote indifferently the full population state either as the vector $N^K(t):=(N^K_1(t),\ldots,N^K_d(t))$, or as the counting measure on the set $\{1,2,\ldots,d\}$, 
\[
  \nu^K_t=\frac{1}{K}\sum_{i=1}^d N^K_i(t)\delta_i,
\]
where $\delta_i$ is the Dirac mass at $i$. One notation or the other will appear as more convenient depending on the situation. Note that the counting measure notation is standard for individual-based models~\cite{fournier-meleard2004}. 

The parameter $K\geq 1$ in the previous definition denotes a scaling parameter, which is related to the carrying capacity of the ecosystem, in the sense that the total population size
\[
  K\langle\nu^K_t,\mathbf{1}\rangle=\sum_{i=1}^d N^K_i(t)
\]
will be of the order of $K$ when $K\to+\infty$. 

We also consider an environmental variable $E^K(t)\in\mathbb{R}^\dimE$ for some $\dimE\geq 1$, that may represent, for example, substrate
concentration, temperature, or even other species that can be regarded as resources, pathogens or predators. \emph{In this work we regard this environmental variable as continuous and evolving continuously}, although our results can be extended to other cases such as discrete components evolving through births and deaths in the vector $E^K(t)$ (see Remark~\ref{rem:extension} below).

We assume that births and deaths in the population occur one by one in continuous time, and that the birth rate per individual of species $i$ at time $t$ is given by $b_i(\frac{1}{K}N^K(t),E^K(t))$ (which we will also denote $b_i(\nu^K_t,E^K(t))$ when needed by abuse of notation) and the death rate per individual of species $i$ at time $t$ by $d_i(\frac{1}{K}N^K(t),E^K(t))$ (or $d_i(\nu^K_t,E^K(t))$), for some functions $b_i$ and $d_i$ from $[0,+\infty)^d\times\mathbb{R}^\dimE$ to $\mathbb{R}_+$. 

We also assume that the environmental variable follows the dynamics (although other dynamics, e.g.\ stochastic with jumps, could be considered and would not change the spirit of our results, see Remark~\ref{rem:extension} below):
\begin{equation}
  \label{eq:environment}
  \dot{E}^K(t)=\frac{1}{\log K}F\left(\frac{1}{K}N^K(t),E^K(t)\right),
\end{equation}
for some vector field $F:[0,+\infty)^d\times\mathbb{R}^\dimE\to\mathbb{R}^\dimE$. This makes the process $(\nu^K_t,E^K(t))_{t\geq 0}$ a
piecewise deterministic Markov process (PDMP)~\cite{davis1984piecewise} under mild assumptions, that will be implied by
Assumption~\ref{hyp:regularity-Ltka-Volterra} below (see Section~\ref{sec:prelim} for more details on the construction of the
process).

Note that we introduce a time-scaling of $1/\log K$ which slows down the dynamics of the environmental variable. This is in line with the discussion in Section~\ref{sec:intro}, namely that niche construction is the result of a slow modification of the environment together with a progressive adaptation of the population to this changing environment. For simplicity, we assume that mutations are unable to introduce new populations. All these restrictions of our model can be extended and generalized using standard techniques; see Remark~\ref{rem:extension} below. 

\bigskip

In this work, we focus on the competitive Lotka-Volterra density-dependence relationship
\[
b_i(\nu,e)=b_i(e)\quad\text{and}\quad d_i(\nu,e)=d_i(e)+\int_{\cal X} c_i(e,y)\nu(dy).
\]
Again, it is possible to generalize this density-dependence, see Remark~\ref{rem:extension} below.
Notice that in this setting, if we evaluate $d_i$ on the empirical measure $\nu_t^K$ and the environment $E^K(t)$, we are understanding that
$$
d_i(\nu^K_t,E^K(t)) 
=d_i(E^K(t))+\frac{1}{K}\sum_{j=1}^d c_i(E^K(t),j)N^K_j(t),
$$
or simply
$$
d_i(\nu,e)=d_i(e)+\frac{1}{K}\sum_{j=1}^d c_{i}(e,j)n_j
$$
for $\nu=\frac{1}{K}\sum_{j=1}^d n_j\delta_j$
\medskip

Our first assumption, that we will assume throughout the article, deals with the regularity of the model parameters.
\begin{hyp}
  \label{hyp:regularity-Ltka-Volterra}
  The maps $e\mapsto b_i(e)$, $e\mapsto d_i(e)$ and $e\mapsto c_i(e,j)$ are $C^1_b(\mathbb{R}^\dimE)$, i.e.\ bounded, continuously
  differentiable with bounded derivatives, and the map $(n,e)\mapsto F(n,e)$, is $C^1_b$ from $[0,+\infty)^d\times\mathbb{R}^\dimE$ to
  $\mathbb{R}^\dimE$. In addition, we assume that, for all $i,j\in\{1,2,\ldots,d\}$ and $e\in\mathbb{R}^\dimE$
  \begin{equation}
  \label{eq:hyp-1}
      b_i(e)-d_i(e)>0\quad \text{and}\quad c_{i}(e,j)\geq \underline{c}
  \end{equation}
  for some constant $\underline{c}>0$.
\end{hyp}

The assumption~\eqref{eq:hyp-1} means that no trait has a tendency to go fast to extinction when the population size is small (i.e.
when the nonlinear term in the death rate is small), regardless of the value of the environmental variable $e$ and of the traits that
are non-extinct in the population. In particular, as we will see in Section~\ref{lem_control_reaching_time_of_0}, with probability
converging to 1, the total population size $N^K_1(t)+\ldots+N^K_d(t)$ is of the order of $K$ for times much longer than $\log K$.

\medskip

Next, we consider an assumption which prevents the long-term coexistence of more than two species (this is the so-called ``invasion implies fixation'' assumption~\cite{geritz2002invasion,geritz2005resident}). In this assumption, the term ``globally asymptotically stable'' for a steady state $x^\star$ of a dynamical system $\dot{x}=g(x)$ taking values in $[0,+\infty)^\ell$ will mean that any solution of this dynamical system with initial condition in $(0,+\infty)^\ell$ converges to $x^\star$ for large times. This is a particular case of the classical notion of global asymptotic stability~\cite{guckenheimer-holmes-2013}.

To state this assumption, we introduce for any $i\in\{1,2,\ldots,d\}$ and $e\in\mathbb{R}^\dimE$ the notation 
\[
n^\star_i(e)=\frac{b_i(e)-d_i(e)}{c_i(e,i)}>0
\]
for the unique globally asymptotically stable equilibrium of the dynamical system (with $e$ fixed)
\[
  \dot{x}_i=x_i\big(b_i(e)-d_i(e)-c_i(e,i)x_i\big),
\]
and we also denote by
\[
f_{i,j}(e):=b_i(e)-d_i(e)-c_i(e,j)n_j^\star(e)
\]
the invasion fitness of species $i$ when species $j$ is at equilibrium and all other species are negligible. The notion of fitness is classical in population genetics~\cite{metz1992}; in our case, it represents the growth rate of individuals in species $i$ when all species are extinct except $j$ which is assumed to be at its equilibrium density $n^\star_j(e)$.
  


\begin{hyp}
  \label{hyp:Lotka-Volterra-max-2}
  For any $i\neq j$ in $\{1,2,\ldots,d\}$ and any $e\in\mathbb{R}^\dimE$ such that $f_{j,i}(e)>0$ we assume that $f_{i,j}(e)<0$. This is
  the so-called \emph{invasion implies fixation (IIF)} condition for the 2-dimensional dynamical system (with $e$ fixed)
  \[
    \begin{cases}
      \dot{x}_i=x_i(b_i(e)-d_i(e)-c_i(e,i)x_i-c_i(e,j)x_j) \\
      \dot{x}_j=x_j(b_j(e)-d_j(e)-c_j(e,i)x_i-c_j(e,j)x_j),    
    \end{cases}
  \]
  and it is well-known that it is equivalent to the fact that when the equilibrium $(n^\star_i(e),0)$ is hyperbolic and linearly
  unstable, then the equilibrium $(0,n^\star_j(e))$ is hyperbolic and globally asymptotically stable (see e.g. \cite{geritz1998}).
\end{hyp}

\begin{rem}
  \label{rem:intra-inter}
Notice that when inter-specific competition is ``stronger'' than intra-specific competition, in the sense that
\[
\forall\,e\in\R^\dimE,\quad c_{i}(e,j)c_{j}(e,i)\geq c_{i}(e,i)c_{j}(e,j),\qquad\forall\, i,j=1,\ldots,d
\]
condition (IIF) is direct. Indeed, if 
\[
f_{j,i}(e) = b_j(e)-d_j(e)-c_j(e,i)\,\frac{b_i(e)-d_i(e)}{c_i(e,i)} >0, 
\]
then
\[
\frac{b_j(e)-d_j(e)}{b_i(e)-d_i(e)}>\frac{c_j(e,i)}{c_i(e,i)}\geq \frac{c_j(e,j)}{c_i(e,j)},
\]
and rearranging terms we find that
\[
0>b_i(e)-d_i(e) -c_i(e,j)\,\frac{b_j(e)-d_j(e)}{c_j(e,j)} = f_{i,j}(e).
\]
\end{rem}

We will also need the following very mild assumption on the model parameters.

\begin{hyp}
    \label{hyp:not-many-zeros}
    We assume that none of the functions $f_{i,j}$ for $i\neq j$ in $\{1,2,\ldots,d\}$ are identically zero on a non-empty open ball in $\R^\dimE$.
\end{hyp}

\medskip

Our final assumption requires that there is no heteroclinic cycle~\cite{may1975} in the dynamical system 
\[
\dot{x}_i=x_i\left(b_i(e)-d_i(e)-\sum_{j=1}^d c_i(e,j)x_j\right),\quad 1\leq i\leq d,
\]
regardless of the value of the environmental variable $e\in\R^\dimE$.

\begin{hyp}
  \label{hyp:no-heteroclinic-loop}
  Assume for all $2\leq n\leq d$, all distinct $i_1,i_2,\ldots,i_n\in\{1,2,\ldots,d\}$ and all $e\in\R^\dimE$, the following $2n$ inequalities are not simultaneously all satisfied: for all $\ell\in\{1,\ldots,n\}$, with the convention $i_{n+1}=i_1$,
  \begin{equation*}
  f_{i_{\ell+1},i_\ell}(e)\geq 0\text{ and }f_{i_{\ell},i_{\ell+1}}(e)\leq 0.
  \end{equation*}
\end{hyp}

This assumption is technical and somewhat restrictive. However, we will see in Section~\ref{sec:main-2} that our main result still
holds true if Assumption~\ref{hyp:no-heteroclinic-loop} is not satisfied, possibly only up to a finite time. Note also that, in the
case where $d=2$ and Assumption~\ref{hyp:Lotka-Volterra-max-2} is satisfied, this assumption is equivalent to assuming that $f_{1,2}$
and $f_{2,1}$ have no common zeros. Note that the property that two functions from $\R^\dimE$ to $\R$ have no common zeroes is generic if
$\dimE=1$ but not if $\dimE\geq 2$, hence Assumption~\ref{hyp:no-heteroclinic-loop} is somewhat restrictive if $d\geq 2$ and $\dimE\geq 2$. However,
by Remark~\ref{rem:intra-inter}, this property is automatic if interspecific competition is strictly stronger than intraspecific
competition. 
From the biological point of view, this assumption means first that, in a fixed environment, two species cannot coexist on the long term and second that once a species is able to invade some dominant species, then it will eventually replace it. The first assumption is somewhat restrictive biologically because it implies competitive exclusion of all species except one, but the second one holds true in general for most model parameters and environment values~\cite{geritz2005}. However, we will see in our examples that a weaker form of coexistence is allowed in our model due to the change of environment. 

\subsection{Main result}
\label{sec:main}

We introduce the so-called exponents of the species sizes: for all $i\in\{1,2,\ldots,d\}$ and all $t\geq 0$, we set
\begin{equation}
\label{eq:exponents}
  \beta^K_i(t)=\frac{\log(1+N^K_i(t\log K))}{\log K},\quad\text{i.e.}\quad N^K_i(t\log K)=K^{\beta^K_i(t)}-1
\end{equation}
and $\beta^K(t)=(\beta^K_1(t),\ldots,\beta^K_d(t))$. Note that the time scaling of $1/\log K$ that we introduced in the definition of
$\beta^K$ is the same as the one introduced in~\eqref{eq:environment}. Note also that $N^K_i(t\log K)=0$ if and only if $\beta^K_i(t)=0$ and that if $\beta^K_i(t)$ converges as $K\to+\infty$ to some value in $(0,1)$, then $N^K_i(t\log K)$ is asymptotically negligible with respect to $K$, but is non-zero, meaning that species $i$ is non extinct at time $t\log K$ but negligible with respect to a total population size of order $K$.

Before stating our result, some topological reminder is in order. Recall that if $X$ is a topological space, a set $A \subset X$ is
said to be nowhere dense if its closure has empty interior. A Baire meager set is the (at most) countable union of nowhere dense
sets. Also recall the Baire Category Theorem, which states that on a complete metric space, every Baire meager set has empty interior
(see  \cite{brezis2011}, Thm. 2.1). The notion of Baire meager sets is a topological notion of negligible sets, different from Lebesgue's
negligible sets, which is well-suited to prove generic properties of dynamical systems~\cite{demazure1999}.

Our main result is the following.

\begin{thm} 
  \label{thm:BDP-max-2-a.e.i.c} 
  Assume that Assumptions~\ref{hyp:regularity-Ltka-Volterra},~\ref{hyp:Lotka-Volterra-max-2},~\ref{hyp:not-many-zeros} and~\ref{hyp:no-heteroclinic-loop} are satisfied. Then, for all $\ell\in\{1,2,\ldots,d\}$, there exists $D\subset(0,1)^{d-1}\times\R^\dimE$, a Baire meager set, with the following properties: whenever  
  \[\left(\beta_1^{(0)},\ldots,\beta^{(0)}_{\ell-1},\beta^{(0)}_{\ell+1},\ldots,\beta_d^{(0)},e_0\right)\in (0,1)^{d-1}\times\R^\dimE\setminus D,
  \]
  if
  \begin{itemize}
  \item $(E^K(0))_{K \ge 1}$ converges in probability to $e_0$;
  \item $(\beta^K(0))_{K \ge 1}$ converges in probability to $\beta^{(0)}= (\beta_1^{(0)},\ldots,\beta^{(0)}_{\ell-1},1, \beta^{(0)}_{\ell+1},\ldots,\beta_d^{(0)})$;
  \item $(\frac{1}{K}N^K_{\ell}(0))_{K \ge 1}$ converges in probability to $n^\star_{\ell}(e_0)$, 
  \end{itemize}
  then:
  \[
    (\beta^K(t), E^K(t \log (K)))_{t\geq 0}\xrightarrow[K\to+\infty]{}(\beta(t),e(t))_{t\geq 0},
  \]
  in probability in $L^\infty_\text{loc}([0,+\infty),\mathbb{R}^{d+\dimE})$, where the limit $(\beta(t),e(t))_{t\geq 0}$ is deterministic
  and defined as explained next. There exists an increasing sequence $(t_k)_{k\geq 0}$ in $[0,+\infty]$ converging to $+\infty$ and a
  sequence $(\ell^\star_k)_{k\geq 0}$ in $\{1,2,\ldots,d\}$ such that $(\beta(t),e(t))_{t\geq 0}$ is inductively defined on the intervals
  $[t_k,t_{k+1})$ as:
  \begin{itemize}
  \item For $k=0$, let $t_0=0$, $\ell^\star_0=\ell$, $\beta(0)=\beta^{(0)}$ and $e(0)=e_0$.
  \item {For $k\geq 0$ assume that $t_0,\ldots,t_k$ and $\ell^\star_0,\ldots,\ell^\star_k$ are constructed, $t_k<+\infty$ and $(\beta(t),e(t))=((\beta_1(t),\dots,\beta_d(t)),e(t))$ is constructed on $[0,t_k]$.
    Then for $t\geq t_k$, we set}
    \begin{equation}
      \begin{cases}
        \beta_{\ell^\star_k}(t)=1, & \\
        \beta_i(t)=0, & \forall i\text{ such that }\beta_i(t_k)=0, \\
        \dot{\beta}_i(t)=f_{i,\ell^\star_k}(e(t)), & \forall i\neq
        \ell^\star_k\text{ such that }\beta_i(t_k)>0, \\
        \dot{e}(t)=F_{\ell^\star_k}(e(t)), &
      \end{cases}
      \label{eq:exponent_system}
    \end{equation}
    where for all $i\in\{1\ldots,d\}$ and $e\in\R^\dimE$,
    \begin{equation}
      \label{eq:def-F_i}
      F_i(e)=F(n^\star_i(e)\delta_i,e),
    \end{equation}
    until the first time $t>t_k$ where 
    \begin{enumerate}
    \item either $\beta_i(t)=1$ for some $i\neq \ell^\star_k$,
    \item or $\beta_i(t)=0$ for some $i\neq \ell^\star_k$ such that $\beta_i(t_k)>0$.
    \end{enumerate}
    This defines the time $t_{k+1}$. {If $t_{k+1}=\infty$ then we stop the induction. If $t_{k+1}<\infty$, then only one of the previous two events occurs at time $t_{k+1}$.} If the first point holds true at time $t_{k+1}$ for some index $i$, then this index is unique
    and we set $\ell^\star_{k+1}=i$  and proceed with the next induction step. If the second point holds true at time $t_{k+1}$ for
    some index $i$, we set $\ell^\star_{k+1}=\ell^\star_k$ and proceed with the next induction step. {In both cases,  $(\beta(t),e(t))$ has a left limit as $t\nearrow t_{k+1}$, denoted $(\beta(t_{k+1}),e(t_{k+1}))$ leading to the construction of the process on $[0,t_{k+1}]$.} 
  \end{itemize}
\end{thm}

This result characterizes the asymptotic behavior of the exponents of all species and proves that there is always, except for one Baire meager set of initial conditions, only a single species with size of the order of $K$, that we will call in the sequel the \emph{dominant species}, except at times of change of dominant species, which are all close to the times $(t_k)_{k\geq 1}$ defined in the last result, at which there is an instantaneous change of slope in the limit dynamics of $\beta(t)$. In addition, due to Assumption~\ref{hyp:Lotka-Volterra-max-2}, we will actually prove that at each time $t_k$ where an exponent hits $1$, then necessarily the exponent of the former dominant species becomes strictly smaller than 1 immediately after $t_k$, which indeed means that there is no coexistence of several non-negligible species in the limit dynamics of Theorem~\ref{thm:BDP-max-2-a.e.i.c}. We will also prove that, when species $i$ is dominant at time $t$, with probability converging to 1, its density $\frac{1}{K}N^K_i(t\log K)$ is close to $n^\star_i(e(t))$, which means that the dominant species density is always at equilibrium. The times $(t_k)_{k\geq 1}$ may be also times when one exponent $\beta_i(t)$ hits zero, and we will also prove below that this corresponds to an actual extinction of species $i$, in the sense that, if $\beta_i(t_k)=0$ and $\beta_i(t)>0$ for $t<t_k$, then, for all $s>t_k$, $N^K_i(s\log K)=0$ with probability converging to 1 (see Theorem~\ref{thm:BDP-max-2} below).

In the last result, we assume that every species initially has some positive exponent. This is because we do not take into account the possibility of mutations, so without loss of generality, we can remove the initially extinct species to apply the theorem. In a model with mutations, some exponents $\beta_i^K(0)$ might be 0, meaning that the species are initially extinct, but can appear afterwards by mutations (see Remark~\ref{rem:extension} below for a possible extension). More generally, the models and assumptions under which we prove Theorem~\ref{thm:BDP-max-2-a.e.i.c} are somewhat restrictive, but we will discuss in the next subsection how the model can be generalized and the result extended to more general settings.

\subsection{Strategy of proof and possible extensions}
\label{sec:main-2}

The proof of Theorem~\ref{thm:BDP-max-2-a.e.i.c} is based on the fact that, when $E^K(0)$ and $\beta^K(0)$ converge in probability to deterministic $e(0)\in\R^\dimE$ and $\beta(0)\in(0,1]^d$ respectively, with exactly one $\ell\in\{1,2,\ldots,d\}$ such that $\beta_{\ell}(0)=1$, then we can prove that   
\[
(\beta^K(t), E^K(t \log (K)))_{t\in [0,t_1)}\xrightarrow[K\to+\infty]{}(\beta(t),e(t))_{t\in [0,t_1)},
\]
in probability in $L^\infty_\text{loc}([0,t_1),\mathbb{R}^{d+\dimE})$, where the limit $(\beta(t),e(t))_{t\geq 0}$ is the deterministic solution of 
\[
  \begin{cases}
    \beta_{\ell}(t)=1, & \forall t\in[0,t_1),\\
    \dot{\beta}_i(t)=f_{i,\ell}(e(t)), & \forall i\neq
    \ell,\ \forall t\in[0,t_1) \\
    \dot{e}(t)=F_\ell(e(t)), & \forall t\in[0,t_1),
  \end{cases}
\]
with 
\[
    t_1 := \inf\{t>0, \, \beta_i(t)\in \{0,1\} \text{ for some $i\neq \ell$}\}.
\]
If such $t_1$ is well defined, then we would like to repeat the argument inductively on some interval $[t_1,t_2)$. However, to do so
we need to be able to define some new initial condition $\beta_1\in[0,1]^d$ after time $t_1$ with exactly one $\ell'$ such that
$(\beta_1)_{\ell'}=1$. This cannot be done directly, because we have defined $t_1$ as the instant for which either $\beta_i(t_1)=0$
or $\beta_i(t_1)=1$ for some $i\neq\ell$. We can distinguish several possibilities:
\begin{enumerate}
\item In the case where $\beta_i(t_1)<1$ for all $i\neq\ell$, there will be no problem in extending the argument by just considering that species $i$ such that $\beta_i(t_1)=0$ become extinct (i.e. we set $\beta_i(t)=0$ for all $t\geq t_1$) and we can apply the first step above by considering a smaller number of initially non-extinct species.
\item There is only one $i\neq \ell$ such that $\beta_i(t_1)=1$ and $\beta_j(t_1)>0$ for all $j\neq \ell,i$. Since the exponent $\beta_i$ hits 1 at time $t_1$ and $\beta_i(t)<1$ for all $t<t_1$, then necessarily its derivative at time $t_1$, given by $f_{i,\ell}(e(t_1))$, is nonnegative.
\begin{enumerate}
    \item If $f_{i,\ell}(e(t_1))>0$, then the invasion implies fixation assumption~\ref{hyp:Lotka-Volterra-max-2} yields
      $f_{\ell,i}(e(t_1))<0$, and hence after a small time $\eta>0$, we will have $\beta_\ell(t_1+\eta)<1$, hence, taking $\eta>0$
      small enough, we can repeat the first step changing $\ell$ with $i$ after time $t_1+\eta$ and with initial non-dominant
      exponents all positive and smaller than 1.
    \item If $f_{i,\ell}(e(t_1))=0$, then it is unclear if species $i$ will indeed be able to invade and replace the former dominant species $\ell$, or if it will remain negligible (its fitness $f_{i\ell}(e(t))$ could become negative after time $t_1$). This is a case where we cannot determine easily what the limit dynamics after time $t_1$ is.
\end{enumerate}
\item There are at least two different indices $i,j$ both different to $\ell$ such that $\beta_i(t_1)=\beta_j(t_1)=1$. In this case,
  we cannot distinguish which species will become the dominant one after $t_1$, or whether the three species $i$, $j$ and $\ell$ will
  coexist. Indeed the 3-dimensional competitive Lotka-Volterra dynamical system is much richer in terms of long time behavior and its
  behavior may depend on the initial condition~\cite{zeeman1993}. In this case, we cannot define the limit population state just
  after time $t_1$ and we cannot continue the inductive construction of the limit.
\item There is only one $i\neq\ell$ such that $\beta_i(t_1)=1$ and there is some $j\neq i,\ell$ such that $\beta_j(t_1)=0$. Note that this implies that $f_{j,\ell}(e(t_1))\leq 0$. However, after the change of dominant species from $\ell$ to $i$, species $j$ has fitness $f_{j,i}(e(t_1))$ which may be positive. In this case, it is not clear if species $j$ has sufficient time to get extinct while species $\ell$ is dominant or if it can survive up to the time where species $i$ becomes dominant and start increasing. Again, in this case, we cannot determine the limit population state just after time $t_1$ and we need to stop the inductive construction of the limit process.
\end{enumerate}

Therefore, our strategy of proof of Theorem~\ref{thm:BDP-max-2-a.e.i.c} is based on an inductive argument, which may stop because of the occurrence of bad events (corresponding to cases 2.(b), 3.\ and 4.\ above). This method of proof gives a general convergence result valid only under Assumptions~\ref{hyp:regularity-Ltka-Volterra} and~\ref{hyp:Lotka-Volterra-max-2}, stated in Theorem~\ref{thm:BDP-max-2} below. Under Assumption~\ref{hyp:not-many-zeros}, we are next able to prove in Theorem~\ref{thm:Baire-d>2} below that no such bad event occurs in the inductive construction of the limit process for all initial conditions except in a Baire meager set. This allows to define the sequence $(t_k)_{k\geq 1}$ of Theorem~\ref{thm:BDP-max-2-a.e.i.c} up to the first index $k$ such that $t_k=+\infty$, or for all $k\geq 1$ with $t_k<+\infty$. This last case may lead to a finite accumulation of times $t_k$ at $S=\sup_k t_k$. We will finally show in Proposition~\ref{prop:S-infty} below that, under Assumption~\ref{hyp:no-heteroclinic-loop}, we must have $S=+\infty$. Hence Theorem~\ref{thm:BDP-max-2-a.e.i.c} is a direct corollary of the next three results.


%

\begin{thm} 
  \label{thm:BDP-max-2} 
  Assume that Assumptions~\ref{hyp:regularity-Ltka-Volterra} and~\ref{hyp:Lotka-Volterra-max-2} are satisfied, that, as $K\to+\infty$, $E^K(0)$ and $\beta^K(0)$ converge in probability to deterministic $e_0$ and $\beta^{(0)}\in(0,1]^d$, respectively, with exactly one index $\ell\in\{1,2,\ldots,d\}$ such that $\beta^{(0)}_{\ell}=1$. We also assume that $\frac{1}{K}N^K_{\ell}(0)$ converges in probability to $n^\star_{\ell}(e_0)$. Then, there exists a deterministic time $S$ such that, in probability in $L^\infty_\text{loc}([0,S),\mathbb{R}^{d+\dimE})$,
  \begin{equation}
      \label{eq:convergence-main}
          (\beta^K(t), E^K(t \log (K)))_{t\geq 0}\xrightarrow[K\to+\infty]{}(\beta(t),e(t))_{t\geq 0},
  \end{equation}
  where there exists a non-decreasing sequence $(t_k)_{k\geq 0}$ and a sequence $(\ell^\star_k)_{k\geq 0}$ in $\{1,2,\ldots,d\}$
  such that $(\beta(t),e(t))_{t\geq 0}$ is inductively defined on the intervals $[t_k,t_{k+1})$ as follows:
  \begin{itemize}
  \item For $k=0$, let $t_0=0$, $\ell^\star_0=\ell$, $\beta(0)=\beta^{(0)}$ and $e(0)=e_0$.
    \item {For $k\geq 0$ assume that $t_0,\ldots,t_k$ and $\ell^\star_0,\ldots,\ell^\star_k$ are constructed, $t_k<+\infty$ and $(\beta(t),e(t))=((\beta_1(t),\dots,\beta_d(t)),e(t))$ is constructed on $[0,t_k]$.
    Then for $t\geq t_k$, we set}
    \[
      \begin{cases}
        \beta_{\ell^\star_k}(t)=1, & \\
        \beta_i(t)=0, & \forall i\text{ such that }\beta_i(t_k)=0, \\
        \dot{\beta}_i(t)=f_{i,\ell^\star_k}(e(t)), 
        & \forall i\neq
        \ell^\star_k\text{ such that }\beta_i(t_k)>0, \\
        \dot{e}(t)=F_{\ell^\star_k}(e(t)), &
      \end{cases}
    \]
    until the first time $t>t_k$ where 
    \begin{enumerate}
    \item either $\beta_i(t)=1$ for some $i\neq \ell^\star_k$,
    \item or $\beta_i(t)=0$ for some $i\neq \ell^\star_k$ such that $\beta_i(t_k)>0$.
    \end{enumerate}
    This defines the time $t_{k+1}$. If $t_{k+1}=+\infty$, we set $S=+\infty$ and stop the induction. Otherwise, if several indices $i$ satisfy either point 1 or point 2 above at the same time $t_{k+1}$, we set $S=t_{k+1}$ and stop the induction. If only the first point holds true at time $t_{k+1}$ and only for a single index $i$, if $f_{i,\ell^\star_k}(e(t_{k+1}))=0$ then we set $S=t_{k+1}$ and stop the induction; otherwise we must have $f_{i,\ell^\star_k}(e(t_{k+1}))>0$ and then we set $\ell^\star_{k+1}=i$ and continue the induction. If only the second point holds true at time $t_{k+1}$ and only for a single index $i$, if $f_{i,\ell^\star_k}(e(t_{k+1}))=0$
    we set $S=t_{k+1}$ and stop the induction; otherwise we must have $f_{i,\ell^\star_k}(e(t_{k+1}))<0$
    and then we set $\ell^\star_{k+1}=\ell^\star_k$ and continue the induction.
    {In both cases when we continue the induction, $t_{k+1}<+\infty$, then $(\beta(t),e(t))$ has a left limit as $t\nearrow t_{k+1}$, denoted $(\beta(t_{k+1}),e(t_{k+1}))$ leading to the construction of the process on $[0,t_{k+1}]$.}
  \item If the induction does not stop, we set $S=\sup_k t_k$.
  \end{itemize}
  In addition, for all $i\in\{1,\ldots,d\}$ and all $t\in(0,S)$ such that $\beta_i(s)=0$ for some $s<t$,
  \[
  \lim_{K\to+\infty}\mathbb{P}\left(N^K_i(t\log K)=0\right)=1.
  \]
\end{thm}




\begin{dfn}
    In the sequel, we will say that \emph{induction never stops in Theorem~\ref{thm:BDP-max-2}} (or more simply \emph{induction never stops}) if, using the notations of Theorem~\ref{thm:BDP-max-2}, either $S=+\infty$ or $t_k$ is well-defined for all $k\geq 1$ and $S=\sup_k t_k$.
\end{dfn}

\noindent Note that the last condition in the previous definition does not necessarily imply that $S=+\infty$.

\begin{thm}
    \label{thm:Baire-d>2}
    Assume that Assumptions~\ref{hyp:regularity-Ltka-Volterra},~\ref{hyp:Lotka-Volterra-max-2} and~\ref{hyp:not-many-zeros} are satisfied.
    Then there exists a subset $D_0$ of 
    \[
    D:=\bigcup_{i=1}^d (0,1)^{i-1}\times\{1\}\times(0,1)^{d-i}\times\mathbb{R}^\dimE
    \]
    such that for all $i\in\{1,2,\ldots,d\}$
    \[
    \left\{((\beta_\ell)_{\ell\neq i},e)\text{ s.t.\ }(\beta_1,\ldots,\beta_{i-1},1,\beta_{i+1},\ldots,\beta_d,e)\in D_0\right\}
    \]
    is a Baire meager subsets of $(0,1)^{d-1}\times\mathbb{R}^\dimE$ and, for all initial conditions $(\beta(0),e(0))\in D\setminus D_0$,  induction never stops in Theorem~\ref{thm:BDP-max-2}, which means that either $S=+\infty$ or the full sequence $(t_k)_{k\geq 0}$ is well-defined, increasing and $S=\sup_{k\geq 1}t_k$.
\end{thm}

\begin{prop}
\label{prop:S-infty}
Assume that Assumptions~\ref{hyp:regularity-Ltka-Volterra},~\ref{hyp:Lotka-Volterra-max-2} and~\ref{hyp:no-heteroclinic-loop} are satisfied. Then, given an initial condition $(\beta(0),e(0))$ such that induction never stops in Theorem~\ref{thm:BDP-max-2}, we
have $S=+\infty$. 
\end{prop}

The method of proof of Theorem~\ref{thm:BDP-max-2} relies on estimates of exit from an attracting domain for a logistic birth and
death processes coupled with a slow environmental dynamics, converging to the logistic ODE $\dot{x}=x(r(e)-c(e)x)$ with $r(e)>0$ and
$c(e)>0$ when the environmental variable is constant equal to $e$, where the attracting domain is a neighborhood of the equilibrium
$r(e)/c(e)$ (Lemma~\ref{lem:time-to-equi}). This implies that, as long as species $i$ is dominant, the environmental dynamics on the
slow timescale is close to the solution of the dynamical system $\dot{e}=F_i(e)$, where $F_i$ is defined in~\eqref{eq:def-F_i}. This
part of the proof is new compared to existing related results (notably~\cite{champagnat2021stochastic}) because of the coupling with
a slowly varying environment. Our method makes explicit use of the bound $\|F\|_\infty/\log K$ on the speed of evolution of the environmental
variable in~\eqref{eq:environment}. Next we use comparison arguments between non-dominant populations and branching processes in a
slowly varying deterministic environment, for which we prove the convergence of the exponents, which is also a new result
(Lemma~\ref{lem:time-inhomo-LBD}). Finally, a precise study of the fast phase of change of dominant species is required, based on
Lemma~\ref{lem:competition} in Section~\ref{sec:competitive.Lotka.volterra}.

The proof of Theorem~\ref{thm:Baire-d>2}, given in Sections~\ref{sec:Baire-d=2} and~\ref{sec:Baire-d>2}, is based on a precise characterization of the initial conditions leading to a stop of
induction at time $t_k$ in Theorem~\ref{thm:BDP-max-2} for all $k\geq 1$, from which we deduce in particular that these sets are the
image by explicit homeomorphisms of nowhere dense sets, hence are nowhere dense, hence their union is Baire meager. Note that Theorem~\ref{thm:Baire-d>2} is the first result of this kind for the parameter scaling we consider here,
showing that for generic model parameters and all initial conditions in a Baire co-meager
set, the limit dynamics is defined for any positive time. The proof of
Proposition~\ref{prop:S-infty}, given in Section~\ref{sec:S-infty}, is based on the analysis of the signs of fitness functions when
time approaches $S$ assuming by contradiction that $S<+\infty$, using that the environmental variable $e(t)$ converges to some finite
limit when $t\to S$. 

In view of the previous discussion on the method of proof, several extensions of our main results can be considered, gathered in the following remark.

\begin{rem}
  \label{rem:extension}
  \begin{itemize}
  \item In Theorems~\ref{thm:BDP-max-2-a.e.i.c} and~\ref{thm:BDP-max-2}, we assume that initially $\frac{1}{K}N^K_\ell(0)$ converges to $n^\star_\ell(e_0)$. Actually, our result would still hold true if there exists $\varepsilon>0$ such that $\P(N^K_\ell(0)\geq\varepsilon K)\to 1$ when $K\to+\infty$, but the convergence would hold in $L^\infty_\text{loc}((0,+\infty),\mathbb{R}^{d+\dimE})$ in Theorem~\ref{thm:BDP-max-2-a.e.i.c}, or in $L^\infty_\text{loc}((0,S),\mathbb{R}^{d+\dimE})$ in Theorem~\ref{thm:BDP-max-2} (cf.\ e.g.~\cite{champagnat2021stochastic}).
  \item Other density-dependencies in the functions $b_i(\nu,e)$ and $d_i(\nu,e)$ may be considered in place of the Lotka-Volterra model. Everything would work the same in Theorem~\ref{thm:BDP-max-2} provided that there is a single globally asymptotically stable positive equilibrium $n^\star_i(e)$ for the one-dimensional macroscopic ODE of species $i$ alone for all $e\in\R^\dimE$ and that the conditions $f_{j,i}(e)>0$ and $f_{i,j}(e)<0$ imply for the two-dimensional macroscopic system of ODEs for species $i$ and $j$ that the equilibrium $(0,n^\star_j(e))$ is globally asymptotically stable.
  \item We could consider mutations with rate of the order of $K^{-\alpha}$ as in~\cite{champagnat2021stochastic,coquille2021stochastic}. We could also consider horizontal transfer as in~\cite{champagnat2021stochastic}. Theorem~\ref{thm:BDP-max-2} would remain true with the appropriate modifications provided we could extend Lemma~\ref{lem:time-inhomo-LBD} to branching processes with time-inhomogeneous immigration.
  \item Instead of Assumption~\ref{hyp:Lotka-Volterra-max-2}, we could assume that for any $I\subset\{1,2,\ldots,d\}$ such that $I\neq\emptyset$ and any $e\in\mathbb{R}^\dimE$, the dynamical system
  \[
    \forall i\in I,\quad \dot{x}_i=x_i\left(b_i(e)-d_i(e)-\sum_{j\in I}c_i(e,j)x_j\right)
  \]
  admits a unique globally asymptotically stable equilibrium, denoted $n^\star_I(e)=(n^\star_i(I,e))_{i\in I}$, assumed non-identically zero. In this case, a dominant population could be composed of several coexisting species in a set $I\subset\{1,2,\ldots,d\}$, and when a new species $j\not\in I$ invades this dominant population (i.e.\ its exponents hits 1) at some time $t$, then the next set of surviving species would be given by the nonzero coordinates of $n^\star_{I\cup\{\j\}}(e(t))$ 
  (see~\cite{champagnat2011,champagnat2014,coquille2021stochastic}).
  \item Instead of a piecewise deterministic dynamics for the environmental variable, we could consider stochastic dynamics evolving on the same timescale, e.g.\ a diffusion process of the form 
  \[
  dE^K(t)=\frac{1}{\log K}\lambda(\nu^K_t,E^K(t))dt+\frac{1}{\sqrt{\log K}}\sigma(\nu^K_t,E^K(t))dB_t
  \]
  for some Brownian motion $(B_t)_{t\geq 0}$. Theorem~\ref{thm:BDP-max-2} could be extended to this case provided we can prove
  appropriate estimates about the time of exit of the dominant population density from a neighborhood of its equilibrium on a long
  timescale, as described above.
  \item We could also assume that parts of the environmental variable evolve on the fast time-scale of population dynamics. In this case, we would need to replace the equilibrium $n^\star_i(e)$ when species $i$ is dominant and the slow environmental variable has state $e$, with an equilibrium $(n^\star_i(e),u^\star_i(e))$ for the joint limit system of ODE for species $i$ density and fast environmental variable. The method of proof could be extended if we have good convergence properties of this system of ODE, typically that there is only one globally asymptotically stable equilibrium and that the convergence is controlled by some Lyapunov function as e.g.\ in~\cite{champagnat2014}. It would also be possible to assume periodic dynamics of the fast environmental variable as in~\cite{esser2025}. Similarly, we could consider a metapopulation model with fast migration between patches, provided similar good convergence and Lyapunov properties of the limit dynamical system would hold true. Metapopulations with slow migration (typically with rate of the order of $K^{-\gamma}$, as for mutations above) could also be considered: in this case, the exponent of each species within each patch should be considered separately.
  \end{itemize}
\end{rem}

\section{Examples} 
\label{sec:examples}

\subsection{A model of indirect mutualism with periodic selection or extinction}
\label{sec:ex-twodim}

We start with a simple model to exemplify the conditions under which a cyclic substitution of dominance emerges. Let $d=2$ and $\dimE=1$, i.e.\ $e(t)\in\R$ and consider the \textit{per capita} birth and death rates
\[
b_1(e) = e^+,\qquad b_2(e) = (1-e)^+,\qquad d_1(e)=d_2(e)=0,
\]
where $x^+=x\vee 0$. We assume that the intra- and inter-specific competition functions are all constant such that
\[
c_1(e,2) = c_2(e,1)=c,\qquad c_1(e,1)=c_2(e,2)=1
\]
with $c\in[1,+\infty)$. We also assume that the environmental dynamic is driven by the field 
\[
F(n_1,n_2,e) = \lambda(-n_1+n_2)
\]
for some parameter $\lambda>0$ allowing to tune the pace of environmental change compared to the dynamics of the exponents $\beta$.
We finally assume that $E^K(0)$ converges in probability to some $e(0)\in(0,1)$ and $(\beta^K_1(0),\beta^K_2(0))$ converges in probability  to $(1,\beta_2(0))$ with $\beta_2(0)\in(0,1)$, i.e. initially population 1 is dominant and population 2 is non-extinct but negligible.

This choice of parameters means that species 1 is favored when the environment variable is larger than $1/2$ and species 2 is favored
when the environment variable is smaller than $1/2$. Interspecific competition is larger than intraspecific competition, which
ensures that invasion implies fixation (see the Remark~\ref{rem:intra-inter} after Assumption~\ref{hyp:Lotka-Volterra-max-2}). The
function $F$ implies that species 1 consumes the environmental variable, whereas species 2 produces it. This is a simple model of the
situation where two specialist species are consuming two resources $S_1$ and $S_2$ that are interacting negatively and
$e(t)=1/2+S_1(t)-S_2(t)$. This is a classical situation of indirect mutualism through competing
resources~\cite{macarthur1972,levine1976} which is a toy model for several real life examples of niche construction, such as for
granivorous species in desert ecosystems where ants and rodents eat seeds of different
sizes~\cite{davidson1984}. 

Our choice of the parameters implies that, for all $e\in[0,1]$,
\begin{gather*}
n_1^\star(e) =e,\quad n_2^\star(e)=1-e,\quad f_{1,2}(e) =(1+c)e-c=1-(1+c)(1-e), \\ f_{2,1}(e)=1-(1+c)e,\quad F_1(e)=-\lambda e,\quad F_2(e)=\lambda(1-e).
\end{gather*}
In particular $f_{2,1}(e)$ is positive if $e<1/(1+c)$ and negative if $e>1/(1+c)$, and $f_{1,2}(e)$ is positive if $e>c/(1+c)$ and negative if $e<c/(1+c)$. Since $c\geq 1$, this indeed corresponds to the situation where invasion implies fixation and the two fitness functions have common zeros only if $c=1$.

This setting does not exactly correspond to the assumptions of Theorem~\ref{thm:BDP-max-2}, because the assumption $b_i(e)>d_i(e)$ is
not satisfied for all $e$ (we only have $b_i\geq d_i$) and the functions $b_1$, $b_2$ and $F$ are not bounded. However, the limit
environmental dynamics either solves $\dot{e}=-\lambda e$ or $\dot{e}=\lambda(1-e)$ depending on the resident species, hence never
exits $(0,1)$ if $e(0)\in(0,1)$. In particular, for all $T>0$, there exists $\varepsilon>0$ such that
$e(t)\in[\varepsilon,1-\varepsilon]$ for all $t\in[0,T]$. Therefore, we may modify $b_i(e)$ out of the interval
$[\varepsilon,1-\varepsilon]$ to make it fulfill the assumptions of Theorem~\ref{thm:BDP-max-2} without modifying the limit dynamics
on $[0,T]$.


\subsubsection{Illustration of the convergence towards the deterministic limit}

\begin{figure}
  \centering
   \mbox{
    \subfigure[$K=100$.]%
{\includegraphics[width=.32\textwidth, trim = 0cm 0cm 1.4cm 1.cm, clip=true]{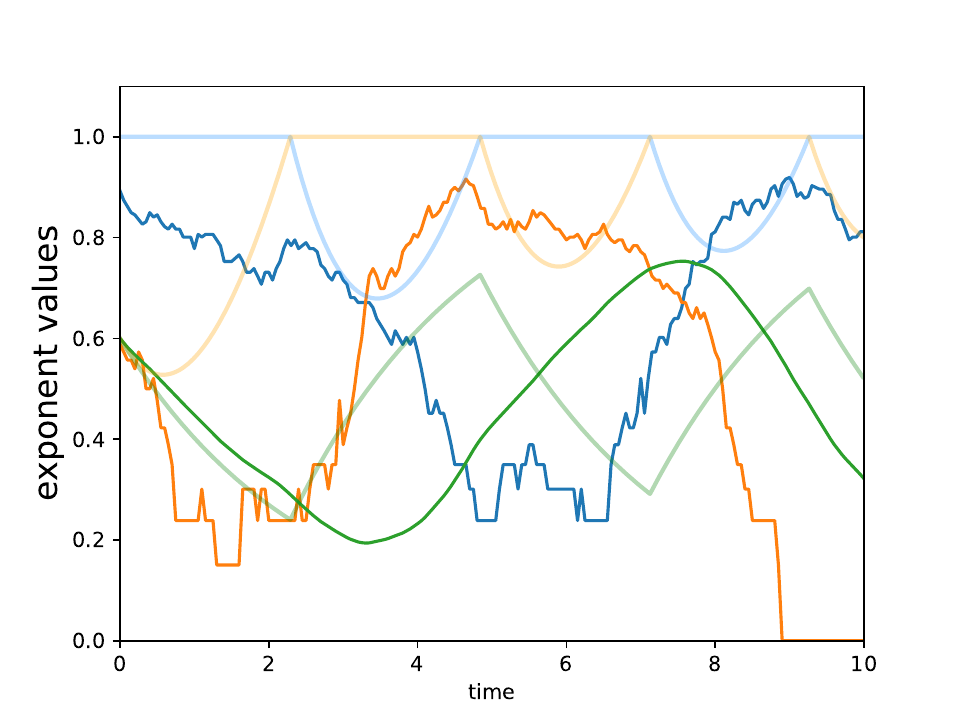}} \quad
    \subfigure[$K=1\,000$.]%
{\includegraphics[width=.32\textwidth, trim = 0cm 0cm 1.4cm 1.cm, clip=true]{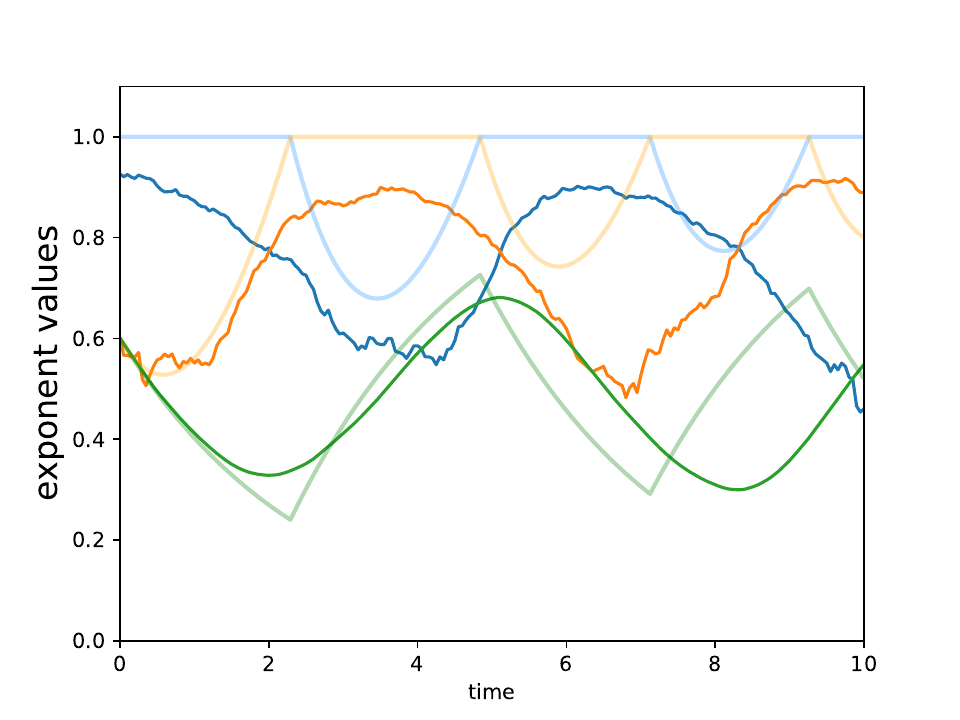}} \quad
  \subfigure[$K=10\,000$.]%
{\includegraphics[width=.32\textwidth, trim = 0cm 0cm 1.4cm 1.0cm, clip=true]{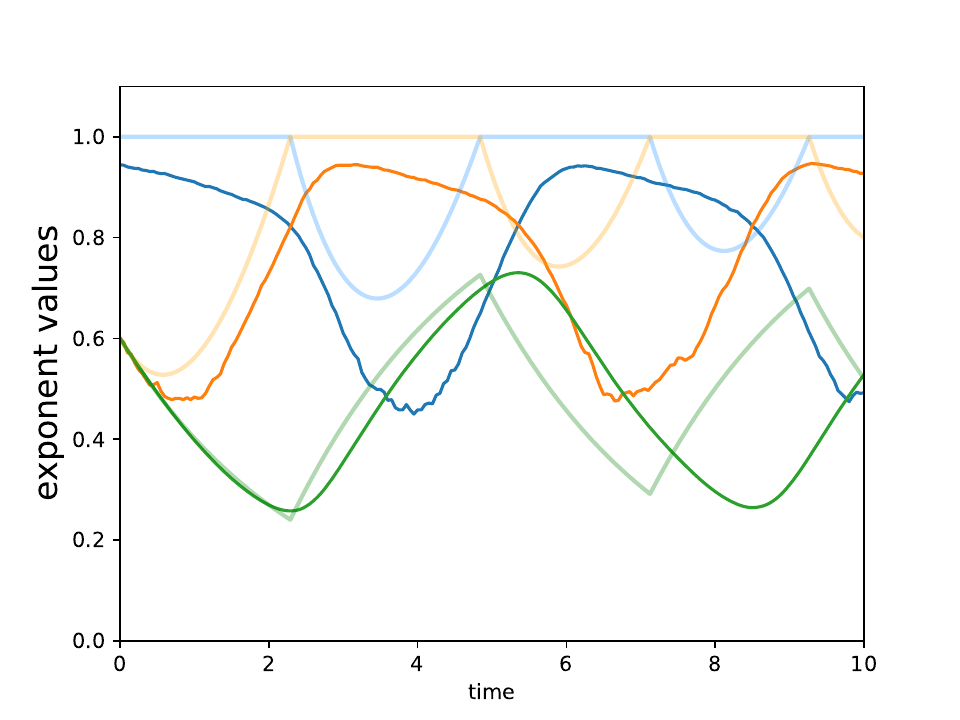}}
    } \\
    \mbox{
    \subfigure[$K=1\,000\,000$.]%
{\includegraphics[width=.32\textwidth, trim = 0cm 0cm 1.4cm 1.cm, clip=true]{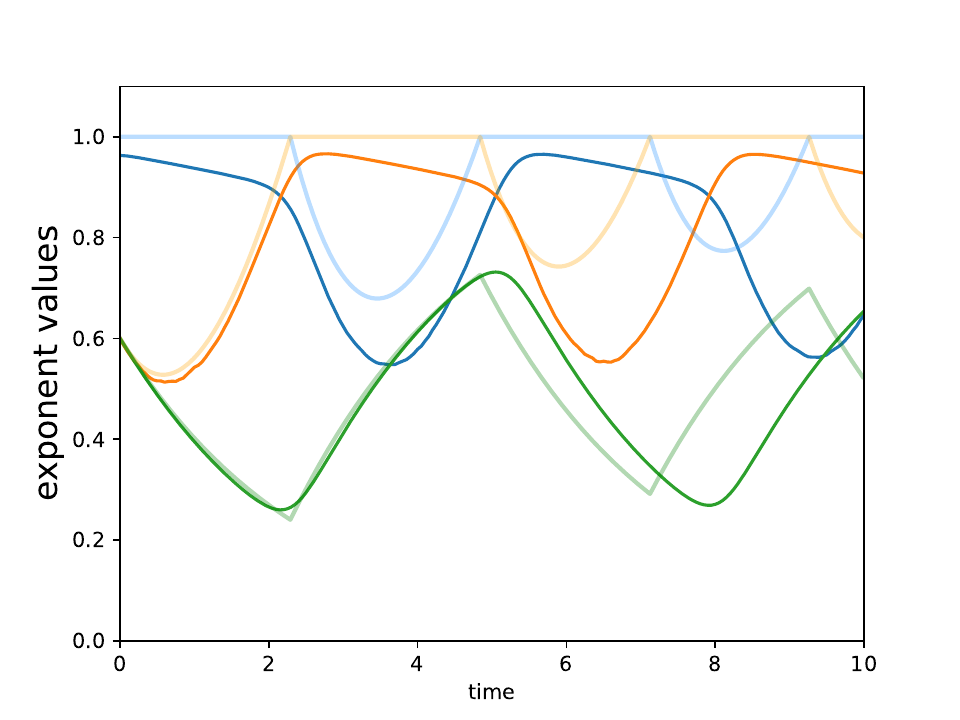}} \quad
  \subfigure[$K=100\,000\,000$.]%
{\includegraphics[width=.32\textwidth, trim = 0cm 0cm 1.4cm 1.cm, clip=true]{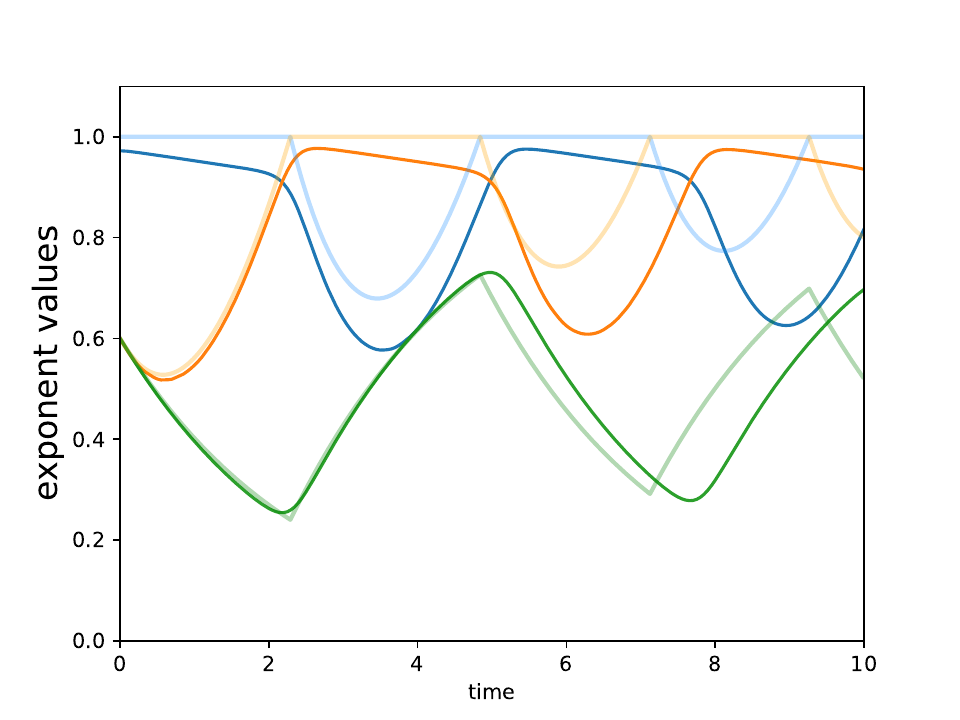}}
    \quad  
    \begin{minipage}{5cm}
    \center
    \vspace{-4cm}
    \includegraphics[width=4cm]{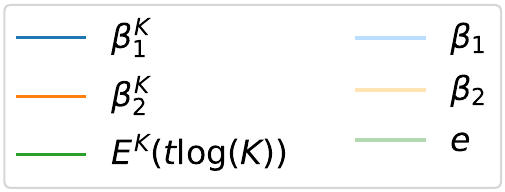}
    \end{minipage}}
\caption{Dynamics of exponents of the first and second species, represented respectively in blue and orange as well as the environmental dynamics, represented in green, for the stochastic process (continuous normal curves) and the asymptotic limit (thick semi-transparent lines) for different values of $K$ and for $c=1.1$, $\lambda=0.4$, $E^K(0)=e(0)=0.6$, $(\beta^K_1(0),\beta^K_2(0))=\left(\frac{\log(1+\lfloor 0.6K \rfloor)}{\log K},\frac{\log(\lfloor K^{0.6} \rfloor)}{\log K}\right)$ and $(\beta_1(0),\beta_2(0))=(1,0.6)$.}
\label{fig.convergence.illustration} 
\end{figure}

Figure~\ref{fig.convergence.illustration} illustrates the convergence \eqref{eq:convergence-main}. The stochastic process $(N^K_1(t),N^K_2(t), E^K(t))_t$ have been simulated by acceptance - rejection Gillespie Algorithm \cite{gillespie1977a} until time $T\log(K)$ to construct the process $(\beta^K_1(t),\beta^K_2(t), E^K(t\log(K)))_t$ on the time interval $[0,T]$, with $T=10$, for different values of $K$. Given our choice of parameters, the limiting process $(\beta_1(t),\beta_2(t), e(t))_t$ does not exhibit population extinction over the time interval $[0, T]$.
In small population size, i.e. with small values of $K$, the behaviour of the process $(\beta^K_1(t),\beta^K_2(t), E^K(t\log(K)))_t$ is stochastic, with potential extinction of the population on the considered time interval. When $K$ increases, we observe a behaviour closer and closer to the deterministic limit process $(\beta_1(t),\beta_2(t), e(t))_t$.  We can notice that the convergence is slow and is not reached, even for $K=100\,000\,000$. However, for large values of $K$, although the stochastic process exhibits a time delay relative to the limit process after a certain time, the limiting process accurately describes the characteristics of the overall long-term behavior of the process in terms of dominant population replacement and resource dynamics.

\subsubsection{Analysis of the deterministic limit}


Let us now study the limit dynamics $(e(t),\beta_1(t),\beta_2(t))$. Before time $t_1$, $e(t)=e(0)\exp(-\lambda t)$ and
\[
\beta_2(t)=\beta_2(0)+t-(1+c)e(0)\frac{1-e^{-\lambda t}}{\lambda}.
\]
Note that the infimum of $t-(1+c)e(0)\left(1-e^{-\lambda t}\right)/\lambda$ is reached at time $t^\star=\ln((1+c)e(0))/\lambda$, which is positive only if $e(0)> 1/(1+c)$. Hence 
\begin{itemize}
    \item if $e(0)> 1/(1+c)$ and $\lambda\beta_2(0)<(1+c)e(0)-\ln((1+c)e(0))-1$, then species 2 goes extinct at time $t_1<t^\star$, and since only species 1 remains, $S=+\infty$ and induction never stops;
    \item if $e(0)> 1/(1+c)$ and $\lambda\beta_2(0)=(1+c)e(0)-\ln((1+c)e(0))-1$, then $t_1=t^\star$, $\beta_2(t_1)=0$ and $f_{2,1}(e(t_1))=0$, so induction stops at time $t_1$;
    \item if $e(0)\leq 1/(1+c)$ or $\lambda\beta_2(0)>(1+c)e(0)-\ln((1+c)e(0))-1$, then $t_1<+\infty$ is the unique root larger than $t^\star$ of $\beta_2(0)-1+t-(1+c)e(0)(1-e^{-\lambda t})/\lambda$, i.e.\ $\beta_2(t_1)=1$, $f_{2,1}(e(t_1))>0$ and (since invasion implies fixation) $f_{1,2}(e(t_1))<0$, hence induction does not stop at time $t_1$.
\end{itemize}

From now on, we assume that $\lambda\beta_2(0)>(1+c)e(0)-\ln((1+c)e(0))-1$ and we denote $u_1:=1-e(t_1)$. Note that $u_1>c/(1+c)>1/2$ since $\beta_2$ was increasing just before time $t_1$. Similarly, we will define for all $k\geq 0$ such that $t_k$ is well-defined,
\[
u_k:=\begin{cases}
    e(t_k) & \text{if }k\text{ is even,} \\
    1-e(t_k) & \text{if }k\text{ is odd.}
\end{cases}
\]

Let us now study the limit dynamics on $[t_1,t_2]$: for all $s\in [0,t_2-t_1]$, $1-e(t_1+s)=u_1 \exp(-\lambda s)$ and
\[
\beta_1(t_1+s)=1+s-(1+c)u_1\frac{1-e^{-\lambda s}}{\lambda}.
\]
Similarly as above, $\beta_1$ reaches its minimum in $t_1+\ln((1+c)u_1)/\lambda>0$, then
\begin{itemize}
    \item if $(1+c)u_1-\ln((1+c)u_1)>1+\lambda$, then species 1 goes extinct at time $t_2<t_1+\ln((1+c)u_1)/\lambda$, $S=+\infty$ and induction never stops;
    \item if $(1+c)u_1-\ln((1+c)u_1)=1+\lambda$, then $\beta_1(t_2)=0$ and $f_{1,2}(e(t_2))=0$, so induction stops at time $t_2:=t_1+\ln((1+c)u_1)/\lambda$;
    \item if $(1+c)u_1-\ln((1+c)u_1)<1+\lambda$, then $t_2<+\infty$ and $t_2-t_1$ is the unique positive root of $t-(1+c)u_1(1-e^{-\lambda t})/\lambda$, i.e.\ $\beta_1(t_2)=1$, $f_{1,2}(e(t_2))>0$ and $f_{2,1}(e(t_2))<0$, hence induction does not stop at time $t_2$.
\end{itemize}

From now on, we assume that $\lambda\beta_2(0)>(1+c)e(0)-\ln((1+c)e(0))-1$ and $(1+c)u_1-\ln((1+c)u_1)<1+\lambda$. We have seen above that $t_2$ satisfies
\[
t_2-t_1=(1+c)u_1\frac{1-e^{-\lambda(t_2-t_1)}}{\lambda},
\]
which implies that
\[
	u_2 = e(t_2) = 1-u_1e^{-\lambda(t_2-t_1)}=1+\frac{\lambda(t_2-t_1)}{1+c}-u_1.
\]
Now, $1-u_2=u_1e^{-\lambda(t_2-t_1)}$, so $\lambda(t_2-t_1)=\ln \frac{u_1}{1-u_2}$ and thus
\[
	u_2 = 1-u_1 + \frac{\ln u_1-\ln(1-u_2)}{1+c}.
\]
Hence, given $u_1$, $u_2$ does not depend on $\lambda$. Note however that $u_1$ does depend on $\lambda$.

Let us define
\[
  \psi(x) = (1+c)x + \ln (1-x)
\]
and observe that $\psi$ is increasing on $[0,c/(1+c)]$ and decreasing on $[c/(1+c),1[$. Since $f_{1,2}(u_2)=f_{1,2}(e(t_2))>0$, this means that $u_2$ is the unique solution to $\psi(u_2)=\psi(1-u_1)$ such that $u_2\in[c/(1+c),1]$.


Hence, proceeding by induction, it is easy to check as above that, as long as $t_k$ is well-defined for some $k\geq 1$ with $\beta_1(t_k)=\beta_2(t_k)=1$ and induction did not stop up to time $t_k$, then $u_{k+1}$ is the unique solution in $[c/(1+c),1]$ of
\begin{equation}
    \label{eq:induction-u_k}
    \psi(u_{k+1})=\psi(1-u_k)
\end{equation}
and
\begin{equation}
\label{eq:t_k_ex1}
    t_{k+1}-t_k=\frac{1+c}{\lambda}(u_k+u_{k+1}-1).
\end{equation}
Let us define $u_k$ and $t_k$ for all $k\geq 1$ by the previous formulas, even when induction has stopped before $t_k$, and, by abuse of notation, let us denote by $(t_k)_{k\geq 1}$ the obtained sequence of times, although it might differ from the sequence defined in Theorem~\ref{thm:BDP-max-2}. Having introduced these notations, we can now state our main result on this example. 

\begin{thm}
\label{thm:ex1}
\begin{enumerate}[label=(\roman*)]
    \item \label{ex.uniq.sol} There exists a unique $u^\star$ solution in $[c/(1+c),1]$ of $\psi(u^*)=\psi(1-u^*)$.
        \item \label{ex.Baire} For the previous model, the Baire meager set $B$ of bad initial conditions in Theorem~\ref{thm:BDP-max-2} is the set of $(e(0),\beta_2(0))\in(0,1)^2$ such that either
        \begin{enumerate}
            \item \label{beta2_value_stop_t1} $(e(0),\beta_2(0))\in A$ with 
            \[A:= \left\{(e,\beta) \quad \text{s.t.} \quad e> 1/(1+c) \quad \text{and} \quad \lambda\beta=(1+c)-\ln((1+c)e-1\right\},\]
            \item 
            $(e(0),\beta_2(0))\notin A$ and $(1+c)u_1-\ln((1+c)u_1)=1+\lambda$,
            \item or 
            $(e(0),\beta_2(0))\notin A$, $(1+c)u_1-\ln((1+c)u_1)<1+\lambda$, $u_1<u^\star$ and there exists $k\geq 1$ such that $(1+c)u_k-\ln((1+c)u_k)=1+\lambda$,
        \end{enumerate}
        where $u_1$ is the unique solution in $[c/(1+c),1]$ of $\psi(u_1)=\psi(1-e(0))-\lambda(1-\beta_2(0))$ and, given $u_k\in[c/(1+c),1]$, $u_{k+1}$ is the unique solution in $[c/(1+c),1]$ of $\psi(u_{k+1})=\psi(1-u_k)$.
        \item \label{ex.extinction} For initial conditions $(e(0),\beta_2(0))\in(0,1)^2\setminus B$, there is extinction of one of the two species, meaning that there exists $k\geq 1$ such that $\beta_1(t_k)=0$ or $\beta_2(t_k)=0$, in the three following cases:
        \begin{enumerate}
            \item \label{ex.extinction.1} 
            $(e(0),\beta_2(0))\in E$ with 
            \[E:= \left\{(e,\beta) \quad \text{s.t.} \quad e> 1/(1+c) \quad \text{and} \quad \lambda\beta<(1+c)-\ln((1+c)e-1\right\},\]
            \item \label{ex.extinction.2} 
            $(e(0),\beta_2(0))\notin E \cup A$, $u_1> u^\star$ and $(1+c)u_1-\ln((1+c)u_1)>1+\lambda$,
            \item \label{ex.extinction.3} or 
            $(e(0),\beta_2(0))\notin E \cup A$, $u_1\leq u^\star$ and $(1+c)u^\star-\ln((1+c)u^\star)> 1+\lambda$.
        \end{enumerate}
        We denote by $C$ the set of such initial conditions $(e(0),\beta_2(0))$.
      \item\label{ex.cycle} For all initial conditions $(e(0),\beta_2(0))\in(0,1)^2\setminus (B\cup C)$, that is under the assumption
        that $\lambda\beta_2(0)>(1+c)e(0)-\ln((1+c)e(0))-1$ and either $u_1\geq u^\star$ and $(1+c)u_1-\ln((1+c)u_1)<1+\lambda$, or
        $u_1<u^\star$ and $(1+c)u^\star-\ln((1+c)u^\star)\leq 1+\lambda$, the times $t_k$ in Theorem~\ref{thm:BDP-max-2} are
        well-defined for all $k\geq 1$ and the sequence $(u_k)_{k\geq 1}$ converges monotonously when $k\to+\infty$ to $u^\star$ and
        $S=+\infty$. If in addition $c>1$ (resp.\ $c=1$), then $t_{k+1}-t_k$ converges to $\frac{1+c}{\lambda}(2u^\star-1)$ (resp.\
        to 0). 
    \end{enumerate}
\end{thm}

In biological terms, in case \ref{ex.cycle} with $c>1$, the niche construction process is quasi-periodic, meaning that the two species
alternate by depleting their own resource and building the niche of the other species by producing the other resource, with a
dynamics converging to a limit cycle. In case \ref{ex.cycle} with $c=1$, the niche construction process also alternates between the two
species but with faster and faster changes of dominant population. In case \ref{ex.extinction}, one of the two species goes extinct at a time of the
order of $T\log K$ for some $T>0$. After this extinction, if the surviving species is the first one, then $e(t)\to 0$ when $t\to\infty$
and hence the density $n_1^\star(e(t))=e(t)$ of species 1 converges to 0, hence species 1 also goes extinct, but on a longer time
scale. We have a similar conclusion if species 2 survives. In conclusion, long term survival of the two species is only possible
through periodic mutual niche construction of the two species, a phenomenon which is often termed indirect
mutualism~\cite{levine1976,davidson1984}. 

\begin{figure}[h]
  \centering
  \mbox{
\subfigure[$\beta_2(0)=0.21$.\label{ex.extinction.before.t1}]%
{\includegraphics[width=.33\textwidth]{
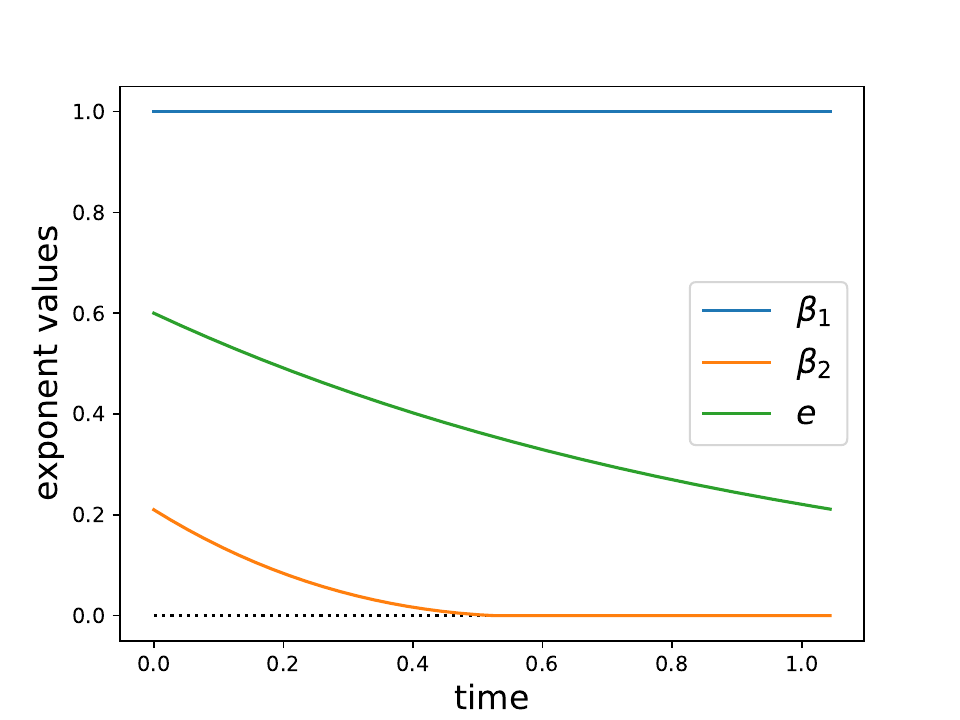}}\quad
\subfigure[$\beta_2(0)=\beta_2^c=0.2122$.\label{ex.stop.at.t1}]%
{\includegraphics[width=.33\textwidth]{
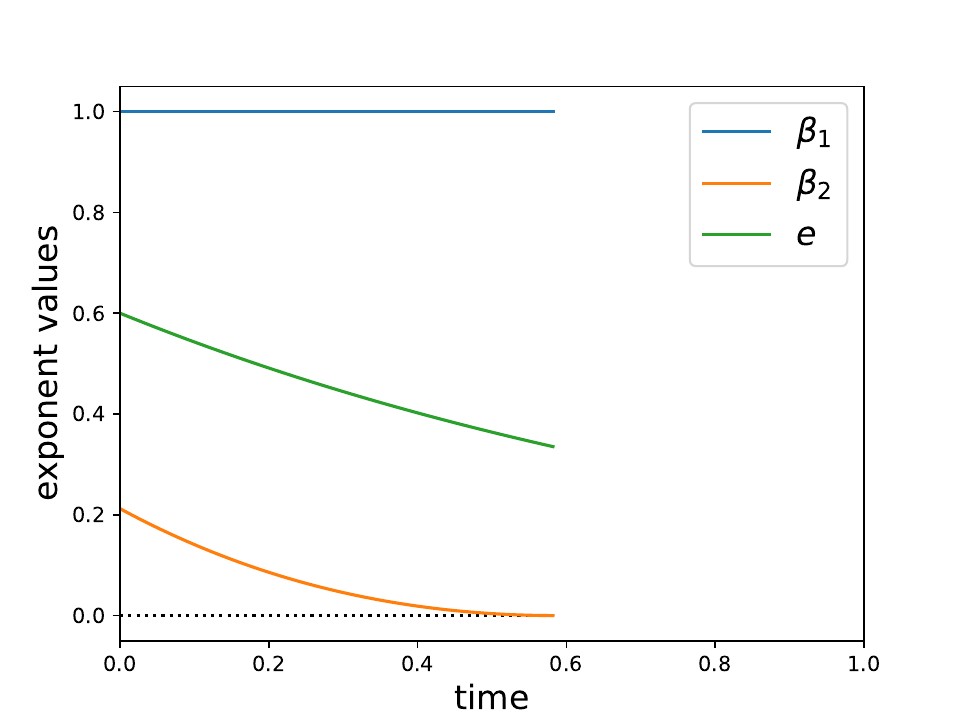}}\quad
    \subfigure[ $\beta_2(0)=0.22$.\label{ex.coexistence.after.t1}]%
{\includegraphics[width=.33\textwidth]{
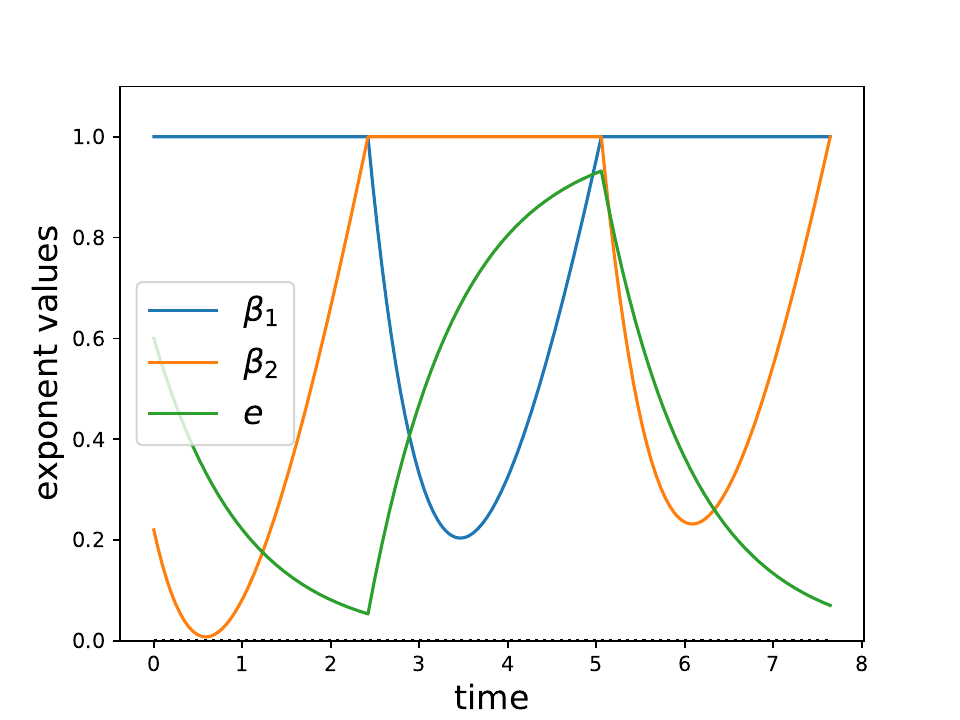}}}
\caption{Dynamics of $(\beta_1,\beta_2)$ solution of~\eqref{eq:exponent_system} for the example of subsection~\ref{sec:ex-twodim} with $c=2$, $\lambda=1$ and $e(0)=0.6$ and different values of $\beta_2(0)$.}
  \label{fig:ex-DA}
\end{figure}
\begin{figure}[h]
  \centering
  \mbox{\subfigure[$c=1.5$, $\lambda=1$, $e(0)=0.3$, $\beta_2(0)=0.1$.\label{ex.c.eq.1}]%
{\includegraphics[width=.47\textwidth]{
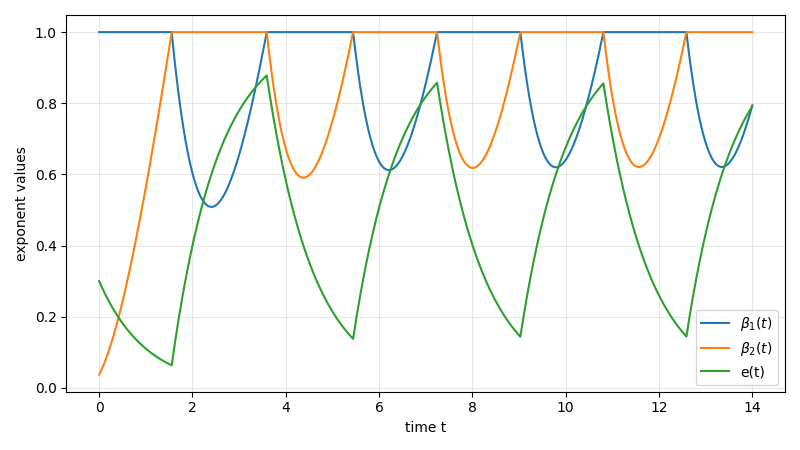}}\quad
    \subfigure[$c=1$, $\lambda=1$, $e(0)=0.3$, $\beta_2(0)=0.1$.\label{ex.c.geq.1}]%
{\includegraphics[width=.47\textwidth]{
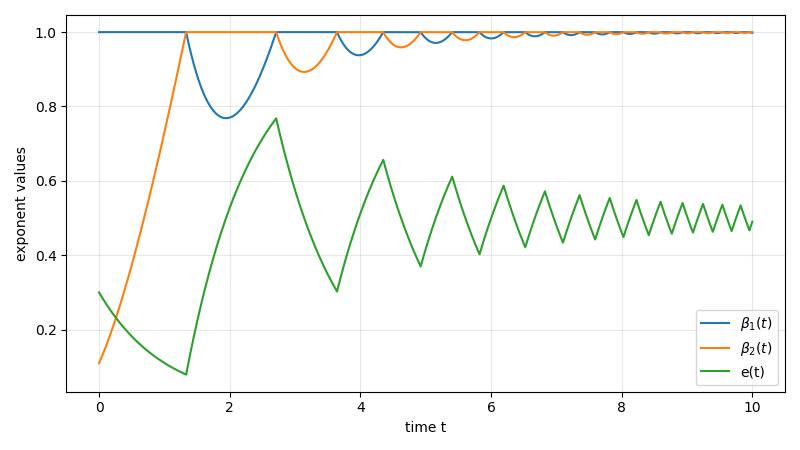}}} \\
\mbox{\subfigure[$c=1.5$, $\lambda=0.3$, $e(0)=0.4$, $\beta_2(0)=0.9$.\label{ex.several.changes}]%
{\includegraphics[width=.47\textwidth]{
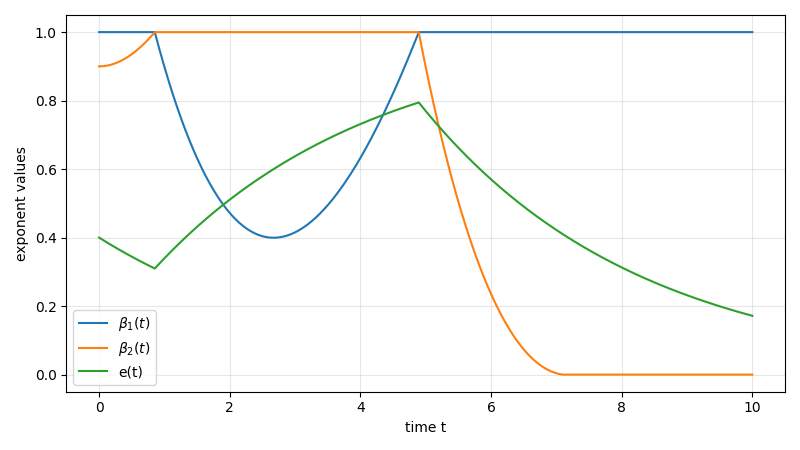}}}
\caption{Numerical illustration of the possible limit dynamics of the exponents $\beta_1$ and $\beta_2$ and the environmental
  variable $e$: we observe (a) 
convergence to a limit cycle; (b) 
accelerating exchanges of dominant populations and convergence of $e(t)$ to $1/2$ when $t\to+\infty$; (c) 
extinction of the population at time $t_3$ rather than time $t_1$ as in Figure~\ref{fig:ex-DA}(c); this requires to slow down the environmental dynamics by taking $\lambda=0.3$ to leave sufficient time to the exponent $\beta_2(t)$ to reach 0 before $e(t)$ becomes smaller than $1/(1+c)$, which is the threshold at which $f_{2,1}$ vanishes. All simulations were implemented in Python 3.12 with a custom explicit Euler scheme and uniform time mesh with $\Delta t = 10^{-3}$.}
\label{fig:ex1-bis}
\end{figure}

Numerical simulations of the limit dynamics $(e(t),\beta_1(t),\beta_2(t))$ are shown in
  Figures~\ref{fig:ex-DA} and~\ref{fig:ex1-bis}. 
  Figure~\ref{fig:ex-DA} shows the different behaviors close to the critical point when induction stops at time $t_1$ (given by Theorem~\ref{thm:ex1}-\ref{beta2_value_stop_t1}). For our choice of parameters ($c=2$, $\lambda=1$ and $e(0)=0.6$), this critical value for $\beta_2(0)$ is $\beta_2^c=0.2122$. Figures~\ref{ex.extinction.before.t1} and \ref{ex.coexistence.after.t1} show the dynamics of $(\beta_1(t),\beta_2(t))$ for a value of $\beta_2(0)$ respectively slightly smaller and slightly larger than $\beta_2^c$. We observe very different
  behaviors, where, respectively, species 2 goes extinct before becoming dominant as stated in
  Theorem~\ref{thm:ex1}-\ref{ex.extinction} in Figure~\ref{ex.extinction.before.t1}, or long term coexistence of the two species occurs, which successive time intervals where each
  species is dominant and convergence to a limit cycle as stated in Theorem~\ref{thm:ex1}-\ref{ex.cycle} in Figure~\ref{ex.coexistence.after.t1}. In Figure~\ref{ex.stop.at.t1}, $\beta_2(0)=\beta_2^c$, then the induction stops at time $t_1$ as stated by Theorem~\ref{thm:ex1}-\ref{beta2_value_stop_t1}. The convergence result of Theorem~\ref{thm:BDP-max-2} does not allow to conclude on the asymptotic behavior of the stochastic process afer the time $t_1$.
  
  Figures~\ref{ex.c.eq.1} and~\ref{ex.c.geq.1}
  show the difference between the asymptotic dynamics when $c>1$ and $c=1$, respectively. We observe a limit cycle in
  Figure~\ref{ex.c.eq.1} where the time length of the successive dominant phases for both species converge to some limit value,
  whereas in Figure~\ref{ex.c.geq.1} the successive dominant phases become shorter and shorter, with an acceleration of the
  changes of dominant species and convergence of $\beta_1(t)$ and $\beta_2(t)$ to 1 and of $e(t)$ to $1/2$, as stated in
  Theorem~\ref{thm:ex1}-\ref{ex.cycle}. Figure~\ref{ex.several.changes} illustrates the fact that, when some species goes extinct, this does not
  necessarily occur at time $t_1$: several changes of dominant species may occur before, as shown in Theorem~\ref{thm:ex1}-\ref{ex.extinction}.

\begin{proof}[Proof of Theorem~\ref{thm:ex1}]
We first prove~\ref{ex.uniq.sol}: we observe that the function $\psi(x)-\psi(1-x)$ is increasing on $[1/2,x^\star(c)]$ and decreasing on $[x^\star(c),1]$, where $x^\star(c)=\frac{1+\sqrt{1-2/(1+c)}}{2}$. Its value at $1/2$ is 0, and it converges to $-\infty$ when $x\to 1$. In addition, it is elementary to check that $x^\star(c)\geq c/(1+c)$ for all $c\geq 1$.  Hence, \ref{ex.uniq.sol} is proved.

Before proving Points~\ref{ex.Baire} to~\ref{ex.cycle}, we first study the sequence $(u_k)_{k\geq 1}$. The previous discussion also shows that $u_{k+1}> u_k$ if $u_k< u^\star$, $u_{k+1}<u_k$ if $u_k> u^\star$, and $u_{k+1}=u^\star$ if $u_k=u^\star$. Since in addition $\hat{\psi}^{-1}(\psi(1-x))$ is increasing on $[c/(1+c),1]$, where $\hat{\psi}$ is the restriction of $\psi$ to $[c/(1+c),1]$, we deduce from~\eqref{eq:induction-u_k} that, if $u_1<u^\star$, then $u_k<u^\star$ for all $k\geq 1$, and similarly if $u_1>u^\star$, then $u_k>u^\star$ for all $k\geq 1$. Combining these two facts, we deduce that $(u_k)_{k\geq 1}$ is strictly monotonic (unless $u_1=u^\star$, in which case it is constant), hence it admits a limit, which must be a fixed point of $\hat\psi^{-1}(\psi(1-x))$ on $[c/(1+c),1]$. By~\ref{ex.uniq.sol}, we deduce that $u_k\to u^\star$ for any initial condition $u_1\geq c/(1+c)$.

Let us also observe that the function $x\mapsto (1+c)x-\ln((1+c)x)$ is increasing on $[c/(c+1),1]$. Hence if $u_1\geq u^\star$ and $(1+c)u_1-\ln((1+c)u_1)<1+\lambda$, then $(1+c)u_k-\ln((1+c)u_k)<1+\lambda$ for all $k\geq 1$, and if $u_1< u^\star$ and $(1+c)u^\star-\ln((1+c)u^\star)\leq 1+\lambda$, then $(1+c)u_k-\ln((1+c)u_k)<1+\lambda$ for all $k\geq 1$. Conversely, if $u_1> u^\star$ and $(1+c)u_1-\ln((1+c)u_1)>1+\lambda$ or $u_1\leq u^\star$ and $(1+c)u^\star-\ln((1+c)u^\star)> 1+\lambda$, then there exists $k\geq 1$ such that $(1+c)u_k-\ln((1+c)u_k)>1+\lambda$.

We can now prove Points~\ref{ex.Baire} to~\ref{ex.cycle}. We start with~\ref{ex.Baire}. 
We already noticed that induction stops at time $t_1$ iff $\lambda\beta_2(0)=(1+c)e(0)-\ln((1+c)e(0))-1$, and never stops iff $\lambda\beta_2(0)<(1+c)e(0)-\ln((1+c)e(0))-1$. If $\lambda\beta_2(0)>(1+c)e(0)-\ln((1+c)e(0))-1$, then induction stops at time $t_2$ if $(1+c)u_1-\ln((1+c)u_1)= 1+\lambda$ and never stop if $(1+c)u_1-\ln((1+c)u_1)> 1+\lambda$. Let consider the remaining case where $(1+c)u_1-\ln((1+c)u_1)< 1+\lambda$. By similar arguments used for the previous cases, induction stops iff there exists $k\geq 1$ such that $(1+c)u_k-\ln((1+c)u_k)=1+\lambda$. Since $(u_k)_{k\geq 1}$ is non-increasing if $u_1\geq u^\star$, this case cannot lead to a stop in induction, so if induction does not stop at times $t_1$ and $t_2$ we can restrict to the case $u_1<u^\star$. Then~\ref{ex.Baire} is clear.

For~\ref{ex.extinction}, we already proved Points~\ref{ex.extinction.1} and~\ref{ex.extinction.2} above. Therefore, we can restrict to the case where $\lambda\beta_2(0)>(1+c)e(0)-\ln((1+c)e(0))-1$ and $u_1\leq u^\star$, in which case $(u_k)_{k\geq 1}$ is non-decreasing. If $u_1=u^\star$, extinction occurs iff $(1+c)u^\star-\ln((1+c)u^\star)> 1+\lambda$. Otherwise, if $u_1<u^\star$, extinction will occur after finitely many induction steps iff $(1+c)u_k-\ln((1+c)u_k)> 1+\lambda$ for some $k\geq 1$. Since this sequence is increasing and converges to $(1+c)u^\star-\ln((1+c)u^\star)$, this will occur iff $(1+c)u^\star-\ln((1+c)u^\star)>1+\lambda$. Hence~\ref{ex.extinction} is proved.

We now come to the proof of~\ref{ex.cycle}. For any initial condition $(e(0),\beta_2(0))\in(0,1)^2\setminus (B\cup C)$, $(1+c)u_k-\ln((1+c)u_k)<1+\lambda$ for all $k\geq 1$ hence induction never stops, and the sequence $(u_k)_{k\geq 1}$ converges monotonously to $u^\star$. In the case where $c>1$, $u^\star\geq c/(1+c)>1/2$ so, by~\eqref{eq:t_k_ex1}, $t_{k+1}-t_k$ converges to $(1+c)(2u^\star-1)/\lambda>0$ and hence $S=\lim_{k\to\infty}t_k=+\infty$. If $c=1$, $u^\star=1/2$ and $t_{k+1}-t_k\to 0$ by~\eqref{eq:t_k_ex1}. It only remains to prove that the sequence $(t_k)_{k\geq 1}$ converges to $+\infty$ when $c=1$. This follows from the next lemma, hence~\ref{ex.cycle} is proved.
\end{proof}


\begin{lem}
In the case $c=1$, we have
\begin{align*}
\lim_{n \to \infty} n(t_n-t_{n-1}) = \dfrac{3}{\lambda}
\end{align*}
where $t_n$ is defined for any $n\geq 1$ by~\eqref{eq:t_k_ex1}. 
\end{lem}

\begin{proof}
  Assume $c=1$, i.e.\ $\psi(x)=2x+\log(1-x)$. Write $e_{n}= e(t_n)$, and consider the function $\varphi(x)= \log(x)-2x$. Equation
  \eqref{eq:induction-u_k} is equivalent to:
\begin{align} \label{eq:phi}
\varphi (e_{2n+1})= \varphi (e_{2n}),  
\end{align}
and
\begin{align}\label{eq:psi}
\psi(e_{2n-1})=\psi(e_{2n}).
\end{align}

Write $q_n=1/2-e_{2n-1}$, $p_n=e_{2n}-1/2$. Perform a Taylor expansion around $1/2$ up to the third order on both sides of equations \eqref{eq:phi} and \eqref{eq:psi}  to get:
\begin{align}\label{eq:key_1}
p_{n}^2-q^2_{n+1}= \dfrac{4}{3}(p_n^3+q_{n+1}^3)+O (p_n^4+q_{n+1}^4),
\end{align}
and
\begin{align}\label{eq:key_2}
q_n^2-p_n^2= \dfrac{4}{3}(p_n^3+q_{n}^3)+O (p_n^4+q_{n}^4).
\end{align}

We make a simple observation. Since $\varphi$ and $\psi$ are concave functions with $1/2$ as unique non-degenerate critical point which is a maximum, we obtain that there exists $\delta >0$ and constants  $0 < c_1 < c_2$ such that for every $x \in [1/2-\delta, 1/2+\delta]$, we have: 
\begin{align*}
c_1(x-1/2)^2 \le \varphi(1/2)-\varphi(x) \le c_2 (x-1/2)^2, 
\end{align*}
as well as:
\begin{align*}
c_1(x-1/2)^2 \le \psi(1/2)-\psi(x) \le c_2 (x-1/2)^2.
\end{align*}
Since $p_n \to 0$ and $q_n \to 0$, this readily implies that $p_n=O(q_n)=O(q_{n+1})$. So, from equations \eqref{eq:key_1} and \eqref{eq:key_2}
 we obtain that:
\begin{align}\label{eq:key_3}
p_{n}^2-q^2_{n+1}= \dfrac{4}{3}(p_n^3+q_{n+1}^3)+O (q_{n}^4),
\end{align}
and
\begin{align}\label{eq:key_4}
q_n^2-p_n^2= \dfrac{4}{3}(p_n^3+q_{n}^3)+O (q_{n}^4).
\end{align}
The last equation yields: 
\begin{align*}
(q_n-p_n) (q_n+p_n)= O(q_n^3)
\end{align*}
and from here we deduce $p_n=q_n+O(q_n^2)$, which in turn implies:
\begin{align*}
(q_n-p_n) (2q_n+O(q_n^2))= \dfrac{8}{3}q_n^3 +O(q_n^4), 
\end{align*}
and then:
\begin{align}\label{eq:key_5}
p_n= q_n-\dfrac{4}{3} q_n^2 + O(q_n^3).
\end{align}
A similar result using equation \eqref{eq:key_3} yields:
\begin{align}\label{eq:key_6}
q_{n+1}=p_n-\dfrac{4}{3}p_n^2 + O(q_n^3). 
\end{align}
Inserting \eqref{eq:key_5} into \eqref{eq:key_6} we obtain:
\begin{align*}
q_{n+1}=q_{n}- \dfrac{8}{3}q_n^2 + o(q_n^2).
\end{align*}
We deduce $\dfrac{1}{q_{n+1}}-\dfrac{1}{q_n} \to \dfrac{8}{3}$. Now, let $\varepsilon \in (0,1)$ be given. There exist $N> 0$ such that for every $n > N$ :
\begin{align*}
(1-\varepsilon) \dfrac{8}{3}(n-N)\le \sum_{k=N}^{n-1} \left(\dfrac{1}{q_{k+1}}- \dfrac{1}{q_{k}} \right) = \dfrac{1}{q_{n}} - \dfrac{1}{q_{N}}\le (1+\varepsilon) \dfrac{8}{3} (n-N), 
\end{align*}
which yields that $\dfrac{3}{8(1+\varepsilon)}\le\liminf nq_n \le \limsup nq_n \le  \dfrac{3}{8(1-\varepsilon)}$. Since $\varepsilon$ is arbitrary, we deduce $\lim nq_n = \dfrac{3}{8}$. Then, from \eqref{eq:key_5}, $\lim n p_n = \dfrac{3}{8}$. Finally, from \eqref{eq:t_k_ex1},
\begin{align*}
\lim n (t_n-t_{n-1})= \dfrac{3}{\lambda}, 
\end{align*}
and this finishes the proof.
\end{proof}

\subsection{An example of immune escape in cancer}


Our goal in this section is to illustrate how our convergence result can provide reduced approximate models of heterogeneous malignant tumor growth in interaction with the immune system of an individual, where the immune system tracks the dominant cancer clone. The interactions between tumor and the immune system are complex~\cite{eftimie2011} and the models from the literature are difficult to analyze. We aim at designing simple approximate deterministic models built from the underlying complex stochastic dynamics, on which immunotherapy effects can be tested numerically.

We model the effect of the immune system on tumor populations by mimicking the behavior of cytotoxic immune cells, such as CD8${}^+$T cells, also called effector cells~\cite{eftimie2011}, and the immune memory effect of~\cite{tanaka2020}. We model the heterogeneity of cancer cells by assuming a finite number of clones initially present, with a single dominant one, each of them being characterized by some trait in an abstract immune state space where the immune system evolves by trying to minimize its distance to the trait of the dominant population. We wish to illustrate by numerical simulations how the parameters describing the rate of recognition and the range of tumor heterogeneity (compared to the range of specificity of effector cells in the immune state space) can influence the outcome of the dynamics and in particular restrict the intratumor heterogeneity and favor clonal dominance, as empirically observed in~\cite{milo2018}.

We assume that the tumor is composed of $d=3$ cell types (often called clones in oncology) interacting with the immune state of the organism modeled as a value $E^K(t)$ in an abstract immune state space---here, taken for simplicity as $\R^2$. We assume that the three tumor types are characterized by traits in the immune state space $x_1=(1,0)$, $x_2=(-1/2,\sqrt{3}/2)$ and $x_3=(-1/2,-\sqrt{3}/2)$ respectively. We consider the birth, death and competition rates defined as:
\begin{gather*}
    b_i(e)=1, \quad d_i(e)=\exp\left(-\frac{\|e-x_i\|_2^2}{2\sigma^2}\right),\\
    c_i(e,j)=2,\quad c_i(e,i)=1
\end{gather*}
for $i,j\in\{1,2,3\}$ with $i\neq j$. The immune system approaches the average immune state of the tumor defined as the barycenter of the $x_i$ weighted according to the population densities, according to the vector field
\[
F(n_1,n_2,n_3,e)=-\lambda\big((n_1+n_2+n_3)e-(x_1n_1+x_2n_2+x_3n_3)\big)=\lambda\big(n_1(x_1-e)+n_2(x_2-e)+n_3(x_3-e)\big). 
\]

\begin{figure}[h]
  \centering
  \mbox{\subfigure[$\lambda=1$, $\sigma=0.5$, $e(0)=(1.3,0.7)$, $\beta_1(0)=1$, $\beta_2(0)=0.2$ and $\beta_3(0)=0.4$.]%
{\includegraphics[width=.7\textwidth]{
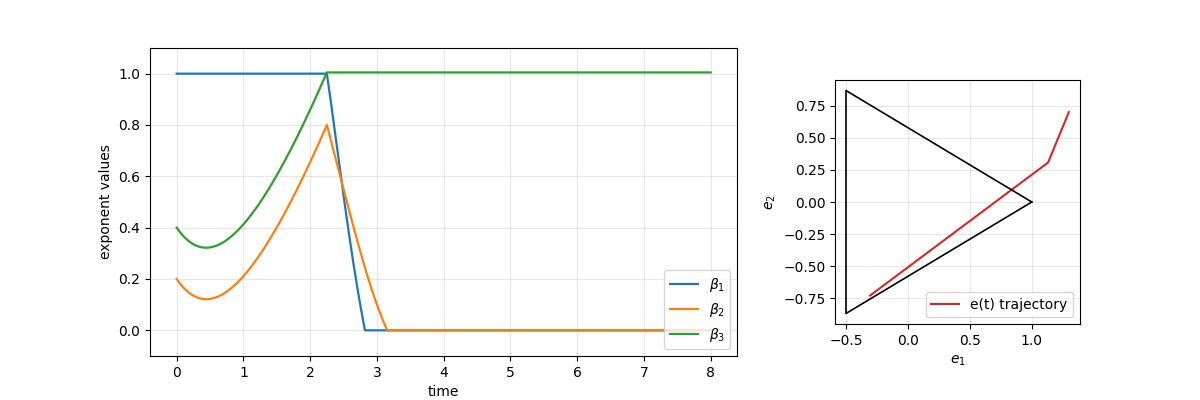}}}\\
    \mbox{\subfigure[$\lambda=1.5$, $\sigma=0.5$, $e(0)=(1.3,0.7)$, $\beta_1(0)=1$, $\beta_2(0)=0.2$ and $\beta_3(0)=0.4$.]%
{\includegraphics[width=.7\textwidth]{
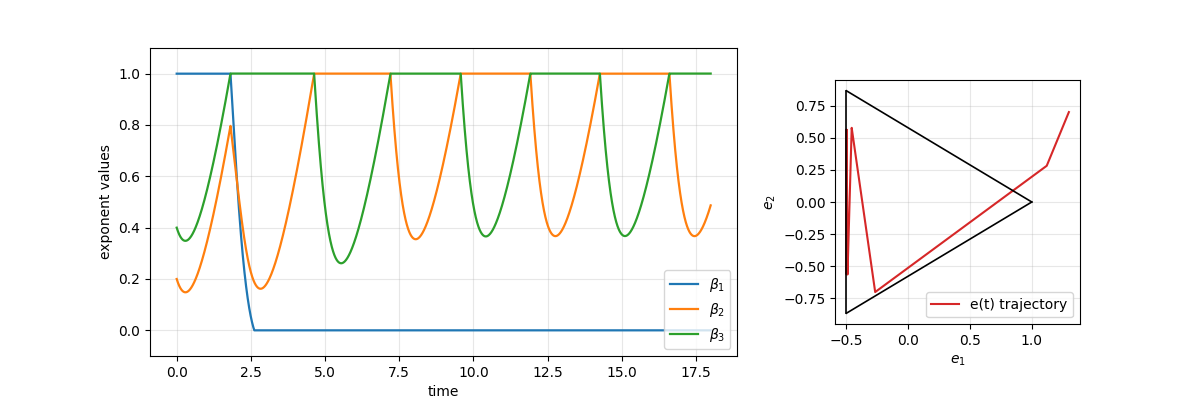}}} \\
\mbox{\subfigure[$\lambda=2$, $\sigma=0.5$, $e(0)=(1.3,0.7)$, $\beta_1(0)=1$, $\beta_2(0)=0.2$ and $\beta_3(0)=0.4$.]%
{\includegraphics[width=.7\textwidth]{
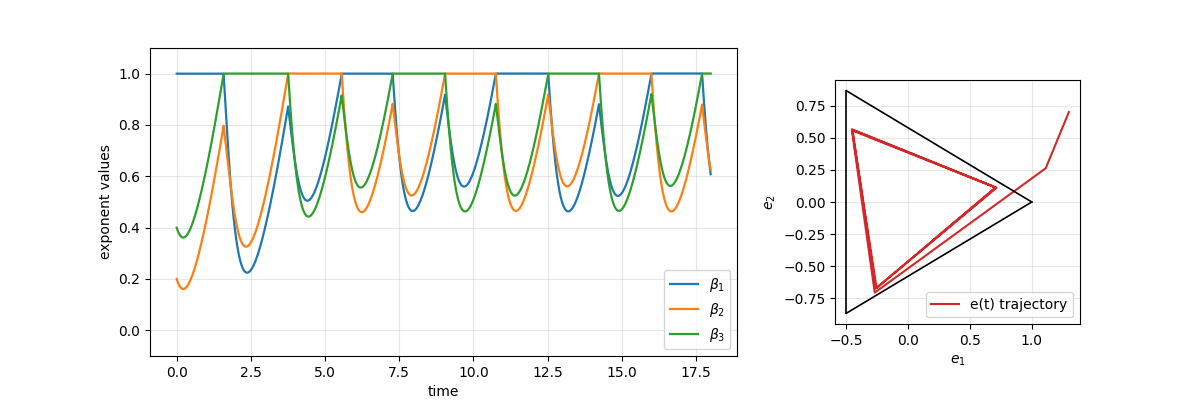}}}
\caption{Influence of the parameter $\lambda$ on the dynamics of the exponents $\beta_1,\beta_2,\beta_3$ defined by \eqref{eq:exponent_system} for the three tumors (left) and the dynamics of the immune system $e=(e_1,e_2)$ (right): (a) small values of $\lambda$ lead to the extinction of tumor types 1 and 2; (b) for intermediate values of $\lambda$, we observe extinction of tumor type 1 and convergence to a limit cycle where tumor types 2 and 3 successively become dominant; (c) large values of $\lambda$ lead to an asymptotic cycle where the 3 tumor types alternatively become dominant.}
\label{fig:ex2-lambda}
\end{figure}

\begin{figure}[h]
  \centering
  \mbox{\subfigure[$\lambda=2$, $\sigma=0.1$, $e(0)=(1.3,0.7)$, $\beta_1(0)=1$, $\beta_2(0)=0.6$ and $\beta_3(0)=0.8$.]%
{\includegraphics[width=.7\textwidth]{
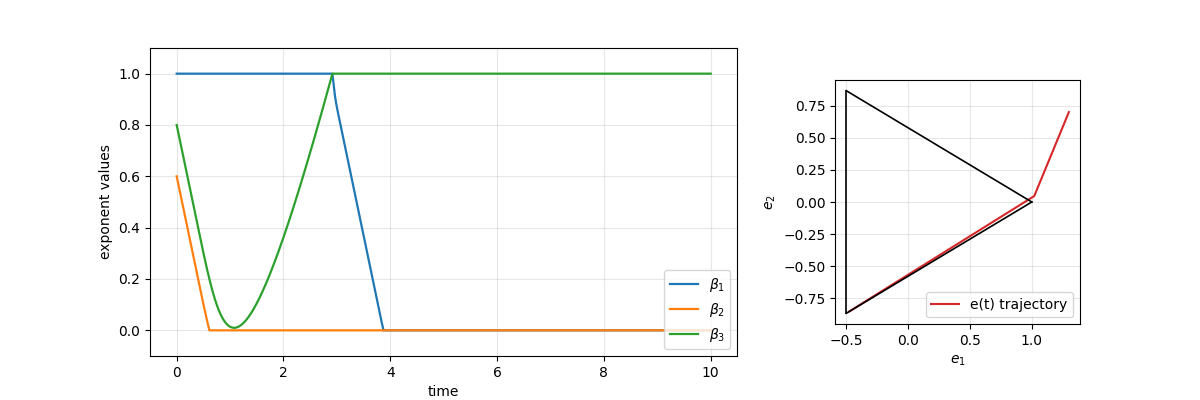}}}\\
    \mbox{\subfigure[$\lambda=2$, $\sigma=0.2$, $e(0)=(1.3,0.7)$, $\beta_1(0)=1$, $\beta_2(0)=0.6$ and $\beta_3(0)=0.8$.]%
{\includegraphics[width=.7\textwidth]{
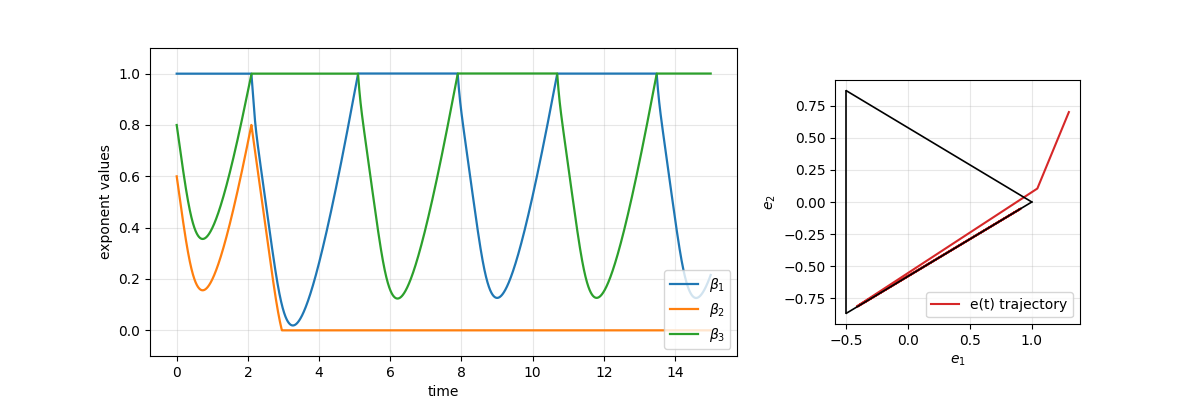}}} \\
\mbox{\subfigure[$\lambda=2$, $\sigma=0.5$, $e(0)=(1.3,0.7)$, $\beta_1(0)=1$, $\beta_2(0)=0.6$ and $\beta_3(0)=0.8$.]%
{\includegraphics[width=.7\textwidth]{
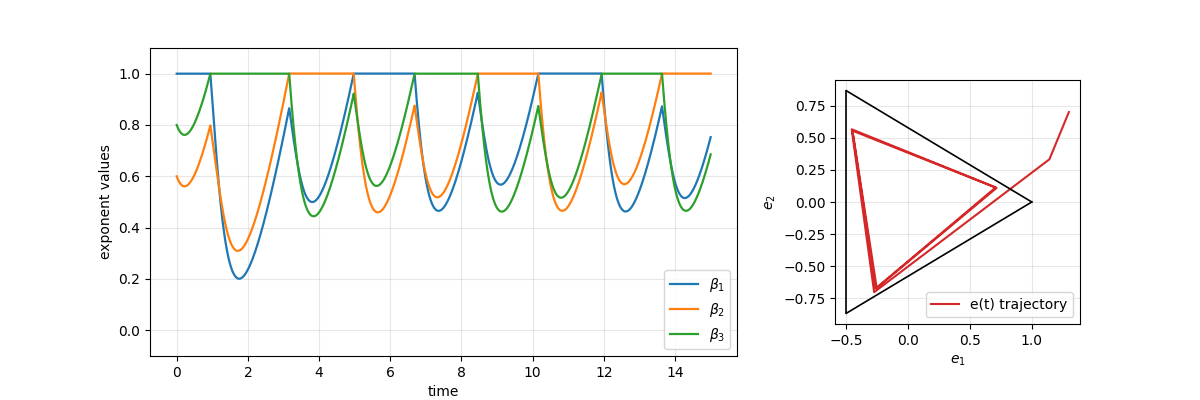}}}
\caption{Influence of the parameter $\sigma$ on the dynamics of the exponents $\beta_1,\beta_2,\beta_3$ defined by \eqref{eq:exponent_system} for the three tumors (left) and the dynamics of the immune system $e=(e_1,e_2)$ (right): (a) small values of $\sigma$ lead to the extinction of tumor types 1 and 2; (b) for intermediate values of $\sigma$, we observe extinction of tumor type 2 and convergence to a limit cycle where tumor types 1 and 3 successively become dominant; (c) large values of $\sigma$ lead to an asymptotic cycle where the 3 tumor types alternatively become dominant.}
\label{fig:ex2-sigma}
\end{figure}

This model mean that the immune system is able to track the barycenter of the three tumor traits $x_1$, $x_2$ and $x_3$ weighted by the relative frequencies of the corresponding tumor types, and the speed of evolution of the immune system is proportional to the total size of the tumor. In particular, $E^K(t)$ does not move if the tumor is extinct. In addition, the parameter $\lambda$ allows to tune the \emph{speed of adaptation of the immune system}. The effect of the immune system on the death rate of a tumor type only depends on its distance to the corresponding tumor trait, in a Gaussian shape with a parameter $\sigma$ which describes \emph{specificity of the immune system}, that is its range of effect in the immune state space. Note also that, by symmetry, there is an equilibrium of the system state with $e=0$ and the 3 species with equal size, although this state cannot be reached by the limit dynamics of Theorem~\ref{thm:BDP-max-2}, as it describes the dynamics as long as there is only one dominant species (except at times of change in dominant species).

With these parameters, we have:
\begin{align*}
f_{i,j}(e)= b_i(e)-d_i(e)-2(b_j(e)-d_j(e)),
\end{align*}
and so the IIF condition of Assumption~\ref{hyp:Lotka-Volterra-max-2} is trivially achieved. 
Moreover, the field $F_i$, associated to the regime where the dominant tumor type is the $i$-th one is given by
\begin{align*}
F_i(e)= -\lambda (b_i(e)-d_i(e)) (e-x_i).
\end{align*}
After the substitution $z_i(t):=e(t)-x_i$, with $e(t)$ the asymptotic dynamics of Theorem~\ref{thm:BDP-max-2}, the associated initial-value problem becomes
\begin{align*}
\dot{z_i} = -\lambda (b_i (e)-d_i(e)) z_i, \quad z_i(0) \in \mathbb{R}^{2}
\end{align*}
and then
\begin{align*}
\dfrac{d}{dt}\Vert z_i \Vert_2^2 = -2\lambda \left( 1-\exp\left(-\frac{\|z_i\|_2^2}{2\sigma^2}\right) \right) \Vert z_i \Vert_2^2,
\end{align*}
so the asymptotic dynamics of the environment driven by the field $F_i$ can be described as the motion of a particle in $\mathbb{R}^2$ which approaches $x_i$ through a straight line with strictly decreasing positive speed, as can be seen in the simulations.

The simulations of Figure~\ref{fig:ex2-lambda} show the influence of the rate of adaptation of the immune system $\lambda$ on the limit dynamics: we observed three different outcomes where either 1, 2 or 3 tumor types survive on the long term. Small values of $\lambda$ lead to the extinction of 2 types as in Figure~\ref{fig:ex2-lambda}(a). In such a situation, the immune state will converge to the corresponding tumor trait, hence the population size of the surviving tumor type converges to 0 when $t$ goes to infinity. Therefore, this case corresponds to a successful control by the immune system of the tumor. Intermediate values of $\lambda$ lead to the extinction of one tumor type and a limit cycle for the two surviving types as in Figure~\ref{fig:ex2-lambda}(b). In this case, the immune system was able to reduce the original heterogeneity of the tumor, but could not eliminate the tumor due to successive outbreaks of different tumor types, to which the immune system need to adapt before eliminating it. During these phases of adaptation, the other surviving tumor type is able to grow again and become dominant. This situation corresponds to the empirical observations of~\cite{milo2018} where the immune system was able to restrict the intratumor heterogeneity and lead to clonal dominance. Finally, for larger values of $\lambda$, the three tumor types are able to survive in an asymptotically cyclic dynamics, as in Figure~\ref{fig:ex2-lambda}(c). Overall, this shows that too large rates of adaptation of the immune system is not good to control the tumor growth because the immune system needs some time to completely eliminate a tumor type before taking care of another emerging tumor type. However, a too small rate of adaptation of the immune system is not good either because it leaves too much time to the tumor to grow and produce adverse effects on the patient. We could incorporate such effects in our model, by allowing for unlimited growth of the tumor as a branching process, in a similar way as in~\cite{champagnat2023,champagnat2026}.

We make similar observations when varying the specificity of the immune system $\sigma$ in Figure~\ref{fig:ex2-sigma}: small values of $\sigma$ lead to the exinction of 2 tumor types as in Figure~\ref{fig:ex2-sigma}(a) and large values of $\sigma$ to the survival of the three tumor types as in Figure~\ref{fig:ex2-sigma}(c). Small specificity of the immune system means that it is unable to affect simultaneously different tumor types, so this finding may seem counter-intuitive. However, another effect is at play in our model: the inter-competition is stronger than the intra-type competition, so small values of $\sigma$ means a longer time before the immune system is able to influence the first dominant tumor type, which leave it sufficient time to eliminate by competition the other tumor types. This explains the extinctions observed in Figures~\ref{fig:ex2-sigma}(a) and~(b).

\begin{figure}
\centering
\includegraphics[width=0.75\textwidth]{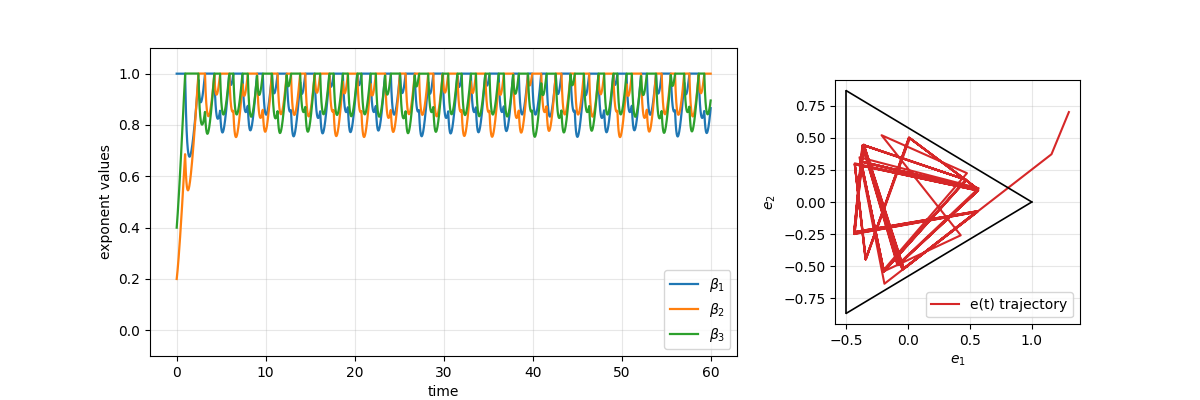}
\caption{Simulation of the limit system for the example of immune escape in cancer for parameters
$\lambda=5$, $\sigma=1$, $e(0)=(1.3,0.7)$, $\beta_1(0)=1$, $\beta_2(0)=0.2$ and $\beta_3(0)=0.4$. The dynamics does not seem to converge to a limit cycle.}
\label{fig:ex2-last}
\end{figure}

Finally, more complex dynamics that those observed in Figures~\ref{fig:ex2-lambda} and~\ref{fig:ex2-sigma} may occur, as can be observed in Figure~\ref{fig:ex2-last}, where the successive dominant tumor types do not seem to follow a cyclic pattern.

\section{Preliminary results}
\label{sec:prelim}

In this section, we gather results about the construction of our niche construction process (in Section~\ref{sec:construction}), the
fact that the total population size is of the order of $K$ (in Section~\ref{sec:total-pop}), the continuity of the exponents
$\beta^K_i$ (in Section~\ref{sec:preliminary}) and technical lemmas needed in the proof of Theorem~\ref{thm:BDP-max-2}. In
Section~\ref{sec:continuity-hitting-0}, we give continuity estimates on the hitting time of 0 by exponents. In
Section~\ref{subsec.LBD}, we study the exponent of one-dimensional time-inhomogeneous linear birth-death processes. In
Section~\ref{sec:logistic}, we prove convergence to equilibrium density of logistic time-inhomogeneous birth death processes.
Finally, in Section~\ref{sec:competitive.Lotka.volterra}, we study a two-dimensional process representing the dynamics of change of
dominant population. 

\medskip

In the rest of the paper, we denote by $|\cdot|$ the standard Euclidean norm in $\R^\dimE$ and $\R^d$, we denote by $\|\cdot\|_\infty$
(respectively $\|\cdot\|_\text{Lip}$) the $L^\infty$ norm (respectively the Lipschitz seminorm) of functions on real vector spaces equipped with the standard Euclidean norm.
By Assumption~\ref{hyp:regularity-Ltka-Volterra}, the functions $b_i$, $d_i$, $c_i(\cdot, j)$, $F$, $n^\star_i$, and $f_{i,j}$ are
bounded and Lipschitz. Note also that since their left and right derivatives are 0 or fitness functions, the functions $\beta_i$ defined in
Theorems~\ref{thm:BDP-max-2-a.e.i.c} or~\ref{thm:BDP-max-2} are also Lipschitz. We denote 
\begin{gather*}
  \bar{b}:=\max_{1\leq i\leq d}\|b_i\|_\infty,\quad \bar{d}:=\max_{1\leq i\leq d}\|d_i\|_\infty, \quad \bar{c}:=\max_{1\leq i,j\leq
    d}\|c_i(\cdot,j)\|_\infty, \quad \bar{f}=\max_{1\leq i\neq j\leq
    d}\|f_{i,j}\|_\infty\\
  \|b\|_\text{Lip} := \max_{1\leq i\leq d} \|b_i\|_\text{Lip}, \quad \|d\|_\text{Lip} := \max_{1\leq i\leq d } \|d_i\|_\text{Lip},
  \quad
  \|c\|_\text{Lip} := \max_{1\leq i,j\leq d} \|c_i(\cdot,j)\|_\text{Lip},\\
  \|n^\star\|_\text{Lip} := \max_{1\leq i\leq d } \|n^\star_i\|_\text{Lip}, \quad \|f\|_\text{Lip} := \max_{1\leq i\neq j\leq d}
  \|f_{i,j}\|_\text{Lip}\quad \text{and}\quad\|\beta\|_\text{Lip}:= \max_{1\leq i\leq d} \|\beta_{i}\|_\text{Lip}.
\end{gather*}

\subsection{Construction of the process}
\label{sec:construction}

As a piecewise deterministic Markov process with locally bounded jump rates and globally Lipschitz vector fields for the piecewise
deterministic dynamics, existence and martingale properties of the process $(N^K(t),E^K(t))_{t\geq 0}$ could be easily deduced from
the general results of~\cite{davis1984piecewise}. However, to contruct couplings between this process and simpler ones, we will need
a more explicit construction, based on classical ideas~\cite{fournier-meleard2004}. We will see examples of such couplings in
Sections~\ref{sec:total-pop} and~\ref{sec:preliminary}.

Consider $2d$ independent Poisson point measures $(P^b_i(ds,dk,dy))_{1\leq i\leq d}$ and $(P^d_i(ds,dk,dy))_{1\leq i\leq d}$ on
$\R_+\times\N\times\R_+$ with intensity measure $\text{Leb}\times\sum_{n\geq 1}\delta_n\times\text{Leb}$, where $\text{Leb}$ denotes
Lebesgue's measure on $\R_+$, on a filtered probability space $(\Omega,({\cal F}_t)_{t\geq 0},\P)$, where the filtration
$({\cal F}_t)_{t\geq 0}$ is right-continuous and such that the random vector $(N^K(0),E^K(0))$ is independent of $P^b_i$ and $P^d_i$
and is ${\cal F}_0$-measurable.

We can then construct the process $(N^K(t),E^K(t))_{t\geq 0}$ as the unique solution of the following stochastic differential equation
driven by the Poisson point measures $P^b_i$ and $P^d_i$: for all $1\leq i\leq d$ and all $t\geq 0$,
\begin{align}
  N^K_i(t) & = N^K_i(0)+\iiint_{[0,t]\times\N\times\R_+}\mathbbm{1}_{k\leq N^K_i(s-)}\mathbbm{1}_{y\leq b_i(E^K(s))}P^b_i(ds,dk,dy)
  \notag \\ & -\iiint_{[0,t]\times\N\times\R_+}\mathbbm{1}_{k\leq N^K_i(s-)}\mathbbm{1}_{y\leq d_i(E^K(s))+\frac{1}{K}\sum_{j=1}^d
  c_i(E^K(s),j)N^K_j(s-)}P^d_i(ds,dk,dy)  \label{eq:construction}
\end{align}
and
\[
  E^K(t)=E^K(0)+\dfrac{1}{\log K}\int_0^t F\left(\frac{1}{K}N^K(s),E^K(s)\right)ds.
\]
Note that this formulation can be understood as an algorithmic construction of the process, where successive jumps are constructed
inductively, solving the ODE~\eqref{eq:environment} between jumps and waiting for the first atom in the Poisson points measures for
which the indicator functions in~\eqref{eq:construction} are non-zero. Almost surely, there is no two such atoms at the same time,
and by standard arguments developed in~\cite{fournier-meleard2004}, under Assumption~\ref{hyp:regularity-Ltka-Volterra} there is almost surely no accumulation of jump times, so
that~\eqref{eq:construction} is valid for all time $t\geq 0$. In addition, the martingale problem associated to the process
$(N^K(t),E^K(t))_{t\geq 0}$ can be easily formulated following the method of~\cite{tran2008,campillo2015} (which
extend~\cite{fournier-meleard2004} to PDMP settings) and it is easy to check following~\cite{fournier-meleard2004} that, if the
initial condition is deterministic, the process $(N^K(t),E^K(t))_{t\geq 0}$ is locally bounded in $L^p$ for all $p\geq 1$.


\subsection{The total population is of the order of $K$}
\label{sec:total-pop}

The next result shows that the total population size remains of the order of $K$ for long times.

\begin{lem}
  \label{lem:uniform-bound} Take Assumption~\ref{hyp:regularity-Ltka-Volterra} for granted.
  \begin{description}
  \item[\textmd{(i)}] If there exists a constant $M>0$ such that, with probability
    converging to 1 when $K\to\infty$,
    \[
      \frac{1}{K}\sum_{i=1}^d N^K_i(0)\leq M,
    \]
    then there exist constants $R, V>0$ such that $\lim_{K\to +\infty}\mathbb{P}(\Gamma^K_{R,V})=1$, where
    \[
      \Gamma^K_{R,V}:=\left\{\frac{1}{K}\sum_{i=1}^d N^K_i(t)\leq R,\ \forall t\in[0,e^{V K}]\right\}.
    \]
  \item[\textmd{(ii)}] If there exist $\varepsilon,M>0$ such that, with probability converging to 1,
    \[
      \frac{1}{K}\sum_{i=1}^d N^K_i(0)\geq \varepsilon\quad\text{and}\quad |E^K(0)|\leq M,
    \]
    then, for all $T>0$, there exists a constant $\eta=\eta(M,\varepsilon,T)>0$ such that
    \[
      \lim_{K\to+\infty}\P\left(\frac{1}{K}\sum_{i=1}^d N^K_i(t)\geq\eta,\ \forall t\in[0,T\log K]\right)=1.
    \]
  \end{description}
\end{lem}

Point~(i) follows from an argument in~\cite{champagnat06} and Point~(ii) can be proved using similar ideas. Therefore, we only give a
sketch of proof below. Theses two points imply in particular that $\max_{1\leq i\leq d}\beta^K_i(t)$ converges to 1 in $L^\infty([0,T])$ in probability. Note that Point (ii) will not be needed to prove our main results (it is actually implied by Theorem~\ref{thm:BDP-max-2}). However we give an independent proof of this result for completeness.

Note that the timescale in (ii) is much shorter than the one in~(i).  This is due to the fact that, in
Assumption~\ref{hyp:regularity-Ltka-Volterra}, we only assume that $b_i(e)-d_i(e)>0$ for all $e\in\R^\dimE$. If instead we had assumed
that $b_i(e)-d_i(e)\geq r>0$ for all $e\in\R^\dimE$, then the proof below would entail~(ii) on the time interval $[0,e^{VK}]$ for some
constant $V>0$ instead of $[0,T\log K]$.

\begin{proof}
  We first prove~(i). Introduce two additional Poisson point measures $\widetilde{P}^b$ and $\widetilde{P}^d$ on $\R_+\times\N\times\R_+$ with intensity $\text{Leb}\times\sum_{n\geq 1}\delta_n\times\text{Leb}$ independent from $(N^K(0),E^K(0))$ and all the Poisson point measures $P^b_i,P^d_i$ for $1\leq i\leq d$. Define the process 
  \begin{align*}
    Y^K_t & =\sum_{i=1}^d N^K_i(0)+\sum_{i=1}^{d}\bigg\{\iiint_{[0,t]\times\N\times\R_+}\mathbbm{1}_{k\leq N^K_i(s-)}\mathbbm{1}_{y\leq \bar{b}}P^b_i(ds,dk,dy)
    \\ & \qquad\qquad\qquad-\iiint_{[0,t]\times\N\times\R_+}\mathbbm{1}_{k\leq N^K_i(s-)}\mathbbm{1}_{y\leq \underline{c}Y^K_{s-}/K}P^d_i(ds,dk,dy)\bigg\} \\
    & +\iiint_{[0,t]\times\N\times\R_+}\mathbbm{1}_{k\leq Y^K_{s-}-N^K_1(s-)-\ldots-N^K_{d}(s-)}\mathbbm{1}_{y\leq \bar{b}}\widetilde{P}^b(ds,dk,dy)
    \\ & -\iiint_{[0,t]\times\N\times\R_+}\mathbbm{1}_{k\leq Y^K_{s-}-N^K_1(s-)-\ldots-N^K_{d}(s-)}\mathbbm{1}_{y\leq \underline{c}Y^K_{s-}/K}\widetilde{P}^d(ds,dk,dy).
  \end{align*}
  It is easy to prove by induction on the successive jumps of $Y^k$ and $N^K_1+\ldots+N^K_d$ that for all $t\geq 0$
  \[
    N^K_1(t)+\ldots+N^K_d(t)\leq Y^K_t.
  \]
  This construction also implies that $Y^K$ is a logistic birth and death process with birth rate $\bar{b}$ and competition rate
  $\underline{c}/K$ (this means that $Y^k$ jumps from $k$ to $k+1$ at rate $k\bar{b}$ and from $k$ to $k-1$ at rate
  $\underline{c}k^2/K$). We can then use the large deviations principle for $Y^K/K$ and the problem of exit from an attracting domain
  to deduce~(i) (see~\cite[Thm.\ 3 (a) and (c)]{champagnat06}).

  To prove~(ii), we fix $T>0$ and we first notice that, on the event $\{|E^K(0)|\leq M\}$, $|E^K(t)|\leq M+\|F\|_\infty T$ for all
  $t\leq T\log K$ and, by continuity of $b$ and $d$, there exists $\underline{r}>0$ such that $b_i(e)-d_i(e)\geq \underline{r}$ for
  all $i$ and all $|e|\leq M+\|F\|_\infty T$. Let us introduce $(\widetilde{N}^K(t),\widetilde{E}^K(t))_{t\geq 0}$ the process obtained from the same construction as $(N^K(t),E^K(t))_{t\geq 0}$ but setting $\widetilde{E}^K(t)=E^K(T\log K)$ for all $t\geq T
  \log K$ and define $(\tau_k)_{k\geq 1}$ and $(\widetilde{\tau}_k)_{k\geq 1}$ the increasing sequence of jump times of the
  processes $N^K_1+\ldots+N^K_d$ and $\widetilde{N}^K_1+\ldots+\widetilde{N}^K_d$ respectively. Define also the induced chains
  \[
    X^K_k=N^K_1(\tau_k)+\ldots+N^K_d(\tau_k)\quad\text{and}\quad \widetilde{X}^K_k=\widetilde{N}^K_1(\tau_k)+\ldots+\widetilde{N}^K_d(\tau_k).
  \]
  Note that $\tau^K_k=\widetilde{\tau}^K_k$ and $X^K_k=\widetilde{X}^K_k$ for all $k$ such that $\tau_k<T\log K$.
  Then, for $k\geq 0$, conditionally on $\tau_k$, denoting by $\phi(t,e)$ the flow of the vector field
  $\dfrac{1}{\log K}F(\frac{1}{K}\widetilde{N}^K(\tau_k),\cdot)$ starting from $e$, $\rho(s)=\lambda(s)-\mu(s)$ and $\gamma(s)=\lambda(s)+\mu(s)$, where
  \[
    \lambda(s)=\sum_{i=1}^d \widetilde{N}^K_i(\widetilde{\tau}_k)b_i(\phi(s,\widetilde{E}^K_{\widetilde{\tau}_k}))
  \]
  and
  \[
    \mu(s)=\sum_{i=1}^d \widetilde{N}^K_i(\widetilde{\tau}_k)\left[d_i(\phi(s,\widetilde{E}^K_{\widetilde{\tau}_k}))+\frac{1}{K}\sum_{j=1}^d c_i(\phi(s,\widetilde{E}^K_{\widetilde{\tau}_k}),j)\widetilde{N}^K_j(\widetilde{\tau}_k)\right],
  \]
  straightforward computations show that
  \begin{align*}
    \P(\widetilde{X}^K_{k+1}=\widetilde{X}^K_k+1\mid \mathcal{F}_{\widetilde{\tau}_{k}})
    & -\P(\widetilde{X}^K_{k+1}=\widetilde{X}^K_k-1\mid\mathcal{F}_{\widetilde{\tau}_{k}}) \\ &
                                                   =\int_0^{+\infty}\rho(s\wedge(T\log K-\widetilde{\tau}_k))e^{-\int_0^s\gamma(u\wedge(T\log K-\widetilde{\tau}_k))du}ds
    \\ & \geq\frac{\underline{r}-\bar{c} \widetilde{X}^K_k/K}{\bar{b}+\bar{d}+\bar{c} \widetilde{X}^K_k/K}.
  \end{align*}
  In particular, on the event $\{\widetilde{X}^K_k\leq K \underline{r}/2\bar{c}\}$,
  \[
    \P(\widetilde{X}^K_{k+1}=\widetilde{X}^K_k+1\mid \mathcal{F}_{\widetilde{\tau}_{k}})
    \geq\frac{1}{2}+\frac{1}{2}\frac{\underline{r}}{2\bar{b}+2\bar{d}+\underline{r}}.
  \]
  This implies that, if $\widetilde{X}^K_k\leq K \underline{r}/2\bar{c}$ for some $k\geq 0$, there exists a biased random walk
  $(Y_\ell)_{\ell\geq k}$ with parameter $\frac{1}{2}+\frac{1}{2}\frac{\underline{r}}{2\bar{b}+2\bar{d}+\underline{r}}$ such that
  $Y_k=\widetilde{X}^K_k$ and $Y_\ell\leq \widetilde{X}^K_\ell$ for all $k\leq \ell\leq \tau$, where
  \[
    \tau=\inf\{\ell\geq k: \widetilde{X}^K_\ell> K \underline{r}/2\bar{c}\}.
  \]
  Setting $\eta=\underline{r}/4\bar{c}$, large deviations estimates on the biased random walk $Y$ show that each excursion of $\widetilde{X}^K$ below
  $2K\eta$ has a probability smaller than $e^{-VK}$ for some constant $V>0$ to exit $[K\eta,2K\eta]$ from below. Point~(ii) then
  follows from standard arguments on the problem of exit from a domain and domination of the number of
jumps of $N^K$ on $[0,T\log K]$.
\end{proof}


\subsection{Uniform continuity of the exponents}
\label{sec:preliminary}


Our goal here is to prove uniform continuity estimates for the processes $\beta^K_i(t)$.

\begin{lem}[Uniform continuity of the exponents]
  \label{lem:continuity-exponents}
  Let $i\in\{1,\dots, d\}$ and take Assumption~\ref{hyp:regularity-Ltka-Volterra} for granted. Assuming that, for some $M>0$, with probability
  converging to 1 when $K\to\infty$, $\frac{1}{K}(N^K_1(0)+\ldots+N^K_d(0))\leq M$ and $N^K_i(0)\geq K^\delta$ for some $0<\delta<1$, then there exists a constant $C_\text{Lip}$ such that, for all $\eta>0$
  small enough,
    \begin{equation}
        \label{eq:Lipschitz}
            \lim_{K\to+\infty}\mathbb{P}\left(\sup_{t\in[0,\eta]}|\beta^K_i(t)-\beta^K_i(0)|\leq C_\text{Lip}\eta\right)=1.
    \end{equation}
\end{lem}

\begin{proof}
  Let $i$ be fixed. We use Lemma~\ref{lem:uniform-bound}~(i) to obtain constants $R,V$ such that $\mathbb{P}(\Gamma^K_{R,V})$
  converges to 1. For all $t\geq 0$, we define the processes $Y^K$ and $Z^K$ as solutions to
  \[
    Y^K(t) = N^K_i(0)-\iiint_{[0,t]\times\N\times\R_+}\mathbbm{1}_{k\leq Y^K(s-)}\mathbbm{1}_{y\leq \bar{d}+\bar{c}R}P^d_i(ds,dk,dy)
  \]
  and
  \[
    Z^K(t) = N^K_i(0)+\iiint_{[0,t]\times\N\times\R_+}\mathbbm{1}_{k\leq Z^K(s-)}\mathbbm{1}_{y\leq \bar{b}}P^b_i(ds,dk,dy).
  \]
  It is then straightforward to check by induction over the successive jumps of $N^K_i$, $Y^K$ and $Z^K$ that, almost surely on the
  event $\Gamma^K_{R,V}$, for all $t\leq e^{KV}$,
  \[
    Y^K(t)\leq N^K_i(t)\leq Z^K(t).
  \]
  The process $Y^K$ is a linear pure death process with individual rate $\bar{d}+\bar{c}R$ and $Z^K$ is a linear pure birth (Yule)
  process with birth rate $\bar{b}$. Conditionally on $N^K_i(0)$, for all $t\geq 0$, it is well-known that $Y^K_t$ follows a binomal
  distribution with parameters $(N^K_i(0),e^{-(\bar{d}+\bar{c}R)t})$ and $Z^K_t$ a negative binomial distribution of parameters
  $(N^K_i(0),e^{-\bar{b}t})$. In particular,
  \begin{gather*}
    \E(Y^K_t\mid N^K_i(0))=N^K_i(0)e^{-(\bar{d}+\bar{c}R)t},\quad\text{Var}(Y^K_t\mid
    N^K_i(0))=N^K_i(0)e^{-(\bar{d}+\bar{c}R)t}\left(1-e^{-(\bar{d}+\bar{c}R)t}\right) \\
    \E(Z^K_t\mid N^K_i(0))={N^K_i(0)e^{\bar{b}t}} \quad\text{and}\quad\text{Var}(Z^K_t\mid
    N^K_i(0))=N^K_i(0)e^{\bar{b}t}\left(e^{\bar{b}t}-1\right).
  \end{gather*}
  Set $\eta_0:=\frac{\delta}{3(\bar{b}\vee(\bar{d}+\bar{c}R))}$. For all $\eta<\eta_0$, by Chebyshev's inequality, for $K$ large enough,
  \[
    \P\left(Y^K_{\eta\log
      K}<K^{-2(\bar{d}+\bar{c}R)\eta}N^K_i(0)\mid N^K_i(0)\right)\leq\frac{K^{(\bar{d}+\bar{c}R)\eta}}{N^K_i(0)(1-K^{-(\bar{d}+\bar{c}R)\eta})}\leq\frac{2K^{\delta/3}}{N^K_i(0)}
  \]
  and
  \[
    \P\left(Z^K_{\eta\log K}>\left(K^{\bar{b}\eta}+1\right)N^K_i(0)\mid N^K_i(0)\right)\leq\frac{K^{\bar{b}\eta}(K^{\bar{b}\eta}-1)}{N^K_i(0)}\leq\frac{K^{2\delta/3}}{N^K_i(0)}.
  \]
  Set $C_\text{Lip}:=2\left(\bar{b}\vee(\bar{d}+\bar{c}R)\right)$. Then, for all $\eta\leq \eta_0$ and $K$ sufficiently large,
  \begin{multline*}
    \mathbb{P}\left(\sup_{t\in[0,\eta]}|\beta^K_i(t)-\beta^K_i(0)|> C_\text{Lip}\eta\right) \\
    \begin{aligned}
      & \leq\mathbb{P}\left((\Gamma^K_{R,V})^c\cup\{N^K_i(0)<K^\delta\}\right) \\
      & +\P\left(\{N^K_i(0)\geq K^\delta\}\cap\left\{\log(1+Z^K_{\eta\log K})>\log(1+N^K_i(0))+2\bar{b}\eta\log K\right\}\right) \\
      & +\P\left(\{N^K_i(0)\geq K^\delta\}\cap\left\{\log(1+Y^K_{\eta\log
            K})<\log(1+N^K_i(0))-2(\bar{d}+\bar{c}R)\eta\log K\right\}\right) \\ &
      \leq\mathbb{P}\left((\Gamma^K_{R,V})^c\cup\{N^K_i(0)<K^\delta\}\right)+\P\left(\{N^K_i(0)\geq K^\delta\}\cap\left\{ Z^K_{\eta\log
            K}>(K^{\bar{b}\eta}+1)N^K_i(0)\right\}\right) \\ & +\P\left(\{N^K_i(0)\geq K^\delta\}\cap\left\{Y^K_{\eta\log
            K}<K^{-2(\bar{d}+\bar{c}R)\eta}N^K_i(0)\right\}\right).
    \end{aligned}
  \end{multline*}
  Hence~\eqref{eq:Lipschitz} follows.
\end{proof}

\subsection{Continuity estimates of the hitting time of 0 for exponents}
\label{sec:continuity-hitting-0}

{The next lemma will be used later to prove the continuity of hitting times of $0$ of approximations (via birth-death processes) of the limit dynamics of the exponents.}

\begin{lem}
\label{lem_control_reaching_time_of_0}
  Let $h_0 \in (0,1)$, consider $f: [0, +\infty) \mapsto \mathbb{R}$ a continuous function, and $D^-,D^+$ two positive constants. We define for all $\varepsilon>0$,
  \begin{align*}
    T_0 &:=\inf \left\{t\geq 0,\:h_0+\int_0^tf(s)ds=0\right\},
    \\
    T_\varepsilon^-&:=\inf\left\{t\geq 0:\:h_0+\int_0^tf(s)ds - \eps (D^- t+1)=0\right\},
    \\
    T_\varepsilon^+&= \inf\left\{t\geq 0:\: h_0+\int_0^tf(s)ds + \eps (D^+t+1)=0\right\}.
  \end{align*}
  If $T_0<\infty$ and if $f(T_0)<0$ then, for all $\varepsilon>0$ sufficiently small, $T_\varepsilon^-< T_0 < T_\varepsilon^+$,
  \[
    |T_0-T_\varepsilon^-|\leq \frac{4((D^+\vee D^-)2T_0+1))}{|f(T_0)|} \eps
    \quad \text{and} \quad
    |T_0-T_\varepsilon^+|\leq \frac{4((D^+\vee D^-)2T_0+1))}{|f(T_0)|} \eps.
  \]
\end{lem}

\begin{proof} 
  For all $t\geq 0$ and $\varepsilon>0$, let $h(t):=h_0+\int_0^tf(s)ds$ and $h^\pm_\varepsilon(t):=h(t)\pm\varepsilon(D^\pm t+1)$ and
  observe that $h^-_\varepsilon< h< h^+_\varepsilon$. Since $f(T_0)<0$, there exists $\delta>0$ such that, for all
  $t\in[T_0-\delta,T_0+\delta]$,
  \[
    f(t)<\frac{f(T_0)}{2}<0
  \]
  and we can assume without loss of generality that $\delta< T_0$.
  Define in addition $\gamma:=\inf\{h(t):t\in[0,T_0-\delta]\}>0$. Since $h(T_0)=0$, we deduce that, for all $t\in[T_0-\delta,T_0]$,
  \[
    h(t)\geq |t-T_0|\frac{|f(T_0)|}{2}
  \]
  and for all $t\in[T_0,T_0+\delta]$,
  \[
    h(t)\leq -|t-T_0|\frac{|f(T_0)|}{2}.
  \]
  Therefore, assuming that $\varepsilon<\frac{\delta}{A}\wedge\frac{\gamma}{D^- T_0+1}$ with
  \[
    A:=\frac{4((D^+\vee D^-)2T_0+1)}{|f(T_0)|},
  \]
  observing that $T_0+A\varepsilon\leq T_0+\delta\leq 2 T_0$, we have $h^+_\varepsilon(T_0+A\varepsilon)<0$ and
  $h^+_\varepsilon(t)>0$ for all $t\leq T_0$, hence $T^+_\varepsilon\in(T_0,T_0+A\varepsilon)$. Similarly, $h^-_\varepsilon(t)>0$ for
  all $t\in[T_0-\delta,T_0-A\varepsilon]$, $h^-_\varepsilon(T_0)<0$ and $h^-_\varepsilon(t)>0$ for all $t\leq T_0-\delta$ (using that
  $(D^- T_0+1)\eps\leq\gamma$), hence $T^-_\varepsilon\in(T_0-A\varepsilon,T_0)$. The conclusion follows.
\end{proof}

\subsection{Linear time inhomogeneous birth-death process}
\label{subsec.LBD}

The goal of this section is to prove convergence of the exponents for linear time inhomogeneous birth-death processes.

For a generic function $\varphi:\R_+\to\R_+$, we introduce the notation $\varphi_K(t)=\varphi(t/\log K)$. Given $x\in\R_+$ and
continuous functions $b$ and $d$ from $[0,+\infty)$ to $[0,+\infty)$, we denote by $LBD(b_K,d_K,x)$ the law of the time-inhomogeneous
linear birth and death process with generator at time $t$
\[
L_t f(k)=kb\left(\frac{t}{\log K}\right)(f(k+1)-f(k))+kd\left(\frac{t}{\log K}\right)(f(k-1)-f(k))
\]
and with initial condition $\lfloor x\rfloor$, the integer part of $x$ at time 0. A process $(Z^K_t)_{t\geq 0}$ with law $LBD(d_K,b_K,x)$ can be constructed as
\begin{align*}
  Z^K_t & = \lfloor x\rfloor+\iiint_{[0,t]\times\N\times\R_+}\mathbbm{1}_{k\leq Z^K_{s-}}\mathbbm{1}_{y\leq b(s/\log K)}P^b(ds,dk,dy)
   \\ & -\iiint_{[0,t]\times\N\times\R_+}\mathbbm{1}_{k\leq Z^K_{s-}}\mathbbm{1}_{y\leq d(s/\log K)}P^d(ds,dk,dy),
\end{align*}
with $P^b(ds,dk,dy)$ and $P^d(ds,dk,dy)$ two independent Poisson point measures on $\R_+\times\N\times\R_+$ with intensity measure
$\text{Leb}\times\sum_{n\geq 1}\delta_n\times\text{Leb}$.

\begin{lem}
  \label{lem:time-inhomo-LBD}
  \begin{description}
  \item[\textmd{(i)}] Let $Z_t^K$ follow the law $LBD(b_K,d_K,K^\beta)$, where $b,d:\R_+\to\R_+$ are continuous, nonnegative and
    bounded functions, and $\beta>0$. Define $r(s)=b(s)-d(s)$, and
    \begin{align*}
      T_0:=\inf\left\{t\geq0,\beta+\int_0^tr(s)\,ds=0\right\}.
    \end{align*}
    Then for all $T<T_0$, $\left(\frac{\log(1+Z_{t\log K}^K)}{\log K}\right)_{t\in[0,T]}$ converges in probability in
    $L^\infty([0,T])$ when $K\to\infty$ to $\left(\beta+\int_0^{t}r(s)\,ds\right)_{t\in[0,T]} $.
  \item[\textmd{(ii)}] If in addition $T_0<\infty$ and $b(T_0)-d(T_0)<0$, then
    $\left(\frac{\log(1+Z_{t\log K}^K)}{\log K}\right)_{t\in[0,T]}$ converges in probability in $L^\infty([0,T])$ when $K\to\infty$
    for all $T>0$ to $\left(\beta+\int_0^{t\wedge T_0}r(s)\,ds\right)_{t\in[0,T]} $ and for all $t > T_0$
    \begin{equation}
      \label{eq:true-extinction-LBP}
      \lim_{K\rightarrow +\infty}\P\left(Z^K_{t\log K}=0\right)=1.
    \end{equation}
  \end{description}
\end{lem}

\begin{proof}
  \noindent\textit{Step 1. Preliminary estimate.}
  In view of the construction above,
  $$
  Z_t^K = \lfloor K^\beta\rfloor + M_t^K + \int_0^t r\left(\frac{s}{\log K}\right)Z_s^K\,ds,
  $$
  where $M^K$ is a square integrable martingale with predictable quadratic variation~\cite{jacod-shiryaev,ikeda-watanabe}
  $$
  \langle M^K\rangle_t = \int_0^t \left(b\left(\frac{s}{\log K}\right)+d\left(\frac{s}{\log K}\right)\right)Z_s^K\,ds.
  $$
  Note also that
  $$
  \mathbb E(Z_t^K) = \lfloor K^\beta\rfloor \exp\left(\int_0^t r(s/\log K)\,ds\right).
  $$
  Denoting for all $t\geq 0$
  \begin{align*}
    r^{K}(t):= \int_{0}^t r \left ( \dfrac{s}{\log K}\right)ds=\log K \int_{0}^{t/\log K} r(s)ds
  \end{align*}
  and using It\^o's formula, then
  $$
  e^{-r^K(t)}Z_t^K = \lfloor K^\beta\rfloor+\int_0^t e^{-r^K(s)}dM_s^K=:\lfloor K^\beta\rfloor+\widehat M_t^K,
  $$
  where $\widehat M_t^K$ is a square integrable martingale with predictable quadratic variation
  $$
  \langle\widehat M^K\rangle_t = \int_0^t\left( b\left(\frac{s}{\log K}\right)+d\left(\frac{s}{\log K}\right)\right)Z^K_se^{-2r^K(s)}ds.
  $$
  Using Doob's inequality, for $0<\eta<\beta$:
  \begin{align*}
    \mathbb P\left(\sup_{t\leq T\log K} \left|e^{-r^K(t)}Z_t^K- K^\beta\right|\geq 2K^\eta\right) & \leq\mathbb P\left(\sup_{t\leq T\log K} \left|e^{-r^K(t)}Z_t^K-\lfloor K^\beta\rfloor\right|\geq K^\eta\right) \\ &
      \leq K^{\beta-2\eta}\int_0^{T\log K}\left( b\left(\frac{s}{\log K}\right)+d\left(\frac{s}{\log K}\right)\right)e^{-r^K(s)}ds \\ &
      \leq K^{2(\beta-\eta)}( \|b\|_\infty+\|d\|_\infty) \int_0^{T\log K}K^{-\beta-\int_0^{s/\log K} r(u)du}ds.
  \end{align*}

\medskip

\noindent\textit{Step 2. Proof of~(i).} 
Assume $T<T_0$. Then there exists $\varepsilon>0$ small enough such that $T<T_\varepsilon$, where
$$
T_\varepsilon:=\inf\left\{t>0,\, \beta+\int_0^t r(s)ds \leq \varepsilon\right\}.
$$
Without loss of generality, we can assume $\beta>\varepsilon>0$.
Let $\eta\in(\beta-\eps/2,\beta)$ and define the set
\begin{equation}
  \label{eq:def-Omega-eps}
  \Omega_{1,\eps}^K:=\left\{\sup_{t\leq (T_\eps\wedge T)\log K}
    \left|e^{-r^{K}(t)}Z_t^K- K^\beta \right|\leq 2K^\eta\right\}.
\end{equation}
We have from Step~1 that
$$
\mathbb P((\Omega_{1,\eps}^K)^c) \leq K^{2(\beta-\eta)}(||b||_\infty+||d||_\infty)K^{-\eps} T\log(K)\xrightarrow{K\rightarrow\infty}0.
$$
In addition,
\begin{multline*}
  \sup_{t\leq T} \left|\frac{\log(1+Z^K_{t\log K})}{\log K}- \left(\beta +\int_0^tr(s)ds\right)\right|\\
  \begin{aligned}
    = &\ \sup_{t\leq T}\left|\frac{\log(1+Z^K_{t\log K})}{\log K} - \frac{\log(1+K^{\beta+\int_0^tr(s)ds})}{\log K} + \frac{\log(1+K^{\beta+\int_0^tr(s)ds})}{\log K}-\frac{\log(K^{\beta+\int_0^tr(s)ds})}{\log K}\right| \\
    \leq &\  \frac{1}{\log K}\sup_{t\leq T} \left[\left|\log\Big(\frac{1+Z^K_{t\log K}}{1+K^{\beta+\int_0^tr(s)ds}}\Big)\right| +\log\big(1+K^{-\beta-\int_0^tr(s)ds}\big)\right]\\
    = &\   \frac{1}{\log K}\sup_{t\leq T} \left[\log\left(1+\frac{\left|Z^K_{t\log K}- K^{\beta+\int_0^tr(s)ds}\right|}{1+Z^K_{t\log K}\wedge
          K^{\beta+\int_0^tr(s)ds}}\right) +\log\big(1+K^{-\beta-\int_0^tr(s)ds}\big) \right] \\
    \leq &\  \frac{1}{\log K}\sup_{t\leq T}\left[\frac{|Z^K_{t\log K}-K^{\beta+\int_0^tr(s)ds}|}{1+Z^K_{t\log K}\wedge
        K^{\beta+\int_0^tr(s)ds} }+K^{-\beta-\int_0^tr(s)ds}\right]
  \end{aligned}
\end{multline*}
Hence, 
recalling that $T<T_\varepsilon$, on the event $\Omega_{1,\eps}^K$, for $K$ large enough,
\begin{multline*}
  \sup_{t\leq T} \left|\frac{\log(1+Z^K_{t\log K})}{\log K}- \left(\beta +\int_0^tr(s)ds\right)\right|\\
    \leq   \frac{1}{\log K}\,\sup_{t\leq T}\frac{2K^{\eta+\int_0^tr(s)ds}}{1+K^{\beta+\int_0^tr(s)ds}-2 K^{\eta+\int_0^tr(s)ds}} + \frac{K^{-\varepsilon}}{\log K} 
    \leq \frac{4K^{-(\beta-\eta)}+K^{-\varepsilon}}{\log K}, 
\end{multline*}
which converges to zero when $K\rightarrow +\infty$. Hence~(i) is proved.

\medskip

\noindent\textit{Step 3. Proof of~\eqref{eq:true-extinction-LBP}.}
Assume that $T_0<+\infty$ and $r(T_0)=b(T_0)-d(T_0)<0$ and let $t>T_0$. Define the constant $a:=2\|r\|_\infty/|r(T_0)|$ and the time
\[
T_{-\varepsilon}=\inf\left\{t\geq 0, \beta+\int_0^tr(s)ds\leq -\eps\right\}.
\]
It follows from Lemma~\ref{lem_control_reaching_time_of_0} (with $D^+=D^-=0$) that there exists a constant $C$ such that, for
$\varepsilon$ small enough,
\begin{equation}
    \label{eq:bound_T_eps}
    0\leq T_{-a\varepsilon}-T_0\leq C\varepsilon\quad \text{and}\quad 0\leq T_0-T_\varepsilon\leq C\varepsilon.
\end{equation}
Reducing $\varepsilon$ if needed, we can also assume that $t>T_{-a\varepsilon}$ and
\[
\bar{b}:=\sup_{s\in[T_\varepsilon,T_{-a\varepsilon}]}b(s)<\inf_{s\in[T_\varepsilon,T_{-a\varepsilon}]}d(s)=:\underline{d}
\]
are such that
\[
  \underline{r}:=\bar{b}-\underline{d}<\frac{r(T_0)}{2}<0.
\]
We apply Step 2 for this value of $\varepsilon$, fixed until the end of Step 3. In particular, given
$\eta\in(\beta-\varepsilon/2,\beta)$, the set $\Omega^K_{1,\varepsilon}$ defined in~\eqref{eq:def-Omega-eps} has probability
converging to 1.

Similar comparison arguments as those we used in the proof of Lemma~\ref{lem:continuity-exponents} allow us to couple the process
$Z^K$ with a process $\bar{Z}^K$ such that $\bar{Z}^K_{0} =Z^K_{T_\varepsilon\log K }$,
$Z^K_t\leq \bar{Z}^K_{t-T_\varepsilon\log K }$ for all $t\in[T_\varepsilon\log K ,T_{-a\varepsilon}\log K ]$ and $\bar{Z}^K$ is a
time-homogeneous linear birth and death process with (constant) individual birth rate $\bar{b}$ and death rate $\underline{d}$. It
follows from the definition of $T_\varepsilon$ that a.s., on the event $\Omega_{1,\eps}^K$,
$\bar{Z}^K_0\leq K^\varepsilon+K^{\eta+\varepsilon-\beta}\leq 2 K^\varepsilon$.

Remind from~\cite[Section 5.4.5, p. 180]{Meleardbook} that, if $\bar{T}_{\text{ext}}$ denotes the first hitting time of 0 by $\bar{Z}^K$,
  \begin{align}
    \label{eq:Sylvie}
    \mathbb{P}(\bar{T}_{\text{ext}}>t\mid \bar{Z}^K_0=1)=\frac{\underline{r} e^{\underline{r}t}}{\bar{b} e^{\underline{r}t}-\underline{d}}.
  \end{align}
 Hence, the branching property for $\bar{Z}^K$ entails that
  \begin{align*}
    \mathbb{P}(\bar{T}_{\text{ext}}>t\mid \bar{Z}^K_0\leq 2 K^\varepsilon)\leq 1-\left(1-\frac{\underline{r} e^{\underline{r}t}}{\bar{b} e^{\underline{r}t}-\underline{d}}\right)^{2K^\varepsilon}.
  \end{align*}
  Applying this inequality for $t=(T_{-a\varepsilon}-T_\varepsilon)\log K$ and observing that
  \[
    (a+1)\varepsilon=-\int^{T_{-a\varepsilon}}_{T_\varepsilon}r(s)ds=\int^{T_{-a\varepsilon}}_{T_\varepsilon}|r(s)|ds\leq\|r\|_\infty(T_{-a\varepsilon}-T_{\varepsilon}),
  \]
  we deduce that, for $K$ large enough,
  \begin{align*}
    \mathbb{P}\left(\bar{T}_{\text{ext}}>(T_{-a\varepsilon}-T_\varepsilon)\log K \mid \bar{Z}^K_0\leq 2 K^\varepsilon\right)
    & \leq 1-\left(1-\frac{2|\underline{r}|}{\underline{d}} K^{\underline{r}(T_{-a\varepsilon}-T_\varepsilon)}\right)^{2K^\varepsilon}
    \\ & \leq 1-\left(1-\frac{2|\underline{r}|}{\underline{d}} K^{-(a+1)|\underline{r}|\varepsilon/\|r\|_\infty}\right)^{2K^\varepsilon},
  \end{align*} 
  which converges to 0 as $K\to+\infty$ since $a$ was chosen such that $(a+1)|\underline{r}|>\|r\|_\infty$. Since
  \[
    \P(Z^K_{t\log K}\geq 1)\leq \P(Z^K_{T_{-a\varepsilon}\log K}\geq 1)\leq
    \P((\Omega^K_{1,\varepsilon})^c)+\P(\Omega^K_{1,\varepsilon}\cap\{\bar{T}_\text{ext}>(T_{-a\varepsilon}-T_\varepsilon)\log K\}),
  \]
  \eqref{eq:true-extinction-LBP} follows from the facts that $\bar{Z}^K_0\leq 2 K^\varepsilon$ on $\Omega^K_{1,\epsilon}$ and
  $\P(\Omega^K_{1,\varepsilon})\to 1$ when $K\to+\infty$.
    
  \medskip
  
  \noindent\textit{Step 4. Conclusion of the proof of~(ii).}
  Since the last two steps were true for any value of $\varepsilon>0$ small enough, in order to complete the proof, it is enough to
  check that, for $\varepsilon>0$ small enough, with probability converging to 1 when $K\to\infty$,
  \[
    \sup_{s\in[T_\varepsilon,T_{-a\varepsilon}]} \frac{\log(1+Z^K_{s\log K})}{\log K}\leq A\varepsilon
  \]
  for some constant $A>0$. Reminding~\eqref{eq:bound_T_eps}, it is enough to prove that
  \[
    \sup_{s\in[T_0-C\varepsilon,T_0+C\varepsilon]} \frac{\log(1+Z^K_{s\log K})}{\log K}\leq A\varepsilon.
  \]
  This is the kind of estimates provided by Lemma~\ref{lem:continuity-exponents} for the processes $N^K_i(t)$. We leave to the reader to
  check that its proof can be easily adapted to the process $Z^K$ considered above.
\end{proof}

\subsection{1-dimensional logistic time inhomogeneous birth-death process}
\label{sec:logistic}

We consider the setting of Section~\ref{sec:construction}. Consider any family of predictable processes $((Y^K_t)_{t\geq 0})_{K\geq 1}$
in $\R^\dimE$ on the filtered probability space $(\Omega,({\cal F}_t)_{t\geq 0},\P)$, which are almost surely continuous.
Given a ${\cal F}_0$-mesurable random variable $Z^K_0\in\N$ and bounded functions $b$, $d$, and $c$ from $\R^\dimE$ to
$[0,+\infty)$, we denote by $LogBD_K(b,d,c, Y^{K})$ the law of the logistic time-inhomogeneous linear birth and death process
$(Z^K_t)_{t\geq 0}$ in the environment $Y^K$ constructed for all $t>0$ as
\begin{align*}
  Z^K_t & = Z^K_0 +\iiint_{[0,t]\times\N\times\R_+}\mathbbm{1}_{k\leq Z^K_{s-}}\mathbbm{1}_{y\leq b(Y^K_s)}P^b(ds,dk,dy)
  \\ & -\iiint_{[0,t]\times\N\times\R_+}\mathbbm{1}_{k\leq Z^K_{s-}}\mathbbm{1}_{y\leq d(Y^K_s)+\frac{1}{K}c(Y^K_s)Z^K_{s-}}P^d(ds,dk,dy).  
\end{align*}
where $P^b$ is one of the Poisson point measures $(P^b_i)_{1\leq i\leq d}$ and $P^d$ is one of the $(P^d_i)_{1\leq i\leq d}$. This
construction means that 
$(Z^K_t)_{t\geq 0}$ represents a population dynamics where, at time $t\geq 0$, each individual gives birth to a new individual at
rate $b(Y^K_t)$ and dies at rate $d(Y^K_t)+\frac{Z^K_{t-}}{K}c(Y^K_t)$.

\begin{lem}
\label{lem:time-to-equi}
Consider a family of predictable processes $((Y^K_t)_{t\geq 0})_{K\geq 1}$ in $\R^\dimE$ almost surely continuous and differentiable on
the left and the right everywhere and such that $\vert \dot Y^{K}\vert \le C_Y/ \log K$ for some constant $C_Y$ and, with probability
converging to 1 when $K\to\infty$, $|Y^K_0|\leq C_Y$. Consider $b$, $d$ and $c$ functions of class $C^1_b$ from $\R^\dimE$ to $\R_+$ such that $b(x)>d(x)$
and $c(x)\geq \underline{c}$ for all $x\in\R^\dimE$ for some constant $\underline{c}>0$. Let $(Z_t^K)_{t\geq 0}$ with law
$LogBD_K(b,d,c, Y^{K})$ be constructed as above.
For all $T>0$, for $\eps>0$ small enough, assuming that, with probability converging to 1, $Z_0^K/K\in[u^\star(Y^{K}_0)-\eps,u^\star(Y^{K}_0)+\eps]$, where
$u^\star(x):=\frac{b(x)-d(x)}{c(x)}$ for all $x\in\R^\dimE$, then for all $T>0$
$$
\lim_{K\to\infty}\mathbb P\left(\forall t\in[0,T\log K],\,\,\left|\frac{Z_{t}^K}{K}-u^\star(Y^{K}_{t})\right|\leq 3\eps\right)= 1.
$$
\end{lem}

\begin{proof}
  %
  We follow the method of~\cite{champagnat2014}. Fix $T>0$. Observing that, for all $x\in\R^\dimE$, the logistic equation (with fixed
  $x$)
  $$
  \dot y = (b(x)-d(x)-c(x)y)y
  $$
  admits the Lyapunov function $V(y) = (y-u^\star(x))^2$, we define
  \[
    T_\text{exit}:=\inf\left\{t\geq 0:
    \left|\frac{Z_{t}^K}{K}-u^\star(Y^{K}_{t})\right|> 3\eps\right\}\wedge (T\log K)
  \]
  and
  \[
    U_t^K = \left(\frac{Z_{t\wedge T_\text{exit}}^K}{K}-u^\star\left(Y^K_{t\wedge T_\text{exit}}\right)\right)^2.
  \]
  Our goal is to prove that $T_\text{exit}=T\log K$ with probability converging to 1.

  By the assumptions on $Y^K$, conditionally on $\{Z_0^K/K\in[u^\star(Y^{K}_0)-\eps,u^\star(Y^{K}_0)+\eps]\}$, the process
  $Z^K$ admits the semimartingale decomposition
  \[
    Z^K_t=Z^K_0+\int_0^t \left(b(Y^K_s)-d(Y^K_s)-c(Y^K_s)\frac{Z^K_s}{K}\right)Z^K_s ds+M^K_t,
  \]
  where $M^K_t$ is a square integrable martingale with predictable quadratic variation
  \[
    \langle M^K\rangle_t=\int_0^t\left(b(Y^K_s)+d(Y^K_s)+c(Y^K_s)\frac{Z^K_s}{K}\right)Z^K_s ds.
  \]
  It then follows from It\^o's formula that, for all $t<T_\text{exit}$, 
  \begin{align*}
    U_t^K &= U_0^K+2\int_0^t\left(\frac{Z_s^K}{K}-u^\star(Y^K_s)\right)\left[\frac{Z_s^K}{K}\left(b(Y^K_s)-d(Y^K_s)-c(Y^K_s)\frac{Z_s^K}{K}\right)-(u^\star)'(Y^K_s)\dot{Y}^K_s\right]ds
    \\
          &+\frac{1}{K^2}\langle M^K\rangle_t+\widehat M_t^K \\ & = U_0^K-2\int_0^t\frac{Z^K_s}{K}c(Y^K_s) U^K_sds-2\int_0^t\left(\frac{Z_s^K}{K}-u^\star(Y^K_s)\right)(u^\star)'(Y^K_s)\dot{Y}^K_s ds+\frac{1}{K^2}\langle M^K\rangle_t+\widehat M_t^K,
  \end{align*}
  where $\widehat M_t^K:=\frac{2}{K}\int_0^t(\frac{Z_s^K}{K}-u^\star(Y^K_s))dM_s^K$ is a martingale with predictable quadratic variation
  satisfying for some constant $C_M$ independent of $K$ and all $t< T_\text{exit}$, with probability converging to 1,
  $$
  \langle \widehat M^K\rangle_t = \frac{4}{K^2}\int_0^t U_s^KZ_s^K\left(b(Y^K_s)+d(Y^K_s)+c(Y^K_s)\frac{Z_s^K}{K}\right)ds\leq \frac{C_M}{K}\int_0^t U_s^K\,ds.
  $$ 
  In the last inequality, we used the fact that, on the event $\{Y^K_0\leq C_Y\}$, $|Y^K_t|\leq C_Y(1+T)$ for all $t\in[0,T\log K]$, hence $u^\star(Y^K_s)$ is bounded and therefore 
  $Z^K_t/K$ is uniformly bounded for $t\leq T_\text{exit}$. 
  
For the same reason, on the event $\{Y^K_0\leq C_Y\}$ there exists a constant
  $\underline{u}^\star>0$ such that $u^\star(Y^K_t)\geq\underline{u}^\star$ for all $t\in[0,T\log K]$. Since $u^\star$ is Lipschitz,
  there exists a constant $C$ such that, for all $0\leq s\leq t< T_\text{exit}\leq (T\log K)$, for $\varepsilon>0$ small enough and
  $K$ large enough,
  \begin{align*}
    U_t^K & \leq U_s^K+C(t-s)\left(\frac{1}{K}+\frac{\varepsilon}{\log K}\right)-2\int_s^t\frac{Z_u^K}{K}c(Y^K_u)U_u^K\,du+\widehat
            M_t^K-\widehat M^K_s \\ & \leq U_s^K+\frac{2C}{\log K}(t-s)-\underline{c}\underline{u}^\star\int_s^t U_u^K\,du+\widehat
                                      M_t^K-\widehat M^K_s
  \end{align*}

  Now, define $\tau_0=0$ and, inductively for all $k\geq 1$
  \[
    \tau'_k:=\inf\left\{t\geq \tau_{k-1}: \left|\frac{Z^K_t}{K}-u^\star(Y^K_t)\right|> 2\varepsilon\right\}\wedge(T \log K)
  \]
  and
  \[
    \tau_k:=\inf\left\{t\geq \tau'_{k}: \left|\frac{Z^K_t}{K}-u^\star(Y^K_t)\right|\not\in [\varepsilon,3\varepsilon]\right\}\wedge(T \log K).
  \]
  Fix $k\geq 1$ such that $\tau'_k<T\log K$.
  Assuming $K$ is large enough such that $\log
  K>4C/(\underline{c}\underline{u}^\star\varepsilon^2)$, using that $\varepsilon^2\leq U^K_s\leq 9\varepsilon^2$ for all $s\in [\tau'_k,\tau_k)$, we have for all $s<t$ in $[\tau'_k,\tau_k]$,
  \begin{align*}
    U_t^K 
     & \leq U_s^K-\frac{\underline{c}\underline{u}^\star\varepsilon^2}{2}(t-s)+\widehat
                                      M_t^K-\widehat M^K_s
  \end{align*}
  and
  \begin{align*}
    \langle \widehat M^K\rangle_t-\langle \widehat M^K\rangle_s\leq \frac{9C_M\varepsilon^2}{K}(t-s).
  \end{align*}
  Setting $\alpha_K=K\frac{\underline{cu}^\star}{9C_M}$, it then follows from It\^o's formula that
  $(e^{\alpha_K U^K_{t\wedge \tau_k}})_{t\geq \tau'_k}$ is a square integrable supermartingale. In particular, reminding that the
  jumps of $Z^K/K$ are bounded by $1/K$, on the event $\{\tau'_k<T\log K\}$,
  \[
    \E\left(e^{\alpha_K U^K_{\tau_k}}\mid \mathcal{F}_{\tau'_k}\right)\leq e^{\alpha_K U^K_{\tau'_k}}\leq e^{\alpha_K(2\eps+K^{-1})^2}.
  \]
  Since
  \[
    \E\left(e^{\alpha_K U^K_{\tau_k}}\mid \mathcal{F}_{\tau'_k}\right)\geq \P\left(U^K_{\tau_k}> 3\varepsilon\mid \mathcal{F}_{\tau'_k}\right)
    e^{9\alpha_K\varepsilon^2}+\left(1-\P\left(U^K_{\tau_k}> 3\varepsilon\mid \mathcal{F}_{\tau'_k}\right)\right) e^{\alpha_K(\varepsilon-K^{-1})^2}.
  \]
  We deduce that, for $K$ large enough, on the event $\{\tau'_k<T\log K\}$,
  \[
    \P\left(U^K_{\tau_k}> 3\varepsilon\mid \mathcal{F}_{\tau'_k}\right)\leq
    \frac{e^{5\alpha_K\eps^2}}{e^{9\alpha_K\eps^2}-e^{\alpha_K\eps^2/2}}\leq
    \exp\left(-3K\frac{\underline{cu}^\star\eps^2}{9C_M}\right)
  \]
  The result then follows from standard arguments for the problem of exit of a domain~\cite{freidlin-wentzell,dembo-zeitouni}, that
  we sketch below for completeness: the previous bound tells us that we need $Geom(e^{-3K\underline{cu}^\star\eps^2/9C_M})$ trials
  before $\tau_k=T_\text{exit}$, and each of them requires at least one birth or death event, which occurs at a rate bounded from
  above by $(\bar b+\bar d)K$, hence $\tau_k-\tau'_k\geq 1/K$ with positive probability given $\mathcal{F}_{\tau'_k}$ and
  $\{\tau'_k<T\log K\}$.
  %
  %
\end{proof}

\subsection{2-dimensional logistic time inhomogeneous birth-death process}
\label{sec:competitive.Lotka.volterra}

We consider the setting of Section~\ref{sec:construction}. Consider three family of predictable process
$((Z^K_t)_{t\geq 0})_{K\geq 1}$ in $\R^\dimE$, continuous, $((\hat{d}^K_1(t))_{t\geq 0})_{K\geq 1}$ in $\R_+$ and
$((\hat{d}^K_2(t))_{t\geq 0})_{K\geq 1}$ in $\R_+$ on the filtered probability space $(\Omega,({\cal F}_t)_{t\geq 0},\P)$.
Given ${\cal F}_0$-mesurable random variables $X^K_0$ and $Y^K_0$ in $\N$ and bounded functions $b_i$, $d_i$ and $c_{ij}$ for
$i,j\in\{1,2\}$, all $\R^\dimE$ to $[0,+\infty)$, we denote by $LV2D_K(b,d,c,\hat{d}^K,Z^K)$, with $b=(b_1,b_2)$, $d=(d_1,d_2)$,
$c=(c_{11},c_{12},c_{21},c_{22})$ and $\hat{d}^K=(\hat{d}^K_1,\hat{d}^K_2)$, the law of the two-dimensional competitive
Lotka-Volterra process $(X^K_t,Y^K_t)_{t\geq 0}$ with birth rate $b$, death rate $d$, competition rate $c$ in the environment $(\hat{d}^K,Z^K)$
constructed for all $t>0$ as
\begin{align}
  X^K_t & = X^K_0 +\iiint_{[0,t]\times\N\times\R_+}\mathbbm{1}_{k\leq X^K_{s-}}\mathbbm{1}_{y\leq b_1(Z^K_s)}P^{X,b}(ds,dk,dy)
  \notag \\ & -\iiint_{[0,t]\times\N\times\R_+}\mathbbm{1}_{k\leq X^K_{s-}}\mathbbm{1}_{y\leq
       d_1(Z^K_s)+\hat{d}^K_1(s)+\frac{1}{K}c_{11}(Z^K_s)X^K_{s-}+\frac{1}{K}c_{12}(Z^K_s)Y^K_{s-}}P^{X,d}(ds,dk,dy), \label{eq:def-X-2D}
  \\
  Y^K_t & = Y^K_0 +\iiint_{[0,t]\times\N\times\R_+}\mathbbm{1}_{k\leq Y^K_{s-}}\mathbbm{1}_{y\leq b_2(Z^K_s)}P^{Y,b}(ds,dk,dy)
  \notag \\ & -\iiint_{[0,t]\times\N\times\R_+}\mathbbm{1}_{k\leq Y^K_{s-}}\mathbbm{1}_{y\leq
       d_2(Z^K_s) +\hat{d}^K_2(s)+\frac{1}{K}c_{21}(Z^K_s)X^K_{s-}+\frac{1}{K}c_{22}(Z^K_s)Y^K_{s-}}P^{Y,d}(ds,dk,dy), \label{eq:def-Y-2D}
\end{align}
where $P^{X,b}$ and $P^{Y,b}$ are two distinct Poisson point measures among $(P^b_i)_{1\leq i\leq d}$ and $P^{X,d}$ and $P^{Y,d}$ are
two distinct Poisson point measures among $(P^d_i)_{1\leq i\leq d}$. This construction means that 
$(X^K_t,Y^K_t)_{t\geq 0}$ represents a two-dimensional population dynamics where, at time $t\geq 0$, each individual in the
population $X^K_t$ (resp.\ $Y^K_t$) gives birth to a new individual in the same population at rate $b_1(Z^K_t)$ (resp.\ $b_2(Z^K_t)$
and dies at rate $d_1(Z^K_t)+\hat{d}^K_1(t)+\frac{X^K_{t}}{K}c_{11}(Z^K_t)+\frac{Y^K_{t}}{K}c_{12}(Z^K_t)$ (resp.\
$d_2(Z^K_t)+\hat{d}^K_2(t)+\frac{X^K_{t}}{K}c_{21}(Z^K_t)+\frac{Y^K_{t}}{K}c_{22}(Z^K_t)$).



\begin{lem}\label{lem:competition}
  Consider families of predictable processes $((Z^K_t)_{t\geq 0})_{K\geq 1}$ in $\R^\dimE$, $((\hat{d}^K_1(t))_{t\geq 0})_{K\geq 1}$ in
  $\R_+$ and $((\hat{d}^K_2(t))_{t\geq 0})_{K\geq 1}$ in $\R_+$ such that $Z^K$ is almost surely continuous and differentiable on the
  left and the right everywhere and such that, for all $T>0$, in probability, 
  \begin{align}
    \label{ass:vanishing_terms}
    \sup_{t\in[0,T]}|\hat{d}_1^K(t)|+|\hat{d}_2^K(t)|+\left|\frac{dZ^K_t}{dt}\right|\xrightarrow[K\to+\infty]{} 0.
  \end{align}
  Consider $b_i$, $d_i$ and $c_{ij}$ for $i,j\in\{1,2\}$, bounded functions of class $C^1_b$ from $\R^\dimE$ to $\R_+$ such that
  $b_i(z)>d_i(z)$ and $c_{ij}(z)\geq \underline{c}$ for all $i,j\in\{1,2\}$ and $z\in\R^\dimE$, for some contant $\underline{c}>0$.
  Assume also that there exists a compact subset $G$ of $\mathbb{R}^\dimE$ such that,
  \begin{equation}
    \label{eq:hyp-2D-Z}
    \lim_{K\to+\infty}\P\left(Z^K_0\in G\right)=1    
  \end{equation}
  and, for all $z\in G$, $b_2(z)-d_2(z)>\delta$,
   \begin{equation}
     \label{eq:IIF1}
     g_{21}(z):=b_2(z)-d_2(z)-c_{21}(z)\frac{b_1(z)-d_1(z)}{c_{11}(z)}\geq\delta
   \end{equation}
   and 
   \begin{equation}
     \label{eq:IIF2}
     g_{12}(z):=b_1(z)-d_1(z)-c_{12}(z)\frac{b_2(z)-d_2(z)}{c_{22}(z)}\leq-\delta
   \end{equation}
   for some $\delta>0$. Let $(X_t^K,Y^K_t)_{t\geq 0}$ with law $LV2D_K(b,d,c,\hat{d}^K,Z^K)$ be constructed as above and assume that
   there exists a compact subset $H$ of $(0,+\infty)^2$ such that
   \begin{equation}
     \label{eq:hyp-CI-lemme-2D}
     \lim_{K\to+\infty}\P\left(\frac{1}{K}(X^K_0,Y^K_0)\in H\right)=1.     
   \end{equation}
   Then, for all $\varepsilon>0$, there exists $\eta=\eta(\varepsilon,G,H)\in(0,\varepsilon)$
   and $T_0=T_0(\varepsilon,G,H)<+\infty$ such that 
    \begin{equation}
      \label{eq:ccl-lemme-2D}
      \lim_{K\rightarrow+\infty}
      \P\left(\frac{X^K_{T_0}}{K}\in[\eta ,\varepsilon ]\text{ and }\ 
      \left|\frac{Y^K_{T_0}}{K}-u^\star(Z^K_0)\right|\leq \varepsilon\right)=1,  
  \end{equation}
  where, for all $z\in\R^\dimE$, $u^\star(z):=\frac{b_2(z)-d_2(z)}{c_{22}(z)}$.
\end{lem}


\begin{proof}
  \textit{Step 1. Comparison of $\frac{1}{K}(X^K,Y^K)$ with deterministic flows.}
  Our goal in this step is to prove that, for all $t>0$, in probability,
  \begin{equation}
    \label{eq:step1-2D}
    \sup_{s\in[0,t]} \left|\frac{1}{K}(X^K_s,Y^K_s)-\phi_t\left(\frac{1}{K}X^K_0,\frac{1}{K}Y^K_0,Z^K_0\right)\right|\xrightarrow[K\to +\infty]{}0,    
  \end{equation}
  where $\phi_t(x,y,z)$ is the flow of the two-dimensional ODE (for $z$ fixed)
  \begin{equation}
    \label{eq:2D-LV}
    \begin{cases}
        \dot{x}(t)=x(t)(b_1(z)-d_1(z)-c_{11}(z)x(t)-c_{12}(z)y(t)), & \\ \dot{y}(t)=y(t)(b_2(z)-d_2(z)-c_{21}(z)x(t)-c_{22}(z)y(t)), &
    \end{cases}   
  \end{equation}
  i.e.\ $t\mapsto \phi_t(x,y,z)$ is the solution of the previous ODE with initial value $\phi_0(x,y,z)= (x,y)$.

  \medskip

  We adapt the general method of~\cite[Ch. 11]{ethier-kurtz-86}. Let us first notice that the argument allowing to prove the upper bound in
  Lemma~\ref{lem:uniform-bound} also applies in this case, so that there exists $V,R>0$ such that
  \begin{equation}
    \label{eq:bound-R-2D}
    \lim_{K\to+\infty}\P\left(\forall t\in[0,e^{VK}],\ X^K_t+Y^K_t\leq RK\right)=1.    
  \end{equation}
  Therefore, the process
  \begin{align*}
    \tilde{X}^K_t & = X^K_0+\iiint_{[0,t]\times\N\times\R_+}\mathbbm{1}_{k\leq \tilde{X}^K_{s-}\wedge (RK)}\mathbbm{1}_{y\leq b_1(Z^K_s)}P^{X,b}(ds,dk,dy)
    \\ & -\iiint_{[0,t]\times\N\times\R_+}\mathbbm{1}_{k\leq \tilde{X}^K_{s-}\wedge(RK)}\mathbbm{1}_{y\leq \tilde{d}^K_1(\frac{\tilde{X}^K_{s-}}{K}, \frac{\tilde{Y}^K_{s-}}{K},s)}P^{X,d}(ds,dk,dy), 
    \\
    \tilde{Y}^K_t & = Y^K_0+\iiint_{[0,t]\times\N\times\R_+}\mathbbm{1}_{k\leq \tilde{Y}^K_{s-}\wedge(RK)}\mathbbm{1}_{y\leq b_2(Z^K_s)}P^{Y,b}(ds,dk,dy)
    \\ & -\iiint_{[0,t]\times\N\times\R_+}\mathbbm{1}_{k\leq \tilde{Y}^K_{s-}\wedge(RK)}\mathbbm{1}_{y\leq \tilde{d}^K_2(\frac{\tilde{X}^K_{s-}}{K}, \frac{\tilde{Y}^K_{s-}}{K},s)
         }P^{Y,d}(ds,dk,dy),
  \end{align*}
  where
  \[
    \tilde{d}^K_i(x,y,t,\omega):=d_i(Z^K_s(\omega))+\hat{d}^K_i(t,\omega)+c_{11}(Z^K_s(\omega))(x\wedge
    R)+c_{12}(Z^K_s(\omega))(y\wedge R),
  \]
  satisfies $\tilde{X}^K_t=X^K_t$ and $\tilde{Y}^K_t=Y^K_t$ for all $t\leq e^{VK}$ with probability converging to 1. Hence,
  it is sufficient to prove~\eqref{eq:ccl-lemme-2D} for the process $(\tilde{X}^K,\tilde{Y}^K)$. The reason of this transformation is
  that the transition rate of the process $(\tilde{X}^K,\tilde{Y}^K)$ with jump size $(-1,0)$ (resp.\ $(0,-1)$) at time $t$ when
  $(\tilde{X}^K_t,\tilde{Y}^K_t)=K(x,y)$ takes the form $K\mu^K_1(x,y,t)$ (resp.\ $K\mu^K_2(x,y,t)$) where the functions
  $\mu^K_i(x_1,x_2,t,\omega):=(x_i\wedge R) \tilde{d}^K_i(x_1,x_2,t,\omega)$ for $i\in\{1,2\}$ are, for large $K$, bounded Lipschitz with respect to
  $(x_1,x_2)$, uniformly with respect to $K$, $t$ and $\omega$. Similarly, the transition rate of the process
  $(\tilde{X}^K,\tilde{Y}^K)$ with jump size $(1,0)$ (resp.\ $(0,1)$) at time $t$ when $(\tilde{X}^K_t,\tilde{Y}^K_t)=K(x,y)$ takes
  the form $K\lambda^K_1(x,y,t)$ (resp.\ $K\lambda^K_2(x,y,t)$) where the functions
  $\lambda^K_i(x_1,x_2,t,\omega):=(x_i\wedge R)b_i(Z^K_t(\omega))$ for $i\in\{1,2\}$ are for large $K$, bounded Lipschitz with respect to
  $(x_1,x_2)$, uniformly with respect to $K$, $t$ and $\omega$. Let $C_\text{Lip}$ (resp.\ $\bar{C}$) be a constant larger than the
  Lipschitz semi-norms (resp.\ $L^\infty$ norm) of all these functions and larger than the Lipschitz semi-norms (resp.\ $L^\infty$
  norm) of $b_i$, $d_i$ and $c_{ij}$ for all $i,j\in\{1,2\}$.
  
  

  Let $e_1= (1, 0)$, $e_2=(0,1)$. The previous discussion means the process $(\tilde{X}^K,\tilde{Y}^K)$ admits the representation: for all $t\geq 0$,
  %
  \begin{align*}
    (\tilde{X}^{K}_t,\tilde{Y}^K_t)= (\tilde{X}^{K}_0,\tilde{Y}^K_0) + \sum_{i=1, 2} & e_i Q^{b,K}_i \left(K \int_{0}^t
                     \lambda^K_i\left(\frac{\tilde{X}^K_s}{K},\frac{\tilde{Y}^K_s}{K},s\right)ds\right) \\ & -e_i Q^{d,K}_i \left(K \int_{0}^t
                     \mu^K_i\left(\frac{\tilde{X}^K_s}{K},\frac{\tilde{Y}^K_s}{K},s\right)ds\right), 
  \end{align*}
  where the stochastic processes $Q^{b,K}_1$, $Q^{d,K}_1$, $Q^{b,K}_2$ and $Q^{d,K}_2$ are independent, rate $1$ Poisson processes,
  independent from $(X^K_0,Y^K_0)$. This is the so-called method of ``random time change'' developed in~\cite{ethier-kurtz-86}. Note
  that, assuming that these processes are nondecreasing stochastic processes in $\N$ with jumps of size $+1$ and no common jumps, the
  above expression allows to represent uniquely the processes $Q^{b,K}_1$, $Q^{d,K}_1$, $Q^{b,K}_2$ and $Q^{d,K}_2$ in terms of the
  Poisson point measures $P^{X,b}$, $P^{X,d}$, $P^{Y,b}$ and $P^{Y,d}$ involved in~\eqref{eq:def-X-2D} and~\eqref{eq:def-Y-2D}.
 
  This implies that
  \begin{align*}
    \frac{1}{K}(\tilde{X}^{K}_t,\tilde{Y}^K_t)= & \frac{1}{K}(X^{K}_0,Y^K_0) + \sum_{i=1,2} e_i\int_0^t
                                  \left(\lambda^K_i\left(\frac{\tilde{X}^K_s}{K},\frac{\tilde{Y}^K_s}{K},s\right)-\mu^K_i\left(\frac{\tilde{X}^K_s}{K},\frac{\tilde{Y}^K_s}{K},s\right)\right)ds
    \\ & +\sum_{i=1, 2} \Bigg[ \frac{e_i}{K} \widetilde{Q}^{b,K}_i \left(K \int_{0}^t
                     \lambda^K_i\left(\frac{\tilde{X}^K_s}{K},\frac{\tilde{Y}^K_s}{K},s\right)ds\right) \\ & \qquad\quad -\frac{e_i}{K} \widetilde{Q}^{d,K}_i \left(K \int_{0}^t
                     \mu^K_i\left(\frac{\tilde{X}^K_s}{K},\frac{\tilde{Y}^K_s}{K},s\right)ds\right)\Bigg], 
  \end{align*}
  where $\widetilde{Q}^{b,K}_i(t):=Q^{b,K}_i(t)-t$ and $\widetilde{Q}^{d,K}_i(t):=Q^{d,K}_i(t)-t$ are the compensated Poisson
  processes. Therefore, by~\eqref{eq:bound-R-2D}, with probability converging to 1, for all $t>0$, there exists a constant $C>0$ depending only on $C_\text{Lip}$ and $R$ such that 
  \begin{multline*}
    \left|\frac{1}{K}(\tilde{X}^{K}_t,\tilde{Y}^K_t)-\phi_t\left(\frac{1}{K}X^K_0,\frac{1}{K}Y^K_0,Z^K_0\right)\right| \\
    \begin{aligned}
      & \leq C\int_0^t
      \left[ 
      \left|\frac{1}{K}(\tilde{X}^{K}_s,\tilde{Y}^K_s)-\phi_s\left(\frac{1}{K}X^K_0,\frac{1}{K}Y^K_0,Z^K_0\right)\right|
        +
        |Z^K_s-Z^K_0|+
        \left(\hat{d}_1^K(s)+\hat{d}_2^K(s)\right)\right]ds \\
      & +\sum_{i=1}^2 \sup_{v\in[0,t]}\left|\widetilde{Q}^{b,K}_i(K\bar{C}v)\right|+\sup_{v\in[0,t]}\left|\widetilde{Q}^{d,K}_i(K\bar{C}v)\right|.
    \end{aligned}
  \end{multline*}
  We then deduce from Gronwall's lemma that, for all $t>0$,
  \begin{multline*} 
    \sup_{s\in[0,t]}\left|\frac{1}{K}(\tilde{X}^{K}_s,\tilde{Y}^K_s)-\phi_s\left(\frac{1}{K}X^K_0,\frac{1}{K}Y^K_0,Z^K_0\right)\right| \\ \leq
    \left[C
    \int_0^t\left|\frac{d
          Z^K_s}{ds}\right|ds+C\sum_{i=1}^2\sup_{s\in[0,T]}\hat{d}^K_i(s)+\sup_{s\in[0,t]}\left|\widetilde{Q}^{b,K}_i(K\bar{C}s)\right|+\sup_{s\in[0,t]}\left|\widetilde{Q}^{d,K}_i(K\bar{C}s)\right|\right]e^{C
          t}.
  \end{multline*}
  Now, observe that, for all $v>0$, for all $\gamma>0$,
  \[
    \lim_{K\to+\infty}\P\left(\sup_{u\in[0,v]}\left|K^{-1}\widetilde{Q}(Ku)\right|>\gamma\right)=0,
  \]
  for $Q$ a compensated Poisson process (cf.\ e.g.~\cite[Eq.\ (2.5), Ch.\ 11]{ethier-kurtz-86}). Combining this
  with~\eqref{ass:vanishing_terms}, we deduce that, for all $t>0$ and all $\gamma>0$, with probability converging to 1
  \[
    \sup_{s\in[0,t]}\left|\frac{1}{K}(\tilde{X}^{K}_s,\tilde{Y}^K_s)-\phi_s\left(\frac{1}{K}X^K_0,\frac{1}{K}Y^K_0,Z^K_0\right)\right|\leq\gamma.
  \]
  This yields~\eqref{eq:step1-2D}.

  \bigskip

  \noindent\textit{Step 2: properties of the flows $\phi$.} Since $b_i$, $d_i$, $c_{ij}$ are $C^1_b$ for $i,j\in\{1,2\}$, it is
  standard to check that the flow $\phi_t(x,y,z)$ is continuous with respect to all variables $t$, $x$, $y$ and $z$.
  For all $z\in G$, since $b(z)-d(z)\geq\delta$ and by ~\eqref{eq:IIF1} and~\eqref{eq:IIF2}, it is well-known (cf.\
  e.g.~\cite{hofbauer-sigmund} or~\cite{champagnat-jabin-raoul-10}) that the 2-dimensional competitive Lotka-Volterra
  system~\eqref{eq:2D-LV} admits the equilibrium $(0,u^\star(z))$ as globally asymptotically stable equilibrium, i.e.
  \begin{equation}
    \label{eq:pf-penible-bis}
    \forall (x,y)\in(0,+\infty)^2,\ \lim_{t\to+\infty} \phi_t(x,y,z)=\left(0,u^\star(z)\right).    
  \end{equation}
  We need uniform stability estimates of this equilibrium with respect to $z\in G$. We do this using uniform local Lyapunov functions.
  For all $x,y\geq 0$ and $z\in G$, set
  \[
    L_z(x,y):=\alpha x^2+\left(y-u^\star(z)\right)^2\quad\text{with}\quad\alpha:=\frac{5\bar{C}^3}{12\delta\underline{c}^2}\vee 1.
  \]
  Then, for all $z\in G$ and all $x,y$ such that
  \begin{equation}
    \label{eq:pf-penible}
    \sqrt{x^2+(y-u^\star(z))^2}\leq \frac{\delta}{4\bar{C}},    
  \end{equation}
  defining $r_i(z)=b_i(z)-d_i(z)$, and setting $\phi_t(x,y,z)=(\phi_1,\phi_2)$
  \begin{align*}
    \frac{d}{dt}(L_z(\phi_t(x,y,z))) & =2\alpha
                                    \phi_1^2\left(r_1(z)-c_{12}(z)\phi_2-c_{11}(z)\phi_1\right)+2\phi_2\left(\phi_2-u^\star(z)\right)\left(r_2(z)-c_{22}(z)\phi_2-c_{21}(z)\phi_1\right) \\
    & \leq -\frac{3\alpha\delta}{2}\phi_1^2-2c_{22}(z)\phi_2\left(\phi_2-u^\star(z)\right)^2-2c_{21}(z)\phi_1\phi_2\left(\phi_2-u^\star(z)\right),
  \end{align*}
  where we used that $r_1(z)-c_{12}(z)\phi_2=g_{12}(z)-c_{12}(\phi_2-u^\star(z))\leq -\delta+\delta/4$ by~\eqref{eq:IIF2}
  and~\eqref{eq:pf-penible}. By Young's inequality,
  \[
    -2\phi_1\left(\phi_2-u^\star(z)\right)\leq\frac{c_{21}}{2c_{22}}\phi_1^2+2\frac{c_{22}}{c_{21}}\left(\phi_2-u^\star(z)\right),
  \]
  hence
  \begin{align*}
    \frac{d}{dt}(L_z(\phi_t(x,y,z))) & \leq
                                    -\left(\frac{\alpha\delta}{4}-\frac{c_{21}^2}{2c_{22}}\phi_2\right)\phi_1^2. 
  \end{align*}
  Since $u^\star(z)\geq\delta/\bar{C}$, by~\eqref{eq:pf-penible}, $y\leq \frac{5}{4}u^\star(z)\leq \frac{5\bar{C}}{4\underline{c}}$.
  Therefore, by definition of $\alpha$, we deduce that $\frac{d}{dt}(L_z(\phi_t(x,y,z)))\leq 0$ for all $x,y$
  satisfying~\eqref{eq:pf-penible}, i.e. $L_z(\phi_t(x,y,z))$ is non-increasing if $x,y$ satisfies~\eqref{eq:pf-penible}.

  Now, $|(x,y)-(0,u^\star(z))|^2\leq L_z(x,y)\leq \alpha |(x,y)-(0,u^\star(z))|^2$. Therefore, if for some $t\geq 0$,
  $|\phi_t(x,y,z)-(0,u^\star(z))|\leq \varepsilon$ for some $\varepsilon\leq \delta/(4\bar{C})$, then for all $s\geq t$,
  $|\phi_s(x,y,z)-(0,u^\star(z))|\leq \sqrt{\alpha}\varepsilon$. This implies that, for all $\varepsilon\leq \delta/(4\bar{C})$,
  $z\in G$ and $x,y> 0$,
  \[
    t^\text{last}_\varepsilon(x,y,z)\leq t^\text{first}_{\varepsilon/\sqrt{\alpha}}(x,y,z),
  \]
  where $t^\text{first}_\varepsilon(x,y,z)$ (resp.\ $t^\text{last}_\varepsilon(x,y,z)$) is the first (resp.\ last) time where
  $\phi_t(x,y,z)$ hits (resp.\ does not belong to) the ball of radius $\varepsilon$ centered at $(0,u^\star(z))$, {with the convention $t^\text{last}_\varepsilon(x,y,z)=0$ if $\phi_t(x,y,z)$ belongs to the ball for any $t\geq 0$}.
  By~\eqref{eq:pf-penible-bis}, both times $t^\text{first}_\varepsilon(x,y,z)$ and $t^\text{last}_\varepsilon(x,y,z)$ are finite.
        
  Let $\varepsilon\leq \delta/(4\bar{C})$. Let us prove that
  \[
    T_0(\varepsilon,G,H):=\sup_{z\in G,(x,y)\in H}t^\text{last}_{\varepsilon/2}(x,y,z)<+\infty.
  \]
  It is enough to prove it for $t^\text{last}_{\varepsilon/2}$ replaced by $t^\text{first}_{\varepsilon/2\sqrt{\alpha}}$. By contradiction
  and since $G$ and $H$ are compact, if it was not true, there would exist a sequence $(x_n,y_n,z_n)\in H\times G$ such that, as
  $n\to+\infty$, $(x_n,y_n,z_n)\to (x,y,z)\in H\times G$ and $t^\text{first}_{\varepsilon/2\sqrt{\alpha}}(x_n,y_n,z_n)\to+\infty$.
  Since $(x,y,z)\in H\times G$, $t^\text{first}_{\varepsilon/4\sqrt{\alpha}}(x,y,z)<+\infty$. By continuity of the flow $\phi$, we
  deduce that, for $n$ large enough, $t^\text{first}_{\varepsilon/2\sqrt{\alpha}}(x_n,y_n,z_n)\leq
  t^\text{first}_{\varepsilon/4\sqrt{\alpha}}(x,y,z)$. This is a contradiction.

  By continuity of the function $\phi$ again and since $\phi^1_t(x,y,z)>0$ for all $x>0$, $y\geq 0$, $z\in\R^\dimE$ and $t<+\infty$, where $\phi^1$ is the first coordinate of $\phi$, 
  \[
    \eta(\varepsilon,G,H):=\frac{1}{2}\inf_{z\in G,\ (x,y)\in H}\phi^1_{T_0(\varepsilon,G,H)}(x,y,z)>0.
  \]
  Combining Step 1 with $t=T_0(\varepsilon,G,H)$ with~\eqref{eq:hyp-CI-lemme-2D} and the previous definitions
  entails~\eqref{eq:ccl-lemme-2D} with the values of $T_0$ and $\eta$ defined above. Therefore, Lemma~\ref{lem:competition} is
  proved.
\end{proof}

\newcommand{\thetainv}{\theta_\text{inv}^K}
\section{Proof of Theorem~\ref{thm:BDP-max-2}}
\label{sec:proof-main}

\subsection{Reduction of the problem and strategy of proof}
\label{sec:reduct-pb}

Define $(t_k)_{k\geq 0}$, $(\ell^\star_k)_{k\geq 0}$, $S$ and $(\beta(t),e(t))_{t\in[0,S)}$ as in the statement of
Theorem~\ref{thm:BDP-max-2}. Given $T<S$, our goal is to prove the convergence of $(\beta^K(t), E^K(t \log (K)))_{t\in[0,T]}$ to
$(\beta(t),e(t))_{t\in[0,T]}$ in $L^\infty([0,T])$ in probability. It is enough to check this for an increasing sequence of times $T$
converging to $S$, that can be assumed without loss of generality different from all the times $t_k$. Hence fix $T\neq t_k$ for all
$k$ such that $t_k$ is well-defined. By definition of $S$, this means that there exists $k^\star<+\infty$ such that $t_k$ is
well-defined and $t_k< T$ for all $k\leq k^\star$, and $T<t_\ell$ for all $\ell>k^\star$ such that $t_\ell$ is well-defined.

The proof proceeds by induction on the time intervals $[t_k,t_{k+1}]$. However, if we consider the state of the process at a time $t_1\log K$, in the case where $t_1$ is a time of change of dominant population, then we have two exponents $\beta^K_i(t_1)$ converging to 1 when $K$ goes to $+\infty$. Similarly, if $t_1$ is
the time of extinction of species $j$, we may not be sure that $N^K_j(t_1\log K)=0$ with probability converging to 1. Hence, in both
cases, the process at time $t_1\log K$ may not satisfy the assumptions about the initial condition in Theorem~\ref{thm:BDP-max-2}.
Therefore, to make the induction work, we need to prove convergence on a time interval slightly longer than $[0,t_1]$, as explained
below.

\bigskip

It is actually sufficient to prove that there exists $0<\gamma<t_2-t_1$ (with the convention $t_2-t_1=\infty$ if $t_1=\infty$ or
$t_2$ is not defined) such that we have the convergence in probability
\begin{equation}
  \label{eq.cv.proba.iter1}
  \sup_{t\in[0,T\wedge (t_1+\gamma)]}|\beta^K(t)-\beta(t)|+|E^K(t \log (K))-e(t)|\xrightarrow[K\to+\infty]{}0
\end{equation}
and to prove what follows.
\begin{itemize}
\item {\bf Case 1.} If $T>t_1+\gamma$ and if there exists a single $j$ such that $\beta_{j}(t_1)=0$,
  $f_{j,\ell_0^\star}(e(t_1))<0$ and $\beta_\ell(t_1)\in(0,1)$ for all $\ell\neq\ell^\star_0,j$, then, in probability
  \begin{equation}
    \label{eq:cv.proba.t1plusgamma.extinction}
    \frac{1}{K}N_{\ell_0^\star}((t_1+\gamma)\log K)\xrightarrow[K\to+\infty]{}n^\star_{\ell_0^\star}(e(t_1+\gamma))
  \end{equation}
  and for all $t>t_1$
  \begin{equation}
    \label{eq:extinction}
    \lim_{K\to\infty}\P\left(N_j^K(t\log K)=0\right)=1.
  \end{equation}
\item {\bf Case 2.} If $T>t_1+\gamma$ and if there exists a single $j\neq \ell^\star_0$ such that $\beta_{j}(t_1)=1$,
  $f_{j,\ell_0^\star}(e(t_1))>0$ and $\beta_\ell(t_1)\in(0,1)$ for all $\ell\neq\ell^\star_0,j$, then, in probability
  \begin{equation}
    \label{eq:cv.proba.t1plusgamma.invasion}
    \frac{1}{K}N_{j}((t_1+\gamma)\log K)\xrightarrow[K\to+\infty]{}n^\star_{j}(e(t_1+\gamma)).
  \end{equation}
\end{itemize}

Indeed, if $T\leq t_1+\gamma$, Theorem~\ref{thm:BDP-max-2} is already proved from~\eqref{eq.cv.proba.iter1}. If $T>t_1+\gamma$, since
$T<S$, this means that induction did not stop at time $t_1$, hence one and only one of the two cases above occurs. In Case 1, since $0<\gamma<t_2-t_1$, $\beta_{i}(t_1+\gamma)\in(0,1)$ for all $i\neq j,\ell_0^\star$ and by~\eqref{eq:extinction}, with
probability converging to 1, $N^K_j(t_1+\gamma)=0$. Hence, reducing the dimension of the process to remove the extinct population
$j$, by \eqref{eq:cv.proba.t1plusgamma.extinction}, the process at time $t_1+\gamma$ satisfies the assumptions on the initial
condition in Theorem~\ref{thm:BDP-max-2}. In Case 2, $j=\ell^\star_1$ and by Assumption~\ref{hyp:Lotka-Volterra-max-2},
$\beta_{i}(t_1+\gamma)\in(0,1)$ for all $i\neq \ell_1^\star$. Hence, by~\eqref{eq:cv.proba.t1plusgamma.invasion}, the process at time
$t_1+\gamma$ satisfies the assumptions on the initial condition in Theorem~\ref{thm:BDP-max-2}.


Therefore, from the Markov property and by induction we deduce from \eqref{eq.cv.proba.iter1} that, for all
$k\in\{1,\ldots, k^\star\}$, there exists $0<\gamma_{k+1}<t_{k+1}-t_k$ (with the convention $t_{k+1}-t_k=\infty$ if $t_k=\infty$ or
$t_{k+1}$ is not defined) such that, in probability,
\begin{equation*}
\sup_{t\in[T\wedge (t_k+\gamma_k),T\wedge (t_{k+1}+\gamma_{k+1})]}|\beta^K(t)-\beta(t)|+|E^K(t \log (K))-e(t)|\xrightarrow[K\to+\infty]{}0,
\end{equation*}
with the convention $\gamma_1=\gamma$. Hence, by definition of $k^\star$,
\begin{equation*}
\sup_{t\in[0,T]}|\beta^K(t)-\beta(t)|+|E^K(t \log (K))-e(t)|\xrightarrow[K\to+\infty]{}0
\end{equation*}
in probability, which concludes the proof.

\bigskip

Therefore, let us fix throughout the whole proof $T\in(0,S)$ different from all the $t_k$, and let us prove~\eqref{eq.cv.proba.iter1}
for some $0<\gamma<t_2-t_1$, \eqref{eq:cv.proba.t1plusgamma.extinction} and~\eqref{eq:extinction} in Case 1,
and~\eqref{eq:cv.proba.t1plusgamma.invasion} in Case 2.

For this purpose, we introduce, for all $t\geq 0$
\[
  \tilde{e}(t)=e(0)+\int_0^t F_{\ell^\star_0}(\tilde{e}(s))ds.
\]
and for all $i\neq \ell_0^\star$ and all $t\geq 0$
\begin{equation}
  \label{eq:def-tilde-beta_i}
  \tilde \beta_i(t):=\beta_i(0)+\int_0^{t\wedge \tilde{T}_i} f_{i,\ell_0^\star}(\tilde{e}(s))ds,  
\end{equation}
where
\begin{equation}
  \label{eq:def-tilde-T_i}
  \tilde{T}_i=\inf\left\{ t\geq 0: \beta_i(0)+\int_0^{t} f_{i,\ell_0^\star}(\tilde{e}(s))ds=0\right\}.
\end{equation}
Remark that $\tilde{e}(t)=e(t)$ and $\tilde \beta_i(t)=\beta_i(t)$ for all $t\in[0,t_1]$ but that these functions may differ after
time $t_1=\inf\{t\geq 0:\: \exists i\neq \ell_0^\star \ \text{ s.t. } \beta_i(t)\in\{0,1\}\}$ ($\beta_i$ and $e$ can even be undefined if
$t_1=S$).

The proof proceeds along the following steps: we fix $\varepsilon>0$ that we will make converge to 0 at the very end of the proof and
introduce $\thetainv$ the first time where the non-dominant population have total size larger than $\varepsilon K$, i.e.
\[
    \thetainv:=T\wedge\inf\left\{t>0,\, \frac{1}{K}\sum_{i\neq\ell_0^\star}N_i^K(t\log K)> \eps \right\}.
\]
Below, we will say that two quantities are $\varepsilon$-close whenever their distance is smaller than a constant independent of $\varepsilon$ times $\varepsilon$ with probability converging to 1 as $K\to\infty$. All the following steps will be proved for a fixed value of $\varepsilon>0$, provided it is small enough.
\begin{enumerate}
\item We first prove that, before $\thetainv$, the dominant species $\ell^\star_0$ remains $\varepsilon$-close to its equilibrium
  (Subsection~\ref{sec.dominant}) and the environment remains $\varepsilon$-close to $\tilde{e}$ (Subsection~\ref{sec.environment}).
  For this, we use Lemma~\ref{lem:time-to-equi}.
\item We then prove in Subsection~\ref{sec.stoch.domination} that the exponent of any non-dominant species $j$ remains
  $\varepsilon$-close to $\tilde{\beta}_j$ until time $\thetainv\wedge\tilde{t}_2$ for some $\tilde{t}_2>t_1$ using comparison
  techniques to upper and lower bound the exponents (in Subsection~\ref{sec.stoch.domination}) and applying
  Lemma~\ref{lem:time-inhomo-LBD} to the bounding processes. To ensure the convergence of the exponents, we need that the two
  bounding processes have a similar behaviour. In view of Lemma~\ref{lem:time-inhomo-LBD}~(ii), this requires that the two bounding
  processes are far from the unsettled case where $b(T_0)= d(T_0)$ (with the notations of Lemma~\ref{lem:time-inhomo-LBD}). This
  explains why we choose
  \begin{equation}
    \label{eq:def-tilde-t2}
    t_1< \tilde t_2 < \inf_{\substack{i \neq \ell_0^\star \text{ s.t. }\\ \beta_i(t_1)\neq 0}} \inf \left\{t>0,\tilde \beta_i(t)=0\right\}    
  \end{equation}
  to avoid having to consider extinctions of other species than those going extinct at time $t_1$. We can now define the parameter
  $\gamma$ in~\eqref{eq.cv.proba.iter1} as
  \begin{equation}
    \label{eq:def-gamma}
    \gamma:=1\wedge\frac{t_2-t_1}{2}\wedge(\tilde{t}_2-t_1),    
  \end{equation}
  so that $t_1+\gamma\leq \tilde{t}_2$ and $t_1+\gamma<t_2$.
\item After these first two steps, the proof is easy to complete in the following cases:
  \begin{itemize}
  \item if $T<t_1$, then all the $\beta^K_j(t)$ for $j\neq\ell^\star_0$ are close to $\tilde{\beta}_j(t)=\beta_j(t)<1$ up to time
    $\thetainv\wedge t_1$, hence at time $\thetainv\wedge t_1$ the total non-dominant population size is $o(K)$, which means that
    $\thetainv=T$, so~\eqref{eq.cv.proba.iter1} is proved (see Subsection~\ref{sec:ccl-T<t_1});
  \item if $T>t_1$ and there exists $j$ such that $\beta_j(t_1)=0$, then by definition of $\gamma$, up to time
    $\thetainv\wedge (t_1+\gamma)$, $\beta^K_j(t)$ is close to $\tilde{\beta}_j(t)<1$ and all the $\beta^K_\ell(t)$ for
    $\ell\neq\ell^\star_0,j$ are close to $\tilde{\beta}_\ell(t)\in(0,1)$. As above, this implies that $\thetainv=T$ if
    $T<t_1+\gamma$ or $\thetainv>t_1+\gamma$ if $T>t_1+\gamma$. Hence $\tilde{\beta}=\beta$ on $[0,(t_1+\gamma)\wedge T]$ since
    $\ell^\star_1=\ell^\star_0$ and $t_1+\gamma<\tilde t_2$, so~\eqref{eq.cv.proba.iter1} is proved (see Subsection~\ref{sec:ccl-case-1}):.
  \end{itemize}
\item In the remaining case where $T>t_1$ and there exists $j\neq\ell^\star_0$ such that $\beta_j(t_1)=1$, assuming
  $\thetainv<t_1-L\varepsilon$ for some explicit $L$ leads to a contradiction because all the non-dominant exponents are
  $\varepsilon$-close to their limit which is smaller than 1. Similarly, assuming $\thetainv>t_1+L\varepsilon$ would lead to a
  contradiction because at least one non-dominant exponent would become larger that 1. Therefore, $\thetainv$ is $\varepsilon$-close
  to $t_1$ (in particular, $\thetainv<T$). We then use Lemma~\ref{lem:competition} to control the change of dominant species after
  time $\thetainv$. Finally, we proceed similarly as in Steps 1 and 2 above to control the exponents after the change of dominant
  population up to time $t_1+\gamma$ (see Subsection~\ref{sec:ccl-case-2}).
\end{enumerate}

\medskip

We conclude this section with a useful estimate: since $|E^K(0)-e_0|<1$ with probability converging to 1, by~\eqref{eq:environment},
for all $t\in[0,T\log K]$, $|E^K(t)|\leq |e_0|+1+\|F\|_\infty T$ with probability converging to 1. Since by
Assumption~\ref{hyp:regularity-Ltka-Volterra}, for all $i\in\{1,\ldots,d\}$, $b_i$ and $d_i$ are continuous,
$b_i(\cdot)>d_i(\cdot)$ and $c_i(\cdot,i)\geq\underline{c}>0$, there exists a constant $\underline{n}^\star>0$ such that, with probability converging to 1,
\begin{equation}
  \label{eq:LB-nstar}
  \forall i\in\{1,\ldots,d\},\ \forall t\in[0,T\log K],\quad n^\star_i(E^K(t))\geq \underline{n}^\star.  
\end{equation}

Recall also the notations for the different norms of the parameters of the model introduced in the beginning of
Section~\ref{sec:prelim}.

\subsection{The size of the dominant species $\ell_0^\star$ stays close to $n^\star_{\ell_0^\star}$ before $\thetainv$}
\label{sec.dominant}

In this section, we prove that, with probability converging to 1 as $K\to+\infty$, the dominant species stays close to the
equilibrium associated to the environmental process $E^K$ before time $\thetainv$, i.e.
\begin{equation}
\label{eq:step-1-initial}
    \sup_{t\in[0,\thetainv]}\left|\frac{N_{\ell_0^\star}^K(t\log K)}{K}-n_{\ell_0^\star}^\star(E^K(t\log K))\right|\leq 4\varepsilon \left(1+\frac{\bar{c}}{\underline{c}}\right).
\end{equation}
To do so, we apply Lemma~\ref{lem:time-to-equi} using two carefully chosen logistic time-inhomogeneous birth-death processes
$(Y^K(t))_{t\geq 0}$ and $(Z^K(t))_{t\geq 0}$ for which the following stochastic domination holds almost surely:
\begin{equation}
\label{eq.stoch.dom.LogBD}
\forall t\in [0,\thetainv\log K],\:Y^K(t)\leq N_{\ell_0^\star}^K(t)\leq Z^K(t).
\end{equation}

With this in mind, let us control the competition term of the process $N_{\ell_0^\star}^K$. 
By the definition of $\thetainv$, for all $t< \thetainv$,  we have $\frac{1}{K}\sum_{i\neq\ell_0^\star}N_i^K(t\log K)\leq
\varepsilon$. 
Hence, almost surely, for all $t< \thetainv\log K$
\[
 \frac{1}{K}c_{\ell_0^\star}(E^K(t),\ell_0^\star)N^K_{\ell_0^\star}(t) 
 \leq \frac{1}{K}\sum_{j=1}^d c_{\ell_0^\star}(E^K(t),j)N^K_j(t)
 \leq \frac{1}{K}c_{\ell_0^\star}(E^K(t),\ell_0^\star)N^K_{\ell_0^\star}(t) + \bar{c}\varepsilon
\]
Then the two processes $Y^K$ and $Z^K$ defined by
\begin{align*}
    Y^K(t) & = N^K_{\ell^\star_0}(0)+\iiint_{[0,t]\times\N\times\R_+}\mathbbm{1}_{k\leq Y^K(s-)}\mathbbm{1}_{y\leq b_{\ell^\star_0} (E^K(s))}P^b_{\ell^\star_0} (ds,dk,dy)
  \\ & -\iiint_{[0,t]\times\N\times\R_+}\mathbbm{1}_{k\leq Y^K(s-)}\mathbbm{1}_{y\leq d_{\ell^\star_0} (E^K(s))+\bar{c}\varepsilon+\frac{1}{K}
    c_{\ell^\star_0} (E^K(s), {\ell^\star_0})Y^K(s-)}P^d_{\ell^\star_0} (ds,dk,dy).  
\end{align*}
and
\begin{align*}
  Z^K(t) & = N^K_{\ell^\star_0}(0)+\iiint_{[0,t]\times\N\times\R_+}\mathbbm{1}_{k\leq Z^K(s-)}\mathbbm{1}_{y\leq b_{\ell^\star_0} (E^K(s))}P^b_{\ell^\star_0} (ds,dk,dy)
  \\ & -\iiint_{[0,t]\times\N\times\R_+}\mathbbm{1}_{k\leq Z^K(s-)}\mathbbm{1}_{y\leq d_{\ell^\star_0} (E^K(s))+\frac{1}{K}
       c_{\ell^\star_0} (E^K(s), {\ell^\star_0})Z^K(s-)}P^d_{\ell^\star_0} (ds,dk,dy)
\end{align*}
satisfy the stochastic domination~\eqref{eq.stoch.dom.LogBD} (this is elementary to prove by induction on the successive jumps of the
three processes). Note that, using the notation of Lemma~\ref{lem:time-to-equi},
\begin{align*}
Y^K&\sim LogBD_K(b_{\ell_0^\star},d_{\ell_0^\star}+\bar c \eps,c_{\ell_0^\star}(\cdot,\ell_0^\star), E^K),\\
Z^K&\sim LogBD_K(b_{\ell_0^\star},d_{\ell_0^\star},c_{\ell_0^\star}(\cdot,\ell_0^\star), E^K).
\end{align*}

By assumption $N^K_{\ell_0^\star}(0)/K$, $\beta^K(0)$ and $E^K(0)$ converge in probability towards $n_{\ell_0^\star}^\star(e_0)$, $\beta(0)$ and $e(0)$ respectively, with $n_{\ell_0^\star}^\star$ continuous and with $\beta_i(0)<1$ for all $i\neq \ell_0^\star$. Then with probability converging to 1 as $K\to+\infty$,
\begin{align*}
  \left|\frac{N^K_{\ell_0^\star}(0)}{K}-n_{\ell_0^\star}^\star\left(E^K(0)\right)\right|\leq {\varepsilon}\,
  \qquad \text{and} \qquad  \frac{1}{K}\sum_{i\neq\ell_0^\star}N_i^K(0)<\varepsilon,
\end{align*}
in particular $\thetainv>0$ with probability converging to 1 and, by~\eqref{eq:LB-nstar}, taking $\varepsilon$ small enough, with
probability converging to 1, $n^\star_Y(E^K(t))\geq \underline{n}^\star/2$ for all $t\leq T\log K$, where
$n_Y^\star(\cdot)=n^\star_{\ell_0^\star}(\cdot)-\frac{\bar c \eps}{c_{\ell_0^\star}(\cdot,\ell_0^\star)}$. Then,
with probability converging to 1,
\[
  \left|\frac{Y^K(0)}{K}-n_Y^\star\left(E^K(0)\right)\right|
  \leq \left(1+\dfrac{\bar c}{\underline{c}}\right)\eps\,
  \quad \text{and} \quad  \left|\frac{Z^K(0)}{K}-n_{\ell_0^\star}^\star\left(E^K(0)\right)\right|\leq {\varepsilon}.
\]
Therefore, since $|\dot E^K|\leq \|F\|_\infty/\log K$, Lemma~\ref{lem:time-to-equi} entails that
\begin{multline*}
\mathbb P\left(\exists t\in[0,T],\,\, \frac{Y^K(t\log K)}{K}< n_{\ell^\star_0}(E^K(t\log K))-4\left(1+\dfrac{\bar
  c}{\underline{c}}\right)\eps \right) \\ \leq \mathbb P\left(\exists t\in[0,T],\,\, \frac{Y^K(t\log K)}{K}< n_Y^\star(E^K(t\log K))-3\left(1+\dfrac{\bar c}{\underline{c}}\right)\eps \right) \xrightarrow[]{K\to \infty} 0.
\end{multline*}
and
\[
\mathbb P\left(\exists t\in[0,T],\,\,\frac{Z^K(t\log K)}{K}> n_{\ell^\star_0}^\star(E^K(t\log K))+3\eps\right) \xrightarrow[]{K\to \infty} 0.
\]
Hence,~\eqref{eq:step-1-initial} follows from the stochastic domination~\eqref{eq.stoch.dom.LogBD}.

\medskip

Let us define the event
\[
  \Gamma^K_0:=\left\{\forall t\in [0,\thetainv],\ \left|\dfrac{N^K_{\ell_0^\star}(t\log K)}{K} - n_{\ell_0^\star}^\star(E^K(t\log K))\right|\leq 4\left(1+\dfrac{\bar{c}}{\underline{c}}\right)\eps\right\},
\]
whose probability converge to 1.

\subsection{The environment $E^K$ stays close to $\tilde{e}$ before $\thetainv$}
\label{sec.environment}

We prove now that, with probability converging to 1 as $K\to+\infty$, for $\eps$ sufficiently small,
\begin{equation}
\label{eq:step-2-initial}
    \sup_{t\in[0,\thetainv]}\left|E^K(t\log K)-\tilde{e}(t)\right|< A\eps
\end{equation}
with 
\begin{equation*}
A:=\sqrt{1+\left(5+4\frac{\bar c}{\underline c}\right)^2\left(3+2\|n^\star\|_\text{Lip}\right)^{-1}}
\exp \left(\|F\|_\text{Lip}\left(\frac{3}{2}+\|n^\star\|_\text{Lip}\right) T\right).
\end{equation*}

\medskip

To this aim, we define $\psi^K(t):=E^K(t\log K) - \tilde{e}(t)$, and we notice that, denoting by $\langle\cdot,\cdot\rangle$ the
Euclidean scalar product of $\R^\dimE$ and reminding that $|\cdot|$ denotes the associated norm, almost surely on the event
$\Gamma^K_0$, for all $t<\thetainv$,
 \begin{align*}
     \frac{d}{dt}|\psi^K(t)|^2 
     &= 2\left\langle\psi^K(t),F(\nu^K(t\log K),E^K(t\log K))-
            F(n^\star_{\ell_0^\star}(\tilde{e}(t))\delta_{\ell_0^\star}, \tilde{e}(t))\right\rangle\\
    & \leq 2|\psi^K(t)|\,\|F\|_{\text{Lip}}\,\eps
    \\
     &\qquad   + 2|\psi^K(t)|\,\left| F\left(\frac{N^K_{\ell_0^\star}(t\log K)}{K}\delta_{\ell_0^\star}, E^K(t\log K)\right) - F\left(n^\star_{\ell_0^\star}(\tilde{e}(t))\delta_{\ell_0^\star}, \tilde{e}(t))\right)\right|\\
    & \leq 2|\psi^K(t)|\,\|F\|_{\text{Lip}} \left(5 + 4\frac{\bar c}{\underline c}\right)\eps
        + 2\|F\|_{\text{Lip}}\left(1+\|n^\star\|_\text{Lip}\right) |\psi^K(t)|^2
        \\
    & \leq \|F\|_{\text{Lip}} \left(5 + 4\frac{\bar c}{\underline c}\right)^2\eps^2
        + \|F\|_{\text{Lip}} \left(3+2\|n^\star\|_\text{Lip}\right) |\psi^K(t)|^2,   
 \end{align*}
 using Young's inequality in the last line. Gr\"onwall's lemma then entails that, almost surely on $\Gamma^K_0$, for all
 $t\leq\thetainv$,
 \[
   \psi^K(t)^2  \leq \left(\psi^K(0)^2+\left(5+4\frac{\bar c}{\underline c}\right)^2\,\eps^2\left(3+2\|n^\star\|_\text{Lip}\right)^{-1}\right)\exp\left(\|F\|_{\text{Lip}}\left(3+2\|n^\star\|_\text{Lip}\right) t\right).
 \]
Since $E^K(0)$ converges in probability towards $e(0)$, $\psi^K(0)$ converges in probability to 0, we deduce~\eqref{eq:step-2-initial} for $K$ large enough.

\medskip

Combining this with~\eqref{eq:step-1-initial}, the event
\begin{equation}
  \label{eq:Gamma-ell}
  \Gamma^K:=\left\{\forall t\in [0,\thetainv],\ \left|E^K(t\log K)-\tilde{e}(t)\right|\leq A\eps\text{ and }\left|\frac{N^K_{\ell_0^\star}(t\log K)}{K} - n_{\ell_0^\star}^\star(\tilde{e}(t))\right|\leq B\eps\right\},
\end{equation}
where
\[
  B:=4\left(1+\dfrac{\bar{c}}{\underline{c}}\right)+A\|n^\star\|_\text{Lip},
\]
satisfies
\begin{equation}
  \label{eq:proba-Gamma->1}
  \lim_{K\to+\infty}\P(\Gamma^K)=1.
\end{equation}
Note also that, on the event $\Gamma^K$, for all $t\leq\thetainv$,
\[
  1+\frac{\log\left(n_{\ell_0^\star}^\star(\tilde{e}(t))-B\eps\right)}{\log K}\leq\beta^K_{\ell^\star_0}(t)\leq 1+\frac{\log\left(K^{-1}+n_{\ell_0^\star}^\star(\tilde{e}(t))+B\eps\right)}{\log K}
\]
and therefore, 
\begin{equation}
  \label{eq:cv-beta-ell_0}
  \sup_{t\in[0,\thetainv]}\left|\beta^K_{\ell^\star_0}(t)-1\right|\xrightarrow[K\to+\infty]{} 0
\end{equation}
in probability.

\subsection{Stochastic domination of the non-dominant species and bounds on their exponents}
\label{sec.stoch.domination}

Recall the definitions~\eqref{eq:def-tilde-T_i} and~\eqref{eq:def-tilde-beta_i} of $\tilde{T}_i$ and $\tilde{\beta}_i$ for all
$i\neq\ell^\star_0$. In this step, our goal is to prove that, for any $i\neq \ell_0^\star$ and any $s_0\in[0,T]$ such that either
$s_0<\tilde{T}_i$ (or, equivalently, $\tilde{\beta}_i(s_0)>0$), or $s_0>\tilde{T}_i$ (which implies $\tilde{\beta}_i(s_0)<0$) and
$f_{i,\ell^\star_0}(\tilde{e}(\tilde{T}_i))<0$, then for $\eps$ small enough, we have with probability converging to 1 as $K\to+\infty$
\begin{equation}
  \label{eq:step-exponents}
  \lim_{K\to+\infty}\P\left(\forall t\in[0,s_0\wedge\thetainv],\ \left|\beta^K_i(t)-\tilde{\beta}_i(t)\right|\leq C\eps\right)=1
\end{equation}
for some explicit constant $C$. Note that, if $\tilde{T}_i<+\infty$, the left derivative of $\tilde{\beta}_i$ at time $\tilde{T}_i$,
given by $f_{i,\ell^\star_0}(\tilde{e}(\tilde{T}_i))$, is necessarily nonpositive. Therefore any $s_0\leq T$ would satisfy the above
condition except in the case $f_{i,\ell^\star_0}(\tilde{e}(\tilde{T}_i))=0$ where we must take $s_0<\tilde{T}_i$.

To do so, we first prove with probability converging to 1 a stochastic domination property on the time interval $[0,\thetainv]$
between the process $N^K_i$ and two linear time inhomogeneous birth-death processes, next we apply Lemma~\ref{lem:time-inhomo-LBD} to
prove the convergence of the exponents of these two bounding processes and finally we use Lemma~\ref{lem_control_reaching_time_of_0}
to prove that the limits are $\varepsilon$-close.

\bigskip

Let $i\neq\ell^\star_0$. By Assumption~\ref{hyp:regularity-Ltka-Volterra}, almost surely on the event $\Gamma^K$, for all $t\in[0,T]$,
$b^Y(t)\leq b_i(E^K(t\log K))\leq b^Z(t)$, where
\begin{equation}
  \label{eq:blabla-def-b}
  b^Y(t):=b_i(\tilde{e}(t)) - A\|b\|_{\text{Lip}}\eps\quad\text{and}\quad  b^Z(t):=b_i(\tilde{e}(t)) + A\|b\|_{\text{Lip}}\eps.  
\end{equation}
Since $b_i(x)\geq c_i(x,i)n^\star_i(x)\geq\underline{c}\underline{n}^\star$, where $\underline{n}^\star$ is defined
in~\eqref{eq:LB-nstar}, we deduce that $b^Y(t)>0$ for all $t\leq T$ provided $\eps$ is small enough. Similarly, on the event
$\Gamma^K$, for all $t\in[0,T]$,
\[
d_i(\tilde{e}(t)) - A\|d\|_{\text{Lip}} \eps \leq d_i(E^K(t\log K))\leq d_i(\tilde{e}(t)) + A\|d\|_{\text{Lip}}\eps,
\]
and
\[
c_i(\tilde{e}(t),\ell_0^\star) - A\|c\|_{\text{Lip}}\eps \leq c_i(E^K(t\log K),\ell_0^\star) \leq c_i(\tilde{e}(t),\ell_0^\star) + A\|c\|_{\text{Lip}}\eps.
\]
Therefore, using that $\frac{1}{K}\sum_{j\neq\ell^\star_0}N^K_j(t\log K)\leq\varepsilon$ for $t<\thetainv$, the death
rate at time $t\log K$ of an individual of species $i$ satisfies for $\eps$ small enough, for $K$ large enough and almost surely on
the event $\Gamma^K$
\[
\forall t<\thetainv,\quad d^Z(t)\leq d_i(E^K(t\log K))+\frac{1}{K}\sum_{j=1}^d  c_i(E^K(t\log K),j) N_j^K(t\log K) \leq d^Y(t),
\]
where
\begin{equation}
  \label{eq:blabla-def-d-Y}
  d^Y(t):=d_i(\tilde{e}(t)) + A\|d\|_{\text{Lip}}\eps + \left(c_i(\tilde{e}(t),\ell_0^\star\right) + A\|c\|_{\text{Lip}}\eps)\left(n_{\ell_0^\star}^\star(\tilde{e}(t))+B\eps\right) +\bar c\eps
\end{equation}
and
\begin{equation}
  \label{eq:blabla-def-d-Z}
  d^Z(t):=d_i(\tilde{e}(t)) - A\|d\|_{\text{Lip}}\eps + \left(c_i(\tilde{e}(t),\ell_0^\star\right) - A\|c\|_{\text{Lip}}\eps)\left(n_{\ell_0^\star}^\star(\tilde{e}(t))-B\eps\right).
\end{equation}
Note that we used here that, for $\eps>0$ small enough, $c_i(\tilde{e}(t),\ell_0^\star)\geq\underline{c} > A\|c\|_{\text{Lip}}\eps$
and $n_{\ell_0^\star}^\star(\tilde{e}(t))\geq\underline{n}^\star>B\eps$. Reducing once more $\eps$, we also have that $d^Y(t)>0$ for
all $t\in[0,T]$.

Fix a value of $\eps$ as above. In view of the previous inequalities, the two processes $Y^K$ and $Z^K$ defined by
\begin{align*}
    Y^K(t) & = \lfloor K^{\beta_i(0)-\eps}\rfloor +\iiint_{[0,t]\times\N\times\R_+}\mathbbm{1}_{k\leq Y^K(s-)}\mathbbm{1}_{y\leq
             b^Y(s/\log K)}P^b_i (ds,dk,dy) 
  \\ & -\iiint_{[0,t]\times\N\times\R_+}\mathbbm{1}_{k\leq Y^K(s-)}\mathbbm{1}_{y\leq d^Y(s/\log K)}P^d_i (ds,dk,dy). 
\end{align*}
and
\begin{align*}
    Z^K(t) & = \lfloor K^{\beta_i(0)+\eps}\rfloor+\iiint_{[0,t]\times\N\times\R_+}\mathbbm{1}_{k\leq Z^K(s-)}\mathbbm{1}_{y\leq
             b^Z(s/\log K)}P^b_i (ds,dk,dy) 
  \\ & -\iiint_{[0,t]\times\N\times\R_+}\mathbbm{1}_{k\leq Z^K(s-)}\mathbbm{1}_{y\leq d^Z(s/\log K)}P^d_i (ds,dk,dy) 
\end{align*}
satisfy for $K$ large enough, almost surely on the event $\Gamma^K\cap\{\lfloor K^{\beta_i(0)-\eps}\rfloor\leq N^K_i(0)\leq \lfloor
K^{\beta_i(0)+\eps}\rfloor\}$,
\begin{align}
  \label{eq.stoch.dom.LinBD}
  \forall t\in [0,\thetainv\log K],\:Y_i^K(t) \leq N_i^K(t)\leq Z_i^K(t).
\end{align}
Note that, since we assumed that $\beta^K_i(0)$ converges in probability to $\beta_i(0)>0$,
\[
  \lim_{K\to+\infty}\P\left(\Gamma^K\cap\left\{\lfloor K^{\beta_i(0)-\eps}\rfloor\leq N^K_i(0)\leq \lfloor
      K^{\beta_i(0)+\eps}\rfloor\right\}\right)=1.
\]
Note also that the processes $Y^K$ and $Z^K$ have laws $LBD(b^{Y},d^{Y}, K^{\beta_i(0)-\eps})$ and
$LBD(b^{Z},d^{Z}, K^{\beta_i(0)+\eps})$, respectively.

\bigskip

Take $s_0\in[0,T]$ such that either $s_0<\tilde{T}_i$, or $s_0>\tilde{T}_i$ and $f_{i,\ell^\star_0}(\tilde{e}(\tilde{T}_i))<0$ and
let us prove that Lemma~\ref{lem:time-inhomo-LBD} can be applied to $Y^K$ and $Z^K$ on $[0,s_0]$. Define
\[
  D^-=A\left(\|b\|_\text{Lip}+\|d\|_\text{Lip}+\bar{n}^\star\|c\|_\text{Lip}\right)+AB\|c\|_\text{Lip}+(B+1) \bar{c}
\]
and
\[
  D^+=A\left(\|b\|_\text{Lip}+\|d\|_\text{Lip}+\bar{n}^\star\|c\|_\text{Lip}\right)+AB\|c\|_\text{Lip}+B\bar{c}.
\]
Assuming $\eps<1$, we have for all $t\in[0,T]$
\begin{equation}
  \label{eq:control-r_Y-r_Z}
  b^Y(t)-d^Y(t)\geq f_{i,\ell^\star_0}(\tilde{e}(t))-D^-\eps\quad \text{and}\quad b^Z(t)-d^Z(t)\leq f_{i,\ell^\star_0}(\tilde{e}(t))+D^+\eps,  
\end{equation}
where we recall that $f_{i,\ell^\star_0}(x)=b_i(x)-d_i(x)-c_i(x,\ell^\star_0)n^\star_{\ell^\star_0}(x)$.

In the case where $s_0<\tilde{T}_i$, $\tilde\beta_i(t)=\beta_i(0)+\int_0^t f_{i,\ell^\star_0}(\tilde e(s))ds>0$ for all $t\leq s_0$, hence,
by~\eqref{eq:control-r_Y-r_Z}, assuming $\eps$ small enough, $\beta_i(0)-\eps+\int_0^t(b^Y(s)-d^Y(s))ds>0$ for all $t\in[0,s_0]$.
Since obviously we also have $\beta_i(0)+\eps+\int_0^t(b^Z(s)-d^Z(s))ds>0$ for all $t\in[0,s_0]$, we can apply
Lemma~\ref{lem:time-inhomo-LBD}~(i) to $Y^K$ and $Z^K$ on $[0,s_0]$.

In the case where $s_0>\tilde{T}_i$, defining
\[
  T_\eps^Y:=\inf\left\{t\geq 0: \beta_i(0)-\eps+\int_0^t(b^Y(s)-d^Y(s))ds=0\right\}
\]
and
\[
  T_\eps^Z:=\inf\left\{t\geq 0: \beta_i(0)+\eps+\int_0^t(b^Z(s)-d^Z(s))ds=0\right\},
\]
we deduce from~\eqref{eq:control-r_Y-r_Z} that $T_\eps^{-}\leq T_\eps^Y\leq T_\eps^Z\leq T_\eps^{+}$ where
\[
  T_\eps^\pm:=\inf\left\{t\geq 0: \beta_i(0)\pm\eps+\int_0^t(f_{i,\ell^\star_0}(\tilde{e}(t))\pm D^\pm\eps)ds=0\right\}.
\]
Recalling the definition~\eqref{eq:def-tilde-T_i} of $\tilde{T}_i$, we can apply Lemma~\ref{lem_control_reaching_time_of_0} to deduce
that
\begin{equation}
  \label{eq:control-T_Y-T_Z}
  \left| T_\eps^Z-T_\eps^Y\right|\leq  \left| T_\eps^+-T_\eps^-\right|\leq C_1\eps
\end{equation}
for an explicit constant $C_1$. In particular, $T_\eps^Y\leq \tilde T_i\leq T_\eps^Z<s_0$ for $\eps$ small enough. By~\eqref{eq:control-r_Y-r_Z},
this also implies that
\[
  b^Z(T_\eps^Z)-d^Z(T_\eps^Z)\leq f_{i,\ell^\star_0}(\tilde{e}(T^Z_\eps))+D^+\eps\leq
  f_{i,\ell^\star_0}(\tilde{e}(\tilde{T}_i))+C_1\|f\|_\text{Lip}\eps+D^+\eps<0
\]
for $\eps>0$ small enough. Similarly, $b^Y(T_\eps^Y)-d^Y(T_\eps^Y)<0$ if $\eps>0$ is small enough. Therefore,
for such $\eps$, Lemma~\ref{lem:time-inhomo-LBD}~(ii) applies to $Y^K$ and $Z^K$ on $[0,s_0]$.

In both cases, we deduce that
\[
  \lim_{K\to+\infty}\P\left(\forall t\in[0,s_0],\ \left|\frac{\log(1+Y^K(t\log K))}{\log K}-\beta_i(0)+\eps-\int_0^{t\wedge
      T_\eps^Y}(b^Y(s)-d^Y(s))ds\right|\leq \eps\right)=1
\]
and
\[
  \lim_{K\to+\infty}\P\left(\forall t\in[0,s_0],\ \left|\frac{\log(1+Z^K(t\log K))}{\log K}-\beta_i(0)-\eps-\int_0^{t\wedge
      T_\eps^Z}(b^Z(s)-d^Z(s))ds\right|\leq\eps\right)=1.
\]
In view of~\eqref{eq.stoch.dom.LinBD}, this implies that, with probability converging to 1, for all $t\in[0,s_0\wedge\thetainv]$,
\begin{align*}
  \beta_i(0)-2\eps+\int_0^{t\wedge
  T_\eps^Y}(b^Y(s)-d^Y(s))ds \leq\beta_i^K(t)\leq \beta_i(0)+2\eps+\int_0^{t\wedge
  T_\eps^Z}(b^Z(s)-d^Z(s))ds. 
\end{align*}
By~\eqref{eq:control-r_Y-r_Z} and~\eqref{eq:control-T_Y-T_Z}, we have
\begin{align*}
  \beta_i(0)+2\eps+\int_0^{t\wedge T_\eps^Z}(b^Z(s)-d^Z(s))ds & \leq \beta_i(0)+\int_0^{t\wedge
    \tilde{T}_i}f_{i,\ell^\star_0}(\tilde{e}(s))ds +(2+D^+T)\eps+| T_\eps^Z-\tilde{T}_i |\bar{f} \\ & \leq\tilde{\beta}_i(t)+(2+D^+T+C_1\bar{f})\eps
\end{align*}
and similarly for $\beta_i(0)-2\eps+\int_0^{t\wedge T_\eps^Y}(b^Y(s)-d^Y(s))ds$. Therefore, we have proved~\eqref{eq:step-exponents}
with $C:=2+(D^+\vee D^-)T+C_1\bar{f}$.

\bigskip

We conclude this section by noticing that, in the case where $s_0>\tilde{T}_i$ and $f_{i,\ell^\star_0}(\tilde{e}(\tilde{T}_i))<0$, it
also follows from Lemma~\ref{lem:time-inhomo-LBD}~(ii) that, for all $t\in(T^Z_\eps, s_0]$, with probability converging to 1,
$Z^K(t\log K)=0$. By~\eqref{eq.stoch.dom.LinBD} and~\eqref{eq:control-T_Y-T_Z}, for all $t\in(\tilde{T}_i,s_0]$, for $\eps$ small
enough,
\begin{equation}
  \label{eq:estimate-extinction}
  \lim_{K\to+\infty}\P\left(N^K_i(t\log K)\geq 1\text{ and }\thetainv\geq t\right)=0.
\end{equation}

\subsection{Conclusion if $T<t_1$}
\label{sec:ccl-T<t_1}

In the case when $T<t_1$, no exponent $\beta_i$ hit 0 before $T$, therefore we can apply~\eqref{eq:step-exponents} with $s_0=T$, so,
by~\eqref{eq:proba-Gamma->1} and~\eqref{eq:cv-beta-ell_0}, the proof of~\eqref{eq.cv.proba.iter1} will be completed if we prove that, for $\varepsilon$ small enough, with probability
converging to 1, $\thetainv=T$, i.e.
\begin{equation}
  \label{eq:ccl-T<t_1}
  \forall t\in[0,T],\quad \frac{1}{K}\sum_{i\neq\ell^\star_0}N^K_i(t\log K)\leq\varepsilon.
\end{equation}

Since $\beta_i(t)=\tilde{\beta}_i(t)<1$ for all $i\neq\ell^\star_0$ and all $t\leq T$, the continuity of the functions $\beta_i$ ensures the existence
of $\eta>0$ such that
\[
  \forall i\neq\ell^\star_0,\quad\forall t\leq T,\quad\beta_i(t)\leq 1-2\eta.
\]
Therefore, it follows from~\eqref{eq:step-exponents} with $s_0=T$ that, for $\varepsilon$ small enough (and in particular for $\varepsilon<\eta/C$), with probability
converging to 1, $\beta^K_i(\thetainv)<1-\eta$, which implies that $N^K_i(\thetainv\log K)<K^{1-\eta}$. Hence, with
probability converging to 1,
\[
  \frac{1}{K}\sum_{i\neq\ell^\star_0} N^K_i(\thetainv\log K)\leq (d-1)K^{-\eta},
\]
which implies~\eqref{eq:ccl-T<t_1} for $K$ large enough by definition of $\thetainv$. This ends the proof when $T<t_1$.

\subsection{Conclusion if $T>t_1$ and there exists $j$ such that $\beta_j(t_1)=0$}
\label{sec:ccl-case-1}

Assume $T>t_1$ and $\beta_j(t_1)=0$ for some $j\neq\ell^\star_0$. Since $T<S$, this implies that induction did not stop at time
$t_1$, so that there is only a single such $j$, $f_{j,\ell^\star_0}(e(t_1))<0$ and $\beta_i(t_1)<1$ for all $i\neq\ell^\star_0,j$ (in
particular, this case contains Case 1 of Section~\ref{sec:reduct-pb}). Hence, recalling the definition~\eqref{eq:def-tilde-t2} of
$\tilde{t}_2$, we can apply~\eqref{eq:step-exponents} with $s_0$ such that $t_1+\gamma\leq s_0 \leq \tilde{t}_2$, with $\gamma$ defined by \eqref{eq:def-gamma}. In particular, as in Section~\ref{sec:ccl-T<t_1}
above,~\eqref{eq.cv.proba.iter1} will be proved if we prove that, for $\eps$ small enough, with probability converging to 1,
$\thetainv\geq T\wedge (t_1+\gamma)$.

Since $\ell^\star_1=\ell^\star_0$ in the case considered here, we have $\tilde{e}(t)=e(t)$ and $\tilde{\beta}_i(t)=\beta_i(t)$ for
all $i$ and all $t<t_2$. Since $t_1+\gamma<t_2$, we deduce that $\tilde{\beta}_i(t)<1$ for all $i\neq\ell^\star_0$ and all
$t\leq t_1+\gamma$, so we can proceed as in Section~\ref{sec:ccl-T<t_1} above to deduce that, for $\eps$ small enough, with
probability converging to 1, $\thetainv\geq T\wedge (t_1+\gamma)$, so~\eqref{eq.cv.proba.iter1} is proved. This ends the proof if
$t_1<T\leq t_1+\gamma$.

In the case where $T>t_1+\gamma$, we also need to prove~\eqref{eq:cv.proba.t1plusgamma.extinction} and~\eqref{eq:extinction}. Since
$\thetainv\geq t_1+\gamma$ with probability converging to 1, the first property follows from~\eqref{eq:proba-Gamma->1} and
the second property follows from~\eqref{eq:estimate-extinction} (notice that it is enough to prove~\eqref{eq:extinction} for
$t_1<t<t_1+\gamma$ because extinction of a species is irreversible). This ends the proof in the case where $T>t_1$ and
$\beta_j(t_1)=0$ for some $j\neq\ell^\star_0$.

\subsection{Conclusion if $T>t_1$ and there exists $j\neq\ell^\star_0$ such that $\beta_j(t_1)=1$}
\label{sec:ccl-case-2}

Assume $T>t_1$ and $\beta_j(t_1)=1$ for some $j\neq\ell^\star_0$. Since $T<S$, this implies that induction did not stop at time
$t_1$, so that there is only a single such $j$, denoted $\ell^\star_1$. In addition, $f_{\ell^\star_1,\ell^\star_0}(e(t_1))>0$,
$f_{\ell^\star_0,\ell^\star_1}(e(t_1))<0$ and $\beta_i(t_1)\in(0,1)$ for all $i\neq\ell^\star_0,j$ (in particular, this case contains
Case 2 of Section~\ref{sec:reduct-pb}).


\subsubsection{Bound on $\thetainv$}
\label{sec:bound-theta-inv}

Since $f_{\ell^\star_1,\ell^\star_0}(e(t_1))>0$, applying Lemma~\ref{lem_control_reaching_time_of_0} (with $D^\pm=0$) to
$1-\tilde{\beta}_{\ell^\star_1}$ and $1-\tilde{\beta}_{\ell^\star_1}\pm (C+1)\varepsilon$ where the constant $C$ is the constant
in~\eqref{eq:step-exponents}, we deduce that, for $\varepsilon$ sufficiently small and for all $t\leq t_1-L\varepsilon$ where $L:=\frac{4(C+1)}{f_{\ell^\star_1,\ell^\star_0}(e(t_1))}$,
$1-\tilde{\beta}_{\ell^\star_1}(t)-(C+1)\varepsilon>0$, i.e.\ $\tilde{\beta}_{\ell^\star_1}(t)<1-(C+1)\varepsilon$.
By~\eqref{eq:step-exponents}, we deduce that, with probability converging to 1,
\[
  \forall t\leq \left( t_1-L\varepsilon\right)\wedge\thetainv,\quad N^K_{\ell^\star_1}(t\log K)\leq
  K^{1-\varepsilon}.
\]
We conclude that $\thetainv> t_1-L\varepsilon$ as in Section~\ref{sec:ccl-T<t_1}.

Similarly, it also follows from the proof of Lemma~\ref{lem_control_reaching_time_of_0} that
\[
  \tilde{\beta}_{\ell^\star_1}\left(t_1+L\varepsilon\right)>1+(C+1)\varepsilon,
\]
so that, by~\eqref{eq:step-exponents} again, on an event of probability converging to 1, we have
\[
  \thetainv\geq t_1+L\varepsilon\quad\text{implies that}\quad N^K_{\ell^\star_1}\left((t_1+L\varepsilon)\log K\right)\geq K^{1+\varepsilon}-1.
\]
Since $K^{1+\varepsilon}-1>\varepsilon K$ for $K$ large enough, we deduce that $\P(\Gamma^K_\text{inv})\to 1$ when $K\to+\infty$,
where
\[
  \Gamma^K_\text{inv}:=\left\{\thetainv\in\left[ t_1-L\varepsilon \, ,\, t_1+L\varepsilon \right]\right\}.
\]

\subsubsection{Change of dominant population after $\thetainv$}
\label{sec:change-dominant}

On the event $\Gamma^K_\text{inv}$, we apply the strong Markov property at time $\thetainv$ and Lemma~\ref{lem:competition} to the
species $\ell^\star_i$ for $i\in\{0,1\}$ and with $Z^K_t=E^K(t)$. Let us check the assumptions of this lemma.

First, the individual birth rate of species $\ell^\star_i$ at time $t$ takes the form $b_{\ell^\star_i}(E^K(t))$, and its individual
death rate at time $t$ takes the form
\[
  d_{\ell^\star_i}(E^K(t))+\hat{d}^K_{\ell^\star_i}(t,\omega)+c_{\ell^\star_i}(E^K(t),\ell^\star_{i})\frac{N^K_{\ell^\star_i}(t)}{K}+c_{\ell^\star_{i}}(E^K(t),\ell^\star_{1-i})\frac{N^K_{\ell^\star_{1-i}}(t)}{K},
\]
where
\[
  \hat{d}^K_{\ell^\star_i}(t,\omega):=\sum_{j\neq\ell^\star_0,\ell^\star_1}c_{\ell^\star_i}(E^K(t),j)\frac{N^K_{j}(t)}{K}.
\]
Observe that, for all $j\neq \ell^\star_0,\ell^\star_1$, $\beta_j(t_1)<1$, hence there exists $\eta>0$ such that
$\beta_j(t_1)<1-2\eta$ for all $j\neq \ell^\star_0,\ell^\star_1$. Hence, for all
$t\in[t_1-\eta/\|\beta\|_\text{Lip},t_1+\eta/\|\beta\|_\text{Lip}]$, $\tilde\beta_j(t)<1-\eta$ for all $j\neq \ell^\star_0,\ell^\star_1$.
In particular, by~\eqref{eq:step-exponents}, if $\varepsilon<\frac{\eta}{(2C)\vee (L\|\beta\|_\text{Lip})}$, conditionally on the
event $\Gamma^K_\text{inv}$, with probability converging to 1, $N^K_j(\thetainv\log K)\leq K^{1-\eta/2}$ for all
$j\neq \ell^\star_0,\ell^\star_1$. Then, applying Lemma~\ref{lem:continuity-exponents} and reducing $\eta$ if necessary, we deduce
that, with probability converging to 1,
\[
  \forall t\in\left[\thetainv,\thetainv+\frac{\eta}{4C_\text{Lip}}\right],\quad N^K_j(t\log K)\leq K^{1-\eta/4}.
\]
We finally deduce that, if $\varepsilon$ is small enough, conditionally on $\Gamma^K_\text{inv}$, with probability converging to 1,
\[
  \sup_{t\in[\thetainv,\thetainv+\frac{\eta}{4C_\text{Lip}}]}|\hat{d}_{\ell^\star_0}^K(t)|+|\hat{d}_{\ell^\star_1}^K(t)|+\left|\frac{dE^K_t}{dt}\right|\xrightarrow[K\to+\infty]{} 0,
\]
which implies (after applying Markov's property at time $\thetainv\log K$) assumption~\eqref{ass:vanishing_terms} for all $T>0$ in
Lemma~\ref{lem:competition}.

Second, conditionally on the event $\Gamma^K_\text{inv}$, it follows from~\eqref{eq:step-2-initial} that
$|E^K(\thetainv\log K)-\tilde{e}(\thetainv)|\leq A\varepsilon$, hence
\[
  |E^K(\thetainv\log K)-e(t_1)|\leq (A+L\|F\|_\infty)\varepsilon.
\]
Therefore, setting $G:=\{z\in\R^\dimE:|z-e(t_1)|\leq (A+2L\|F\|_\infty)\varepsilon\}$, assumption~\eqref{eq:hyp-2D-Z} in
Lemma~\ref{lem:competition} is satisfied (after applying Markov's property at time $\thetainv\log K$ and conditionally on
$\Gamma^K_\text{inv}$).

Third, since $f_{\ell^\star_0,\ell^\star_1}(e(t_1))<0$ and $f_{\ell^\star_1,\ell^\star_0}(e(t_1))>0$, reducing $\varepsilon>0$ if
necessary, we deduce by continuity of these functions that there exists $\delta>0$ such that for all $z\in G$, $b_{\ell^\star_1}(z)-d_{\ell^\star_1}(z)\geq\delta$,
$f_{\ell^\star_0,\ell^\star_1}(e(t_1))\leq-\delta$ and $f_{\ell^\star_1,\ell^\star_0}(e(t_1))\geq\delta$, hence
assumptions~\eqref{eq:IIF1} and~\eqref{eq:IIF2} of Lemma~\ref{lem:competition} are satisfied.

Fourth, since the jump size of the process $N^K$ is almot surely bounded by $1/K$ which is smaller that $\varepsilon K$ for $K$ large
enough, by~\eqref{eq:Gamma-ell} and by definition of $\thetainv$, conditionally on $\Gamma^K_\text{inv}$, with probability converging
to 1, $\frac{1}{K}(N^K_{\ell^\star_0}(\thetainv\log K), N^K_{\ell^\star_1}(\thetainv\log K))\in H$ where
\[
  H:=\left[n^\star_{\ell^\star_0}(e(t_1))-(B+\|n^\star\|_\text{Lip}\|F\|_\infty L)\varepsilon,
      n^\star_{\ell^\star_0}(e(t_1))+(B+\|n^\star\|_\text{Lip}\|F\|_\infty L)\varepsilon\right]\times[\varepsilon/2,2\varepsilon],
\]
hence, reducing $\eps$ if necessary, assumption~\eqref{eq:hyp-CI-lemme-2D} of Lemma~\ref{lem:competition} is satisfied.

\bigskip

We conclude that there exists $\eta_\varepsilon>0$ and $T_{\varepsilon}<+\infty$ such that, with probability converging to 1,
\begin{equation}
  \label{eq:ccl-new-dominant-1}
  \frac{1}{K}N^K_{\ell^\star_0}(\thetainv\log K+ T_\varepsilon)\in[\eta_\varepsilon ,\varepsilon ]
  \end{equation}
 \begin{equation}
  \label{eq:ccl-new-dominant-1-ter}
  \left|\frac{N^K_{\ell^\star_1}(\thetainv\log K+T_\varepsilon)}{K}-n^\star_{\ell^\star_1}(E^K(\thetainv\log
      K+T_\varepsilon))\right|\leq\varepsilon +\frac{\|n^\star\|_\text{Lip}\|F\|_\infty T_\varepsilon}{\log K}\leq 2\varepsilon, \\
 \end{equation}
 \begin{equation}
  \label{eq:ccl-new-dominant-1-bis}
  \left|E^K(\thetainv\log K+T_\eps)-e(t_1)\right|\leq 2(A+L\|F\|_\infty)\varepsilon
   \end{equation}
   \begin{equation}
  \label{eq:ccl-new-dominant-1-quater}
  \left|\frac{N^K_{\ell^\star_1}(\thetainv\log K+T_\varepsilon)}{K}-n^\star_{\ell^\star_1}(e(t_1))\right|\leq
  2(1+\|n^\star\|_\text{Lip}(A+L\|F\|_\infty))\varepsilon.  
\end{equation}
In addition, applying Lemma~\ref{lem:continuity-exponents} with $\eta=\varepsilon$ and after applying Markov's property at time
$\thetainv$, we deduce from~\eqref{eq:step-exponents} that, for all $j\neq\ell^\star_0,\ell^\star_1$, with probability converging to 1,
\begin{equation}
  \label{eq:ccl-new-dominant-2}
  \left|\beta^K_j(\thetainv+T_\varepsilon/\log K)-\beta_j(t_1)\right|\leq \left(C_\text{Lip}+C+L\|\beta\|_\text{Lip}\right)\varepsilon
\end{equation}

\subsubsection{Study of the exponents after time $\thetainv\log K+T_\varepsilon$}
\label{sec:conclusion-proof}

To conclude the proof, we need to check~\eqref{eq.cv.proba.iter1} and~\eqref{eq:cv.proba.t1plusgamma.invasion}. Define
\[
  \widetilde{\theta}^K_\text{inv}:=(t_1+\gamma)\wedge\inf\left\{t\geq
    \thetainv+T_\varepsilon/\log K:\frac{1}{K}\sum_{j\neq\ell^\star_1}N^K_j(t\log K)> 2\varepsilon\right\}.
\]
Define also the functions $\tilde{e}^{(\varepsilon)}(t)$ for $t\geq\thetainv+T_\varepsilon/\log K$ such that
\[
  \tilde{e}^{(\varepsilon)}(t)=e(t_1)+\int_{\thetainv+T_\varepsilon/\log K}^t F_{\ell^\star_1}(\tilde{e}^{(\varepsilon)}(s))ds
\]
and for all $i\neq \ell^\star_1$, $\tilde{\beta}^{(\varepsilon)}(t)$ for
$t\geq\thetainv\log K+T_\varepsilon$ such that
\[
  \tilde{\beta}_i^{(\varepsilon)}(t)=\beta_i(t_1)+\int_{\thetainv+T_\varepsilon/\log K}^t f_{i,\ell^\star_1}(\tilde{e}^{(\varepsilon)}(s))ds.
\]
Observe that, since $\thetainv\in[t_1-L\varepsilon,t_1+L\varepsilon]$ with probability converging to 1 and $t_1+\gamma<t_2$, the
functions $\tilde{e}^{(\varepsilon)}(t)$ and $\tilde{\beta}^{(\varepsilon)}_i(t)$ for $i\neq\ell^\star_1$ converge in probability in
$L^\infty([\thetainv+T_\varepsilon/\log K,t_1+\gamma])$ to $e(t)$ and $\beta_i(t)$.

Therefore, in view of~\eqref{eq:ccl-new-dominant-1},~\eqref{eq:ccl-new-dominant-1-bis},~\eqref{eq:ccl-new-dominant-1-quater} and~\eqref{eq:ccl-new-dominant-2}, it is
sufficient to prove the following: for all $\varepsilon>0$ small enough, with probability converging to 1, for all
$t\in \left[\thetainv+T_\varepsilon/\log K,\, \widetilde{\theta}^K_\text{inv}\right]$,
\[
  |E^K(t\log K)-\tilde{e}^{(\varepsilon)}(t)|\leq
  A'\varepsilon,\quad \left|\frac{N^K_{\ell^\star_1}(t\log K)}{K}-n^\star_{\ell^\star_1}\left(\tilde{e}^{(\varepsilon)}(t)\right)\right|\leq
  B'\varepsilon,\quad 
|\beta^K(t)-\tilde{\beta}^{(\varepsilon)}(t)|\leq C'\varepsilon
\]
and $\widetilde{\theta}^K_\text{inv}=t_1+\gamma$.

Applying Markov's property at time $\thetainv\log K+T_\varepsilon$, these estimates are very similar to~\eqref{eq:proba-Gamma->1},
\eqref{eq:step-exponents} (in the simplest case where $s_0=\gamma<\tilde{T}_i$) and Section~\ref{sec:ccl-T<t_1} (since
$\gamma<t_2-t_1$) and they indeed will be obtained following very similar arguments. The main differences are that we do not have
convergence of the total population size and environmental variable at initial time (instead, we only have the
dominations~\eqref{eq:ccl-new-dominant-1-bis}-\eqref{eq:ccl-new-dominant-1-quater}), nor convergence of the exponents at initial time to their limits (instead, we only
have the domination~\eqref{eq:ccl-new-dominant-2}) and neither the exponent $\beta^K_{\ell^\star_0}$ at initial time does converge to
a constant smaller than 1 (instead, by~\eqref{eq:ccl-new-dominant-1}, it converges to 1). Since the proof is based on very similar
comparison arguments, we only explain below the main differences.

\bigskip

In view of~\eqref{eq:ccl-new-dominant-1-ter}, the argument of Section~\ref{sec.dominant} is exactly the same, except that the bounding
processes $Y^K$ and $Z^K$ must be taken with $Y^K_0/K=n^\star_{\ell^\star_1}(E^K(\thetainv\log K+T_\varepsilon))-2\varepsilon$ and
$Z^K_0/K=n^\star_{\ell^\star_1}(E^K(\thetainv\log K+T_\varepsilon))+2\varepsilon$. The conclusion of Section~\ref{sec.dominant} then
becomes that, for $\varepsilon>0$ small enough, with probability converging to 1,
\[
  \sup_{t\in[\thetainv+T_\varepsilon/\log K,\widetilde{\theta}^K_\text{inv}]}\left|\frac{N_{\ell_1^\star}^K(t\log K)}{K}-n_{\ell_1^\star}^\star(E^K(t\log K))\right|\leq 8\left(1+\frac{\bar{c}}{\underline{c}}\right) \varepsilon.
\]

The argument of Section~\ref{sec.environment} is also the same: defining for all $t\in[\thetainv+T_\varepsilon/\log K,\,t_1+\gamma]$,
$\psi^K(t):=E^K(t\log K) - \tilde{e}^{(\varepsilon)}(t)$, the computation of Section~\ref{sec.environment} gives, for all $t<\widetilde{\theta}^K_\text{inv}$
\[
  \left|\psi^K(t)\right|^2  \leq \left(\psi^K(0)^2+\left(2+8\left(1+\frac{\bar c}{\underline
        c}\right)\right)^2\,\eps^2\left(3+2\|n^\star\|_\text{Lip}\right)^{-1}\right)\exp\left(\|F\|_{\text{Lip}}\left(3+2\|n^\star\|_\text{Lip}\right) t\right).
\]
Using~\eqref{eq:ccl-new-dominant-1} and~\eqref{eq:ccl-new-dominant-1-bis}, we deduce as in Section~\ref{sec.environment} that, for $\varepsilon$ small enough, the
event
\begin{align*}
  \widetilde{\Gamma}^K:=\left\{\forall t\in [\thetainv+T_\varepsilon/\log K,\, \widetilde{\theta}^K_\text{inv}],\ \left|E^K(t\log
      K)-\tilde{e}^{(\varepsilon)}(t)\right|\leq  A'\eps, \text{ and }\left|\frac{N^K_{\ell_1^\star}(t\log K)}{K} -
      n_{\ell_1^\star}^\star(\tilde{e}^{(\varepsilon)}(t))\right|\leq B'\eps\right\}
\end{align*}
satisfies
\[
  \lim_{K\to+\infty}\P\left(\widetilde{\Gamma}^K\right)=1,
\]
 where
\[
  A':=\sqrt{4(A+L\|F\|_\infty)^2+\left(10+8\frac{\bar{c}}{\underline{c}}\right)^2\left(3+2\|n^\star\|_\text{Lip}\right)^{-1}}\exp
  \left(\|F\|_\text{Lip}\left(\frac{3}{2}+\|n^\star\|_\text{Lip}\right) (t_1+\gamma)\right)\]
and
\[
  B':=8\left(1+\frac{\bar{c}}{\underline{c}}\right)+A'\|n^\star\|_\text{Lip}.
\]

Concerning the argument of Section~\ref{sec.stoch.domination}, we need to be more careful for the exponent $\beta^K_{\ell^\star_0}$.
Let us begin with the exponents $\beta^K_j$ for $j\neq\ell^\star_0,\ell^\star_1$: since $t_1+\gamma<t_2$, $\beta_i(t)\in(0,1)$ for
all $i\neq\ell^\star_0,\ell^\star_1$ and all $t\in[0,t_1+\gamma]$, hence there exists $\eta>0$ such that, for all $\varepsilon>0$
small enough and $i\neq\ell^\star_0,\ell^\star_1$,
\begin{equation}
  \label{eq:blabla-joie}
  \forall t\in[\thetainv\log K+T_\varepsilon,\ t_1+\gamma]\quad \tilde{\beta}^{(\varepsilon)}_i(t)\in[\eta,1-\eta].
\end{equation}
In particular, the same comparison argument as in Section~\ref{sec.stoch.domination} applies to all the exponents $\beta^K_i$ for
$i\neq\ell^\star_0,\ell^\star_1$ replacing $A$ by $A'$, $B$ by $B'$ and $\tilde{e}$ by $\tilde{e}^{(\varepsilon)}$, with the
difference that, by~\eqref{eq:ccl-new-dominant-2}, the initial values of the bounding processes $Y^K$ and $Z^K$ must be taken as
\[
  Y^K(\thetainv\log K+T_\varepsilon)=\lfloor K^{\beta_i(t_1)-(C_\text{Lip}+C+L\|\beta\|_\text{Lip})\varepsilon}\rfloor
\]
and
\[
  Z^K(\thetainv\log K+T_\varepsilon)=\lfloor K^{\beta_i(t_1)+(C_\text{Lip}+C+L\|\beta\|_\text{Lip})\varepsilon}\rfloor.
\]
The argument of Section~\ref{sec.stoch.domination} then gives that, for $\varepsilon>0$ small enough, with probability converging to
1,
\begin{equation}
  \label{eq:blabla-ccl}
  \forall t\in[\thetainv+T_\varepsilon/\log K,\, \widetilde{\theta}^K_\text{inv}],\ \left|\beta^K_i(t)-\tilde{\beta}^{(\varepsilon)}_i(t)\right|\leq C'\eps,
\end{equation}
here 
\[
  C':=1+C_\text{Lip}+C+L\|\beta\|_\text{Lip}+D'(t_1+\gamma)
\]
and
\[
 D':=A'\left(\|b\|_\text{Lip}+\|d\|_\text{Lip}+\bar{n}^\star\|c\|_\text{Lip}\right)+A'B'\|c\|_\text{Lip}+(B'+2)\bar{c}.
\]

To study the exponent $\beta^K_{\ell^\star_0}$, we use the same comparison arguments as in Section~\ref{sec.stoch.domination} between
$N^K_{\ell^\star_0}$ and the linear birth-death process $Y^K$ and $Z^K$ with time inhomogeneous birth and death rates $b^Y$, $d^Y$
and $b^Z$ and $d^Z$ respectively obtained by replacing $A$ by $A'$, $B$ by $B'$, $\tilde{e}$ by $\tilde{e}^{(\varepsilon)}$, $\ell_0^\star$ by $\ell_1^\star$ and $\bar c \varepsilon$ by $2\bar c \varepsilon$
in~\eqref{eq:blabla-def-b},~\eqref{eq:blabla-def-d-Y} and~\eqref{eq:blabla-def-d-Z}, and with the initial conditions
\begin{equation}
  \label{eq:blabla-CI}
  Y^K(\thetainv\log K+T_\varepsilon)=\lfloor \eta_\varepsilon K\rfloor\quad\text{and}\quad Z^K(\thetainv\log K+T_\varepsilon)=\lceil
  \varepsilon K\rceil.  
\end{equation}
By~\eqref{eq:ccl-new-dominant-1}, with probability converging to 1,
\[
  \forall t\in[\thetainv\log K +T_\varepsilon,\widetilde{\theta}^K_\text{inv}],\quad Y^K(t)\leq N^K_{\ell^\star_0}(t)\leq Z^K(t).
\]
We can use Lemma~\ref{lem:time-inhomo-LBD} to study the convergence of the exponents of the processes $Y^K$ and $Z^K$ as in
Section~\ref{sec.stoch.domination}. Note however that the initial conditions of the processes $Y^K$ and $Z^K$ are slightly different
as in Lemma~\ref{lem:time-inhomo-LBD}. The reader can easily check from the proof of this lemma that its conclusions remain valid for
initial conditions of the form~\eqref{eq:blabla-CI}. This entails~\eqref{eq:blabla-ccl} for $i=\ell^\star_0$.

It remains to prove that $\tilde{\theta}^K_\text{inv}=t_1+\gamma$ with probability converging to 1. For that, we need more precise estimates on $Z^K$. By continuity of $e(t)$ and the convergence of $\tilde{e}^{(\varepsilon)}$ to $e$, reducing $\varepsilon>0$ is necessary, we have
\[
  \forall t\in\left[\thetainv+T_\varepsilon/\log K,\, \thetainv+T_\varepsilon/\log K+ \frac{4C'\varepsilon}{\vert f_{\ell^\star_0,\ell^\star_1}(e(t_1))\vert}\right],\quad
  f_{\ell^\star_0,\ell^\star_1}(\tilde{e}^{(\varepsilon)}(t))<\frac{f_{\ell^\star_0,\ell^\star_1}(e(t_1))}{2}<0. 
\]
This implies that
\begin{equation}
  \label{eq:blabla-youpi}
  \tilde{\beta}^{(\varepsilon)}_{\ell^\star_0}\left(\thetainv+T_\varepsilon/\log
    K+\frac{4C'\varepsilon}{\vert f_{\ell^\star_0,\ell^\star_1}(e(t_1))\vert}\right)\leq 1-2C'\varepsilon.  
\end{equation}

Let us prove that, with probability converging to 1, $Z^K(t)/K\leq 3\varepsilon/2$ for all
$t\in[\thetainv\log K+T_\varepsilon, \thetainv\log K+T_\varepsilon+4C'\varepsilon\log K/|f_{\ell^\star_0,\ell^\star_1}(e(t_1))|]$.
Reducing $\varepsilon>0$ again if necessary, we can assume that, for all $t\in[\thetainv+T_\varepsilon/\log K,
\thetainv+T_\varepsilon/\log K+4C'\varepsilon/\vert f_{\ell^\star_0,\ell^\star_1}(e(t_1))\vert ]$,
\[
  b^Z(t)\leq\bar{b}^Z,\quad d^Z(t)\geq\underline{d}^Z \quad\text{and}\quad \bar{b}^Z<\underline{d}^Z,
\]
for two constants $\bar{b}^Z$ and $\underline{d}^Z$. This means that there exists a (time-homogeneous) linear birth-death process
$\tilde{Z}^K$ with individual birth rate $\bar{b}^Z$ and individual death rate $\underline{d}^Z$ such that
$\tilde{Z}^K(\thetainv\log K+T_\varepsilon)=Z^K(\thetainv\log K+T_\varepsilon)=\lceil \varepsilon K\rceil$ and
\[
  \forall t\in\left[\thetainv\log K+T_\varepsilon,\, \thetainv\log K+T_\varepsilon+
    \frac{4C'\varepsilon\log K}{\vert f_{\ell^\star_0,\ell^\star_1}(e(t_1))\vert}\right],\quad Z^K(t)\leq\tilde{Z}^K(t).
\]
Standard properties of subcritical linear birth-death imply that, with probability converging to 1, $\tilde{Z}^K(t)\leq 3\varepsilon K/2$ for all
$t\geq \thetainv\log K+T_\varepsilon$ (e.g.\ one may use the induced Markov chain constructed from $\tilde{Z}^K$ is a biased random
walk, with negative bias). We deduce that, with probability converging to 1,
\begin{equation}
  \label{eq:blablabla-last}
  \forall t\in\left[\thetainv\log K +T_\varepsilon,\, \widetilde{\theta}^K_\text{inv}\wedge\left(\thetainv\log K +T_\varepsilon+
    \frac{4C'\varepsilon\log K}{\vert f_{\ell^\star_0,\ell^\star_1}(e(t_1))\vert}\right)\right],\quad \frac{N^K_{\ell^\star_0}(t)}{K}\leq\frac{3}{2}\varepsilon.  
\end{equation}

We can now conclude using the argument of Section~\ref{sec:ccl-T<t_1}. First, we deduce from this argument and~\eqref{eq:blabla-joie} that, with probability converging to 1,
\[
  \forall t\in\left[\thetainv\log K +T_\varepsilon,\, \widetilde{\theta}^K_\text{inv}\wedge\left(\thetainv\log K +T_\varepsilon+
    \frac{4C'\varepsilon\log K}{\vert f_{\ell^\star_0,\ell^\star_1}(e(t_1))\vert}\right)\right],\quad
\frac{1}{K}\sum_{i\neq\ell^\star_0,\ell^\star_1}N^K_i(t)\leq (d-2)K^{-\eta/2}.
\]
Combining this with~\eqref{eq:blablabla-last}, we deduce that if $\thetainv\log K +T_\varepsilon+
  \frac{4C'\varepsilon\log K}{\vert f_{\ell^\star_0,\ell^\star_1}(e(t_1))\vert}<t_1+\gamma$ then, with probability converging to 1,
\[
  \widetilde{\theta}^K_\text{inv}> \thetainv\log K +T_\varepsilon+
  \frac{4C'\varepsilon\log K}{\vert f_{\ell^\star_0,\ell^\star_1}(e(t_1))\vert}.
\]
Then, after time $\thetainv\log K +T_\varepsilon+ \frac{4C'\varepsilon\log K}{\vert f_{\ell^\star_0,\ell^\star_1}(e(t_1))\vert}$,
using~\eqref{eq:blabla-youpi}, that~\eqref{eq:blabla-ccl} is satisfied for all $i\neq\ell^\star_1$, the convergence of $\tilde{\beta}^{(\varepsilon)}_{\ell^\star_0}$ to $\beta_{\ell^\star_0}$ and that $t_1+\gamma<t_2$,  we deduce from the argument of
Section~\ref{sec:ccl-T<t_1} that, with probability converging to 1,
\[
  \forall t\in\left[\thetainv\log K +T_\varepsilon +\frac{4C'\varepsilon\log K}{\vert f_{\ell^\star_0,\ell^\star_1}(e(t_1))\vert},\, t_1+\gamma\right],\quad
\frac{1}{K}\sum_{i\neq\ell^\star_1}N^K_i(t)\leq (d-1)K^{-\frac{\tilde\eta}{2}},
\]
with $\tilde \eta \leq \eta$ such that 
\[
  \forall t\in\left[\thetainv\log K +T_\varepsilon +\frac{4C'\varepsilon\log K}{\vert f_{\ell^\star_0,\ell^\star_1}(e(t_1))\vert},\, t_1+\gamma\right],\quad
\tilde\beta_{\ell^\star_0}(t)<1-\tilde \eta.
\]
This implies that $\tilde{\theta}^K_\text{inv}=t_1+\gamma$ with probability converging to 1. Hence~\eqref{eq.cv.proba.iter1}
and~\eqref{eq:cv.proba.t1plusgamma.invasion} follow from the previous estimates and the proof of Theorem~\ref{thm:BDP-max-2} is
complete.

\section{Characterization of the initial conditions for which induction never stops}
\label{sec:induction-never-stops}

The purpose of this section is to study the sequence of times $(t_k)_k$ defined in Theorem~\ref{thm:BDP-max-2}. We are particularly interested by conditions on the model parameters and the initial conditions $\beta(0)$ and $e(0)$ ensuring either that the full sequence $(t_k)_k$ is well defined or $S=+\infty$. In the first case, convergence to the limit dynamics is ensured up to a possibly finite accumulation point of the times of change of resident populations. In the second case, the convergence of Theorem~\ref{thm:BDP-max-2} holds true on the complete time interval $\R_+$.

In Subsection~\ref{sec:Baire-d=2}, we show that under mild conditions, when $d=2$, induction never stops for all initial conditions out of a Baire meager set. In Subsection~\ref{sec:Baire-d>2}, we extend this result to $d>2$. These two results entail Theorem~\ref{thm:Baire-d>2}. Finally, we give conditions ensuring that $S=+\infty$ when induction never stops in Subsection~\ref{sec:S-infty}, hence proving Proposition~\ref{prop:S-infty}. Recall that we say that induction never stops in Theorem~\ref{thm:BDP-max-2} if either $S=+\infty$ or $t_k$ is well-defined for all $k$ and $S=\sup_k t_k$. Recall also the definition of the vector fields $F_i$ in~\eqref{eq:def-F_i}.

\subsection{The case $d=2$}
\label{sec:Baire-d=2}

Before stating our result, some topological reminder is in order. Recall that if $X$ is a topological space, a set $A \subseteq X$ is said to be nowhere dense if its closure has empty interior. A (Baire) meager is the (at most) countable union of nowhere dense sets. Of course, if $g: X\mapsto X'$ is a homeomorphism between topological spaces, then $g$ and $g^{-1}$ preserve nowhere dense sets and meager sets. Finally, recall the Baire Category Theorem, which states that on a complete metric space, every meager set has empty interior (see e.g. \cite{brezis2011}, Chapter 2, Thm. 2.1).

\begin{thm}
    \label{thm:Baire-d=2-k=1}
    Assume $d=2$, $\dimE\in\mathbb{N}^*$ and that $f_{1,2}$, $f_{2,1}$, $F_1$ and $F_2$ are Lipschitz bounded functions. Assume also that the functions $f_{1,2}$ and $f_{2,1}$ are non-identically 0 on any open ball of $\R^\dimE$. Assume finally that invasion implies fixation, in the sense that, for all $e\in\R^\dimE$, if $f_{2,1}(e)>0$, then $f_{1,2}(e)<0$, and if $f_{1,2}(e)>0$, then $f_{2,1}(e)<0$. Then there exists a subset $D_0$ of $D:=\{1\}\times(0,1)\times\mathbb{R}^\dimE\cup(0,1)\times\{1\}\times\mathbb{R}^\dimE$ such that
    \[
    \left\{(\beta_2,e)\in(0,1)\times\mathbb{R}^\dimE,\,(1,\beta_2,e)\in D_0\right\}
    \text{ and }
     \left\{(\beta_1,e)\in(0,1)\times\mathbb{R}^\dimE,\,(\beta_1,1,e)\in D_0\right\}
    \]
    are Baire meager subsets of $(0,1)\times\mathbb{R}^\dimE$ and, for all initial conditions $(\beta_1(0),\beta_2(0),e(0))\in D\setminus D_0$, induction never stops in Theorem~\ref{thm:BDP-max-2}, which means that either $S=+\infty$ or 
    the full sequence $(t_k)_{k\geq 0}$ is well-defined, increasing and $S=\sup_{k\geq 1}t_k$.
\end{thm}

The remaining of this section is devoted to the proof of this result.

\paragraph{Definitions and notations:}
From now on $f_2=f_{2,1}, f_1=f_{1,2}$. Set $D_1:=\{1\} \times (0,1) \times \R^\dimE$ and note that, by symmetry, it is enough to prove the result for initial conditions in $D_1$ only, so we will only consider this case below. 

For $i\in\{1,2\}$, we define $\phi_i:\R_+\times\R^\dimE\to \R^\dimE$ as the flow of the vector field $F_i$, i.e.\ for $e(0)\in \R^\dimE$, we set $\phi_i(t,e(0))=e(t)$ where $(e(t),t\geq 0)$ solves the ODE $\dot{e}=F_i(e)$.

For $(\beta,e)= (1, \beta_2,e) \in D_1$ we define
\[
t_1(\beta_2, e)=
\inf\left\{ t\ge 0: \beta_2+\int_{0}^t f_2(\phi_1(s, e))ds \in \{0,1\}\right\}.
\]
Notice that this definition coincides with the time $t_1$ defined in Theorem~\ref{thm:BDP-max-2}. 

We also introduce the following notations.
\begin{enumerate}
\item We define the sets
\begin{align*}
A_1&:= \{(\beta_2, e)\in(0,1)\times\R^\dimE: t_1(\beta_2, e) < \infty\} \\
\hat A_1&:= \left\{(\beta_2, e) \in A_1: \beta_{2}+\int_{0}^{t_1(\beta_2, e)} f_{2}(\phi_1(s,e))ds =0\right\}.
\end{align*}
\item Let $\Psi_1:A_1\to \R_+\times \R^\dimE$ be defined for all $(\beta_2,e)\in A_1$ by 
\[
\Psi_1(\beta_2, e)= (t_1(\beta_2, e), \phi_1(t_1(\beta_2, e), e)).
\]
\item Let $B=B_1 \cup B_2$ where $B_1:= \R_+ \times \{e\in\R^\dimE: f_1(e)> 0, f_2(e)< 0\}$, $B_2:= \R_+ \times \{e\in\R^\dimE: f_2(e)> 0, f_1(e)< 0\}$. Notice that $B=\R_+ \times \{e\in\R^\dimE: f_1(e)f_2(e)< 0\}$
\item Let $C=B\cap(\R_+ \times \{e \in \R^\dimE: s_1(e)< \infty\})$, where, for all $e$ such that $(t,e)\in B$ for some (hence all) $t\in\mathbb{R}_+$,
\begin{align*}
s_1(e):= \begin{cases}
\inf\{t > 0: 1 + \int_{0}^t f_2(\phi_1(s,e))ds \in \{0,1\}\}, \quad \text{ if } f_2(e)< 0, f_1(e)>0\\
\inf\{t > 0: 1 + \int_{0}^t f_1(\phi_2(s,e))ds \in \{0,1\}\}, \quad \text{ if } f_1(e)< 0, f_2(e)>0\\
\end{cases}
\end{align*}
Note that $s_1(e)>0$.
We also set $C_i=C\cap B_i$ for $i\in\{1,2\}$.
\item We define
\begin{align*}
    \hat C_1&:= \left\{(t, e) \in C_1: 1+ \int_{0}^{s_1(e)} f_2(\phi_1(s,e))ds =0\right\} \\
    \hat C_2&:= \left\{(t, e) \in C_2: 1+ \int_{0}^{s_1(e)} f_1(\phi_2(s,e))ds =0\right\}
\end{align*}
and $\hat C= \hat C_1 \cup \hat C_2$.
\item Let $\Upsilon:C\to\R_+\times\R^\dimE$ be defined for all $(t,e)\in C$ by 
\[\Upsilon(t,e)=\begin{cases}
    (t+s_1(e), \phi_1(s_1(e),e)) & \text{if }(t,e)\in C_1, \\ (t+s_1(e), \phi_2(s_1(e),e)) & \text{if }(t,e)\in C_2.
\end{cases}
\]
We denote by $\Upsilon_i$ the restriction of $\Upsilon$ to $C_i$ for $i\in\{1,2\}$.
\item Finally, we introduce the set
\[
G_0= \{e \in \R^\dimE: f_1(e)f_2(e)=0\}.
\]
\end{enumerate}

\paragraph{How we shall use the previous definitions:}
Our goal here is to explain how the previous definitions can be used to study the times $t_k$ and the different outcomes in Theorem~\ref{thm:BDP-max-2}.
\begin{enumerate}
    \item For $(\beta,e)\in D_1$, $t_1=t_1(\beta_2,e)$.
    \item For $(\beta_2,e)\in ((0,1)\times\R^\dimE)\setminus A_1$, $S=t_1=+\infty$ and induction never stops.
    \item For $(\beta,e)\in A_1\setminus\hat{A}_1$ such that $\phi_1(t_1(\beta_2,e),e)\not\in B$, we have $S=t_1<+\infty$ and induction stops at time $t_1$.
    \item For $(\beta,e)\in A_1\setminus\hat{A}_1$ such that $\phi_1(t_1(\beta_2,e),e)\in B$ (note that this actually implies that $\phi_1(t_1(\beta_2,e),e)\in B_2$), 
    \begin{enumerate}
        \item we have $t_2=t_1(\beta_2,e)+s_1(\phi_1(t_1(\beta_2,e),e))$,
        \item if $\phi_1(t_1,e)\in B\setminus C_2$, then $S=t_2=+\infty$ and induction never stops,
        \item otherwise, if $\phi_1(t_1,e)\in C_2\setminus\hat{C}_2$ and $\phi_2(s_1(\phi_1(t_1,e)),\phi_1(t_1,e))\in G_0$, or if $\phi_1(t_1,e)\in \hat{C}_2$ and $f_1(\phi_2(s_1(\phi_1(t_1,e)),\phi_1(t_1,e)))=0$, then $t_2(\beta,e)=S$ and induction stops at time $t_2$,
        \item otherwise, either $S>t_2(\beta,e)$ and induction does not stop at time $t_2$.
    \end{enumerate}
    \item Otherwise, for $(\beta,e)\in \hat{A}_1$, observe that $f_2(\phi_1(t_1,e))\leq 0$, so
    \begin{enumerate}
        \item if $f_2(\phi_1(t_1,e))=0$, then $S=t_1<+\infty$ and induction stops at time $t_1$,
        \item if $f_2(\phi_1(t_1,e))<0$, then $t_2=S=+\infty$ and induction never stops.
    \end{enumerate}
\end{enumerate}
Cases 3 and 4 correspond to the cases where at time $t_1$ the two populations have exponent 1. Case 5 corresponds to the case where a species 2 goes extinct at time $t_1$. We conclude that induction stops in Theorem~\ref{thm:BDP-max-2} in cases 3, 4(c) and 5(a). Induction never stops in cases 2, 4(b) and 5(b). In case 4(d), to decide if induction may or may not stop after time $t_2$, we can proceed recursively as above. The proof below is based on these arguments.

\paragraph{Definition of the sets of initial conditions for which induction stops and first properties:}
Now we define inductively the sets
\[
G_1:= \Psi_1^{-1} (\R_+ \times G_0)
\]
and, given $G_1,\ldots,G_n$ for some $n\geq 1$, we define
\begin{equation}
    \label{eq:F_2p+1}
    G_{n+1}:=
    \begin{cases}
        \Psi_1^{-1}\left(\Upsilon_2^{-1}\left(\ldots\Upsilon_1^{-1}(\R_+  \times G_0)\ldots\right)\right)\setminus (G_1\cup\ldots\cup G_n), & \text{if $n$ is even,} \\ \Psi_1^{-1}\left(\Upsilon_2^{-1}\left(\ldots\Upsilon_2^{-1}(\R_+\times G_0)\ldots\right)\right)\setminus (G_1\cup\ldots\cup G_n), & \text{if $n$ is odd,}
    \end{cases}
\end{equation}
where the functions $\Upsilon_1$ and $\Upsilon_2$ appear $n$ times in total. 

The above sets have the following interpretation: the set of $(\beta_1,\beta_2,e)=(1,\beta_2,e)$ such that $(\beta_2,e)\in G_{n}$ is the set of initial conditions such that induction did not stop until time $t_{n-1}$, $t_{n}<+\infty$ and induction may stop at time $t_n$ in Theorem~\ref{thm:BDP-max-2}. 

\begin{rem}
Note that, in the previous definition, we cannot write that $\Psi_1^{-1}(\Upsilon_2^{-1}(\R_+\times G_0))=(\Upsilon_2\circ\Psi_1)^{-1}(\R_+\times G_0)$ because the function $\Upsilon_2$ is not defined on the full image of $\Psi_1$, even restricted to $D\setminus G_1$. Indeed, the image of $\Psi_1$ may contain points $(t,e)$ such that $s_1(e)=+\infty$ or such that $f_1(e)<0$ and $f_2(e)<0$, which do not belong to the domain of definition of $\Upsilon_2$. However, this will not be a problem since in the first case, induction never stops, and the second case can only occur when one trait goes extinct at time $s_1(e)$; in such case induction also never stops. This remark is important to understand the proof of Proposition~\ref{prop:F_k} below.
\end{rem}

\begin{prop}
\label{prop:F_k}
We have the following properties:
\begin{description}
    \item[\textmd{(i)}] For all $(\beta,e)\in D_1$ such that $(\beta_2,e)\not\in G_1\cup(A_1\setminus\hat{A}_1)$, induction never stops.
    \item[\textmd{(ii)}] For all $(\beta,e)\in D_1$ such that $(\beta_2,e)\in A_1\setminus (\hat{A}_1\cup G_1\cup G_2 )$ and $\Psi_1(\beta_2,e)\not\in C\setminus\hat{C}$, induction never stops.
    \item[\textmd{(iii)}] For all $n\geq 3$, induction never stops for all $(\beta,e)\in D_1$ such that $(\beta_2,e)\not\in (G_1\cup\ldots\cup G_n)$, $(\beta_2,e)\in A_1\setminus\hat{A}_1$, $\Psi_1(\beta_2,e)\in C\setminus\hat{C}$, $\Upsilon\circ\Psi_1(\beta_2,e)\in C\setminus\hat{C}$, $\ldots$, $\Upsilon^{\circ (n-3)}\circ\Psi_1(\beta_2,e)\in C\setminus\hat{C}$ and $\Upsilon^{\circ (n-2)}\circ\Psi_1(\beta_2,e)\not\in C\setminus\hat{C}$.
\end{description}
In particular, for all $(\beta,e)\in D_1$ such that $(\beta_2,e)\not\in \bigcup_{n\geq 1} G_n$, induction never stops.
\end{prop}

Note that in Point (ii) above, all the quantities $\Upsilon^{\circ \ell}\circ\Psi(\beta,e)$ are well-defined because $\Upsilon^{\circ (\ell-1)}\circ\Psi(\beta,e)$ is assumed to belong to $C\setminus\hat{C}$.

\begin{proof}
    We start by proving (i). Assume that $(\beta_2,e)\not\in G_1\cup(A_1\setminus\hat{A}_1)$. The fact that $(\beta,e)\not\in A_1\setminus\hat{A}_1$ implies that (a) either $t_1(\beta_2,e)=+\infty$ (i.e. $(\beta_2,e)\not\in A_1$), (b) or $t_1(\beta_2,e)<\infty$ and $\beta_2+\int_{0}^{t_1(\beta_2, e)} f_2(\phi_1(s,e))ds =0$. In case (a), $S=+\infty$ and induction never stops. In case (b), the fact that $(\beta_2,e)\not\in G_1$, means that $f_2(\phi_1(t_1(\beta_2,e),e))\neq 0$ and since $\beta_2+\int_{0}^{t_1(\beta_2, e)} f_2(\phi_1(s,e))ds$ hits 0 from positive values, this implies that $f_2(\phi_1(t_1(\beta_2,e),e))<0$. By Theorem~\ref{thm:BDP-max-2}, this means that $\beta_2(s)=0$ for all $s\geq t_1(\beta_2,e)$. Since there are only two species, this implies that $\beta_1(s)=1$ for all $s\geq t_1(\beta_2,e)$ and so $S=t_2=+\infty$, so induction never stops. This ends the proof of~(i).

    Let us come to the proof of~(ii). Take $(\beta,e)$ as in~(ii). The fact that $(\beta_2,e)\in A_1\setminus(\hat{A}_1\cup G_1)$ means that $t_1(\beta_2,e)<+\infty$, $\beta_1(t_1(\beta_2,e))=\beta_2(t_1(\beta_2,e))=1$ and $f_i(e_1)\not=0$ for $i\in\{1,2\}$ where $e_1:=\phi_1(t_1(\beta_2,e),e)$. This implies that $f_2(e_1)>0$ and hence, by the invasion implies fixation assumption, we must have $f_1(e_1)<0$, so $(t_1(\beta_2,e),e_1)\in B_2$. In particular, induction does not stop at time $t_1$, and we can apply the function $s_1$ to $e_1$ (with possibly infinite value). The assumption that $\Psi_1(\beta_2,e)\not\in C\setminus\hat{C}$ implies that either $t_2=t_1(\beta_2,e)+s_1(e_1)$ is infinite, or it is finite and $\beta_1(t_2)=0$. In the last case, as above, $\beta_1(s)=0$ and $\beta_2(s)=1$ for all $s\geq t_2$. Hence induction never stops in both cases.
    
    Now, we assume $n\geq 3$ and prove~(iii). The fact that $\Psi_1(\beta_2,e)\in C\setminus\hat{C}$ and $(\beta_2,e)\not\in G_2$ means that $t_2=t_1(\beta_2,e)+s_1(e_1)<\infty$, $\beta_1(t_2)=\beta_2(t_2)=1$ and $f_i(e_2)\neq 0$ for $i\in\{1,2\}$ where $e_2:=\phi_2(s_1(e_1),e_1)$.
    Again this implies that induction does not stop at time $t_2$, and we can apply the function $s_1$ to $e_2$. Inductively, we deduce that $t_{n-1}=t_{n-2}+s_1(e_{n-2})$, $t_{n-1}<\infty$, $\beta_1(t_{n-1})=\beta_2(t_{n-1})=1$, induction does not stop up to time $t_{n-1}$ and $(t_{n-1},e_{n-1})\in B$ where 
    \[
    e_{n-1}=\begin{cases}
        \phi_1(s_1(e_{n-2}),e_{n-2}) & \text{if $n$ is even}, \\ \phi_2(s_1(e_{n-2}),e_{n-2}) & \text{if $n$ is odd}.
    \end{cases}
    \] 
    Therefore, we can apply the function $s_1$ once again to $e_{n-1}$, with possibly infinite value.

    It remains to use the assumption that $\Upsilon^{\circ (n-2)}\circ\Psi_1(\beta_2,e)\not\in C\setminus\hat{C}$ and $(\beta_2,e)\not\in G_n$. Two possibilities can occur: (a) $t_n:=t_{n-1}+s_1(e_{n-1})=+\infty$, in which case induction never stops, or (b) $t_n<+\infty$ and $\beta_i(t_n)=0$ with $i=1$ if $n$ is even, or $i=2$ if $n$ is odd. Since $(\beta_2,e)\not\in G_n$, we deduce as above that $f_{i}(e_n)<0$, where
    \[
    e_{n}=\begin{cases}
        \phi_2(s_1(e_{n-1}),e_{n-1}) & \text{if $n$ is even}, \\ \phi_1(s_1(e_{n-1}),e_{n-1}) & \text{if $n$ is odd}.
    \end{cases}
    \] 
    Hence, again induction never stops. This completes the proof of~(ii).

    It only remains to check the last claim of Proposition~\ref{prop:F_k}. Given $(\beta,e)\in D_1$ such that $(\beta_2,e)\not\in \bigcup_{n\geq 1} G_n$, either (a) $(\beta_2,e)\not\in A_1\setminus\hat{A}_1$, or (b) there exists $n\geq 2$ such that $(\beta_2,e)\in A_1\setminus\hat{A}_1$, $\Psi_1(\beta_2,e)\in C\setminus\hat{C}$, $\Upsilon\circ\Psi_1(\beta_2,e)\in C\setminus\hat{C}$, $\ldots$, $\Upsilon^{\circ (n-3)}\circ\Psi_1(\beta_2,e)\in C\setminus\hat{C}$ and $\Upsilon^{\circ (n-2)}\circ\Psi_1(\beta_2,e)\not\in C\setminus\hat{C}$, or (c) $(\beta_2,e)\in A_1\setminus\hat{A}_1$ and for all $\ell\geq 0$, $\Upsilon^{\circ \ell}\circ\Psi_1(\beta_2,e)\in C\setminus\hat{C}$. Point (i) and (ii) imply that induction never stops in cases (a) and (b). In case (c), the full sequence $t_n$ is well defined and finite, and for all $n\geq 1$ induction does not stop at time $t_n$. Therefore $S=\sup_n t_n$ and again, induction never stops.
\end{proof}


\paragraph{The functions $\Psi_1$ and $\Upsilon_i$ are homeomorphisms:}

The goal here is to prove the following lemma.

\begin{lem}\label{lemma:homeo}
\begin{description}
    \item[\textmd{(i)}] The function $\Psi_1$ is a homeomorphism on the set $A_1\setminus G_1$. 
    \item[\textmd{(ii)}] For all $i\in\{1,2\}$, $\Upsilon_i$ is a homeomorphisms on the set $C_i\setminus \Upsilon_i^{-1}(\R_+\times G_0)$.
    \item[\textmd{(iii)}] {The set $\Psi^{-1}_1(\R_+\times G_0)$ is nowhere dense.}
    \item[\textmd{(iv)}] {The set $\Upsilon^{-1}_i(\R_+\times G_0)$ is nowhere dense for $i=1,2$.}
\end{description}
\end{lem}


\begin{proof}
    The proof relies on the fact that all the functions have explicit inverse on the sets given in the statement. In particular, the proofs {of (i) and (ii)} reduce to checking the continuity properties of $\Psi_1$ and $\Upsilon_i$ and their explicit inverses. Since the proofs of (i) and (ii) are similar, we only prove (i) and leave the proof of (ii) to the reader. 

    Recall that for all $(\beta_2,e)\in A_1$, $\Psi_1(\beta_2,e)=(t_1(\beta_2,e),\phi_1(t_1(\beta_2,e),e))$. Let $(t,e')$ be in the image of $A_1\setminus\hat{A}_1$ by the map $\Psi_1$. Then there exists $(\beta_2,e)\in A_1\setminus\hat{A}_1$ such that $t_1(\beta_2,e)=t$, $\phi_1(t,e)=e'$,
    \[
    \beta_2+\int_0^{t}f_2(\phi_1(s,e))ds=1
    \]
    and for all $s<t$
    \[
    0<\beta_2+\int_0^{s}f_2(\phi_1(s,e))ds<1.
    \]    
    We introduce $\phi^{-}_1$ the flow of the dynamical system $\dot{y}=-F_1(y)$, i.e.\ $\phi^{-}_1(t,y)$ is the value at time $t$ of the solution of this dynamical system started at time 0 from $y$. Then
    \[
    e=\phi^{-}_1(t,e')\quad\text{and}\quad \beta_2=1-\int_0^t f_2(\phi^{-}_1(s,e'))ds.
    \]
    Therefore, $\Psi_1$ is one-to-one from $A_1\setminus\hat{A}_1$ to its image and for all $(t,e')\in \text{Im}(A_1\setminus\hat{A}_1)$
    \[
    \Psi_1^{-1}(t,e')=\left(1-\int_0^t f_2(\phi^{-}_1(s,e'))ds \, ,\, \phi^{-}_1(t,e')\right),
    \]
    and, similarly, $\Psi_1$ is one-to-one from $\hat{A}_1$ to its image and for all $(t,e')\in \text{Im}(\hat{A}_1)$
    \[
    \Psi_1^{-1}(t,e')=\left(-\int_0^t f_2(\phi^{-}_1(s,e'))ds \, ,\, \phi^{-}_1(t,e')\right).
    \]
    Now, for $(\beta_2,e)\in A_1\setminus G_1$, we have $f_2(\phi_1(t_1(\beta_2,e),e))\not=0$. In addition, if $(\beta_2,e)\in A_1\setminus(\hat{A}_1\cup G_1)$, we must have $f_2(\phi_1(t_1(\beta_2,e),e))>0$, and conversely if $(\beta_1,e)\in \hat{A}_1\setminus G_1$. Therefore, on the image of $A_1\setminus G_1$, the inverse function of $\Psi_1$ is given by
    \[
    \Psi_1^{-1}(t,e')=\begin{cases}
        \left(1-\int_0^t f_2(\phi^{-}_1(s,e'))ds \, ,\, \phi^{-}_1(t,e')\right) & \text{if }f_2(e')>0, \\
        \left(-\int_0^t f_2(\phi^{-}_1(s,e'))ds \, ,\, \phi^{-}_1(t,e')\right) & \text{if }f_2(e')<0.
    \end{cases}  
    \]
    In view of the previous expression for $\Psi_1^{-1}$, this function is clearly continuous on the image of $A_1\setminus G_1$.

    Therefore, to conclude, it only remains to check that $\Psi_1$ is continuous on $A_1\setminus G_1$.
    This obviously follows from the next lemma, whose proof is given after the present proof.
    \begin{lem}
    \label{lem:t_1-cont}
        The function $(\beta_2,e)\mapsto t_1(\beta_2,e)$ is continuous on $A_1\setminus G_1$. Also, the function $e \mapsto s_1(e)$ is continuous on $\{e \in \mathbb{R}^q: (0,e)\in C_1\text{ and }f_1(\phi_2(s_1(e), e)> 0\}\cup\{e \in \mathbb{R}^q: (0,e)\in C_2\text{ and }f_2(\phi_1(s_1(e), e)> 0\}$.
    \end{lem}
    This ends the proof of (i).

    Let us now prove (iii). Assume by contradiction that there exists an open set $V$ in the adherence of $\Psi^{-1}_1(\R_+\times G_0)$ and consider $(\beta_2,e)\in \Psi^{-1}_1(\R_+\times G_0)$ and $r>0$ such that $B_r:=B((\beta_2,e),r)\subset V$. For simplicity, we shall denote below by $t_1$ the time $t_1(\beta_2,e)$. Observe that the flow $\phi_1$ is Lipschitz, so there exists a constant $C$ such that, for all $e'\in \R^\dimE$,
    \[
    \sup_{s\in[0,t_1]}|\phi(s,e')-\phi(s,e)|\leq C|e'-e|.
    \]
    Either 
    \[
    \beta_2+\int_0^{t_1}f_{2}(\phi_1(s,e))ds=1\text{ and } \forall t\in[0,t_1],\ 0<\beta_2+\int_0^{t}f_{2}(\phi_1(s,e))ds<1
    \]
    or
    \[
    \beta_2+\int_0^{t_1}f_{2}(\phi_1(s,e))ds=0\text{ and } \forall t\in[0,t_1],\ 0<\beta_2+\int_0^{t}f_{2}(\phi_1(s,e))ds<1.
    \]
    Let us assume we are in the first case. The other case can be dealt with similarly.
    Reducing $r>0$ if necessary, we deduce that for all $(\beta_2',e')\in B'_r$, where
    \[
    B'_r:=\left\{(\beta'_2,e')\in B_r:\beta'_2\in(\beta_2,1)\text{ and }|e'-e|\leq C^{-1}t_1^{-1}\|f_2\|_\text{Lip}^{-1}(\beta'_2-\beta_2)\right\},
    \]
    we have
    \[
    \beta'_2+\int_0^{t_1}f_{2}(\phi_1(s,e'))ds\geq 1
    \]
    and for all $t\in[0,t_1]$
    \[
    \beta'_2+\int_0^{t}f_{2}(\phi_1(s,e'))ds>0.
    \]
    In particular, $t_1(\beta_2',e')\leq t_1$. Since $\beta_2+\int_0^{t}f_{2}(\phi_1(s,e))ds<1$ for all $t<t_1$, it is easy to deduce that
    \[
    \lim_{(\beta_2',e')\to (\beta_2,e),\ (\beta'_2,e')\in B'_r}t_1(\beta_2',e')=t_1(\beta_2,e).
    \]
    Therefore, $t_1$ is continuous on
    \[
    B''_r:=\left\{(\beta'_2,e')\in B_r:\beta'_2\in[\beta_2,1)\text{ and }|e'-e|\leq C^{-1}t_1^{-1}\|f_2\|_\text{Lip}^{-1}(\beta'_2-\beta_2)\right\}\subset V.
    \]
    Using the same argument as for point (i) above, this implies that $\Psi_1$ is continuous on $B_r''$ and one-to-one from $B_r''$ to its image, with continuous inverse that we denote $\widetilde{\Psi}_1^{-1}$ (to avoid confusion with the inverse defined in the proof of (i)). Since $\R_+\times G_0$ is nowhere dense, it then follows that $\widetilde{\Psi}_1^{-1}(\R_+\times G_0)$ is also nowhere dense. This implies that $\Psi_1^{-1}(\R_+\times G_0)\cap B''_r$ is nowhere dense. In particular, $B''_r$ is not included in the adherence of $\Psi_1^{-1}(\R_+\times G_0)$. Since $B''_r\subset V$, we obtain a contradiction by definition of $V$.  
    
   The proof of~(iv) follows the same general idea: let us fix $e\in\R^\dimE$ such that $(t,e)\in\Upsilon^{-1}_1(\R_+\times G_0)\subset C$ for some (hence all) $t\geq 0$. Our goal is to prove that in any neighborhood of $e\in\R^\dimE$, there exists a non-empty open ball $V$ such that $\{(0,e):e\in V\}\subset C$ 
   such that, for all $e'\in V$, $\Upsilon_1(0,e)\not\in \R_+\times G_0$ (and hence $\Upsilon_1(t,e)\not\in \R_+\times G_0$ for all $t\geq 0$). Let us assume that 
    \[
    1+\int_0^{s_1(e)} f_2(\phi_1(s,e)) ds=1\quad\text{and}\quad f_2(\phi_1(s_1(e),e)=0.
    \]
    The case $f_2(\phi_1(s_1(e),e)\neq 0$ (and hence $f_1(\phi_1(s_1(e),e)=0$) cannot occur because this would imply $f_2(\phi_1(s_1(e),e)> 0$ and by the invasion implies fixation assumption, this would implies that $f_1(\phi_1(s_1(e),e)<0$. The case where $1+\int_0^{s_1(e)} f_2(\phi_1(s,e)) ds=0$ can be handled similarly so we leave the details to the reader.

    Let us make the abuse of notation, valid until the end of this proof, of (re-)defining for all $t\in[0,s_1(e)]$ 
    \[
    \beta_2(t):=1+\int_0^t f_2(\phi_1(s,e)) ds.
    \]
    Let $\beta_*:=\inf_{t\in[0,s_1(e)]}\beta_2(t)$ and $t^*:=\sup\{t\in[0,s_1(e)]:\beta_2(t)=\beta_*\}$. For all $\eta\in[\beta^*,1]$, we define
    \[
    t_\ell(\eta):=\sup\left\{t\in[0,t_*]:\beta_2(t)=\eta\right\}\quad\text{and}\quad
    t_r(\eta):=\inf\left\{t\in[t_*,s_1(e)]:\beta_2(t)=\eta\right\}.
    \]
    Note that, by continuity of $\beta_2$ and because $\beta_2(t)<1$ for all $t\in(0,s_1(e))$, $t_\ell(\eta)\to 0$ and $t_r(\eta)\to s_1(e)$ when $\eta\to 1$. Finally, let us also define for all $\eta\in[\beta_*,1]$ $e_\eta:=\phi_1(t_\ell(\eta),e)$. Then $s_1(e_\eta)=t_r(\eta)-t_\ell(\eta)$ and $\phi_1(t_r(\eta),e)=\phi_1(t_r(\eta)-t_\ell(\eta),e_\eta)=\phi_1(s_1(e_\eta),e_\eta)$.

    We first prove that there exists an increasing sequence $(\eta_n)_{n\in\N}$ converging to 1 such that, denoting for short $e_{\eta_n}$ by $e_n$, then for all $n\geq 0$, $f_2(\phi_1(s_1(e_n),e_n))\neq 0$ (hence $f_2(\phi_1(s_1(e_n),e_n))>0$ since $t_r(\eta)$ is a first hitting time from below of a function whose derivative is $f_2(\phi_1(s,e_n))$). Let us proceed by contradiction: if this was not true, there would exist $\eta_0<1$ such that, for all $\eta\in[\eta_0,1)$, $f_2(\phi_1(s_1(e_\eta),e_\eta))=f_2(\phi_1(t_r(\eta),e))= 0$. Let us define for all $\eta\in[\beta_*,1]$
    \[
    t'_r(\eta):=\inf\{t>t_r(\eta):\beta_2(t)>\eta\}.
    \]
    Note that there exists an at most denumerable sequence $(\eta_n)_{n\in I}$ in $(t_*,1)$ (with $I$ finite or denumerable) such that $t_r(\eta_n)<t'_r(\eta_n)$ and for all $\eta\in[t_*,1]\setminus\{\eta_n\}_{n\in I}$, $t_r(\eta)=t'_r(\eta)$ (this is because the function $t_r$ and $t'_r$ are both increasing). Now,
    \begin{align*}
        \beta_2(s_1(e))-\beta_2(t'_r(\eta_0)) & =\int_{t'_r(\eta_0)}^{s_1(e)} f_2(\phi_1(s,e))ds \\ & =\int_{t'_r(\eta_0)}^{s_1(e)} f_2(\phi_1(s,e))\mathbbm{1}_{\R_+\setminus\left(\bigcup_{n\in I} [t_r(\eta_n),t'_r(\eta_n)]\right)}ds
    \end{align*}
    because for all $n\in I$
    \[
    \int_{t_r(\eta_n)}^{t'_r(\eta_n)}f_2(\phi_1(s,e)) ds=\beta_2(t'_r(\eta_n))-\beta_2(t_r(\eta_n))=0.
    \]
    Now, for all $s\not\in \bigcup_{n\in I} [t_r(\eta_n),t'_r(\eta_n)]$, $t_r(\beta_2(s))=t'_r(\beta_2(s))=s$ and hence $f_2(\phi_1(s,e))=f_2(\phi_1(t_r(\beta_2(s)),e))=0$. This implies that $\beta_2(t'_r(\eta_0))=\beta_2(s_1(e))$, which is a contradiction.

    Now, $e_n\to e$ when $n\to+\infty$, and since $f_2(\phi_1(s_1(e_n),e_n))>0$, by Lemma~\ref{lem:t_1-cont}, by Assumption~\ref{hyp:regularity-Ltka-Volterra} and by the continuity of the flow of $\phi_1$, we deduce that, for all $n\geq 0$ there exists a radius $r_n>0$ such that, for all $e'$ in the ball of radius $r_n$ centered at $e_n$, $s_1(e')<+\infty$ and $f_2(\phi_1(s_1(e'),e'))> 0$. By the invasion implies fixation assumption, this implies that $f_1(\phi_1(s_1(e'),e'))< 0$, hence $\Upsilon_1(0,e')\not\in\R_+\times G_0$. This implies that $\Upsilon_1^{-1}(\R_+\times G_0)$ is nowhere dense and Point~(iv) is proved. 
\end{proof}

\begin{proof}[Proof of Lemma~\ref{lem:t_1-cont}]
    Below, we give the proof of the continuity of $(\beta_2,e)\mapsto t_1(\beta_2,e)$. The continuity of $e\mapsto s_1(e)$ can be proved similarly (actually, with simpler arguments). 
    
    Let $(\beta_2,e)\in A_1 \setminus(\hat{A}_1\cup G_1)$ and set $T=t_1(\beta_2,e)+1$. For all $t\geq 0$, $\beta'_2\in[0,1]$ and $e'\in\R^\dimE$, we define
    \[
    \hat{\beta}_2(t,\beta'_2,e')=\beta'_2+\int_0^t f_2(\phi_1(s,e'))ds.
    \]
    Note that $\hat{\beta}_2(t_1(\beta_2,e),\beta_2,e)=1$ and, since $(\beta_2,e)\not\in G_1$, $f_2(\phi_1(t_1(\beta_2,e), e)>0$. This implies, that, for $a:=f_2(\phi_1(t_1(\beta_2,e), e)/2$, there exists $\delta_0>0$ such that, 
    \[
    \forall t \in[t_1(\beta_2,e)-\delta_0,t_1(\beta_2,e)], \quad \hat{\beta}_2(t,\beta_2,e)\leq 1-a(t_1(\beta_2,e)-t)
    \]
    and 
    \[
    \forall t \in[t_1(\beta_2,e),t_1(\beta_2,e)+\delta_0], \quad \hat{\beta}_2(t,\beta_2,e)\geq 1+a(t-t_1(\beta_2,e)).
    \]
    Now, the bounded Lipschitz assumptions for $F_1$ and $f_2$ and usual application of Gronwall's lemma implies the continuity of the flow $\phi_1$ and thus of $\hat{\beta}_2$: there exists a constant $C(T)$ such that, for all $t\in[0,T]$, $e',e''\in\R^\dimE$ and $\beta'_2,\beta''_2\in[0,1]$,
    \[
    \left|\phi_1(t,e')-\phi_1(t,e'')\right|\leq C(T)|e'-e''|
    \]
    and
    \[
    \left|\hat{\beta}_2(t,\beta'_2,e')-\hat{\beta}_2(t,\beta''_2,e'')\right|\leq C(T)(|\beta'_2-\beta''_2|+|e'-e''|).
    \]
    This implies that
    \[
    \forall t \in[t_1(\beta_2,e)-\delta_0,t_1(\beta_2,e)], \quad \hat{\beta}_2(t,\beta'_2,e')\leq 1-a(t_1(\beta_2,e)-t)+C(T)(|\beta'_2-\beta_2|+|e'-e|)
    \]
    and 
    \[
    \forall t \in[t_1(\beta_2,e),t_1(\beta_2,e)+\delta_0], \quad \hat{\beta}_2(t,\beta'_2,e')\geq 1+a(t-t_1(\beta_2,e))-C(T)(|\beta'_2-\beta_2|+|e'-e|).
    \]
    Now, $\sup_{t\in[0,t_1(\beta_2,e)-\delta_0]}\hat{\beta}_2(t,\beta_2,e)<1$, so there exists $\eta>0$ such that, for all $(\beta'_2,e')$ in the ball centered at $(\beta_2,e)$ with radius $\eta$, $\sup_{t\in[0,t_1(\beta_2,e)-\delta_0]}\hat{\beta}_2(t,\beta'_2,e')<1$ and hence $t_1(\beta'_2,e')>t_1(\beta_2,e)-\delta_0$. Setting $\delta:=(a\delta_0/C(T))\wedge\eta$, we deduce that, for all $(\beta'_2,e')$ in the ball centered at $(\beta_2,e)$ with radius $\delta$,
    \[
    t_1(\beta'_2,e')\in\left[t_1(\beta_2,e)-\frac{C(T)}{a}(|\beta'_2-\beta_2|+|e'-e|)\, ,\, t_1(\beta_2,e)+\frac{C(T)}{a}(|\beta'_2-\beta_2|+|e'-e|)\right].
    \]
    The proof is similar in the case $(\beta_2,e)\in \hat{A}_1 \setminus G_1$. Hence Lemma~\ref{lem:t_1-cont} is proved. 
\end{proof}

\paragraph{Conclusion of the proof of Theorem \ref{thm:Baire-d=2-k=1}:}
Let $D_0 := \bigcup_{n \ge 1} G_n$ be the set of initial conditions for which induction might stop. By Lemma~\ref{lemma:homeo}, each $G_n$ is a Baire meager set in $(0, 1) \times \mathbb{R}^\dimE$. The Baire Category Theorem states that a countable union of meager sets is itself a meager set. Therefore, $D_0$ is a Baire meager set. 

Consequently, by Proposition~\ref{prop:F_k}, for all initial conditions $(\beta_1(0), \beta_2(0), e(0)) \in D \setminus D_0$, induction never stops in Theorem~\ref{thm:BDP-max-2}. This completes the proof of Theorem~\ref{thm:Baire-d=2-k=1}.

\subsection{General case: proof of Theorem~\ref{thm:Baire-d>2}}
\label{sec:Baire-d>2}

The goal of this section is to prove Theorem~\ref{thm:Baire-d>2}. We already proved this result for $d=2$. The proof follows a similar argument as for Theorem~\ref{thm:Baire-d=2-k=1}, although combinatorially more intricate. Since most of the arguments are the same, we only explain below  how to modify the definitions of the sets $G_k$ and how to deal with the possible extinction of species.

We first introduce some notation:
\begin{enumerate}
    \item Let $I\subset\{1,2,\ldots,d\}$ be fixed; for all $i\in I$ and $A\subset [0,1]$, we denote by $A^{I\setminus\{i\}}$ the set of vectors $\beta=(\beta_\ell)_{\ell\in I\setminus\{i\}}$ 
    with $\beta_\ell\in A$ for all $\ell\in I\setminus\{i\}$; for $i\neq j$, we define similarly $A^{I\setminus\{i,j\}}$.
    \item For all $i\in I$, we denote by $\phi_i$ the flow associated to the vector field $F_i$.
    \item For all $i\in I$, $\beta\in(0,1)^{I\setminus\{i\}}$ and $e\in\R^\dimE$, we define
    \[
    t_1^i(\beta,e):=\inf\left\{t\geq 0:\exists\ell\neq i,\ \beta_\ell+\int_0^tf_{\ell, i}(\phi_i(s,e))ds\in\{0,1\}\right\}.
    \]
    \item For all $i\neq j$ in $I$, let
    \[
    A_{ij}:=\left\{(\beta,e)\in(0,1)^{I\setminus\{i\}}\times\R^\dimE: t_1^i(\beta,e)<\infty\text{ and }\beta_j+\int_0^{t^i_1(\beta,e)}f_{j, i}(\phi_i(s,e))ds\in\{0,1\}\right\}
    \]
    and
    \[
    \hat{A}_{ij}=\left\{(\beta,e)\in A_{ij}:\beta_j+\int_0^{t^i_1(\beta,e)}f_{j,i}(\phi_i(s,e))ds=0\right\}.
    \]
    \item For all $i\neq j$ in $I$, let $\Psi_{i,j}:A_{ij}\to \R_+\times[0,1]^{I\setminus\{i,j\}}\times\R^\dimE$ be defined for all $(\beta,e)\in A_{ij}$ by
    \[
    \Psi_{i,j}(\beta,e)=\left(t_1^i(\beta,e),\left(\beta_\ell+\int_0^{t_1^i(\beta,e)}f_{\ell, i}(\phi_i(s,e))ds\right)_{\ell\in I\setminus\{i,j\}}, \phi_i\big(t_1^i(\beta,e),e\big)\right).
    \]
    \item For all $i\neq j$ in $I$, $\beta\in(0,1)^{I\setminus\{i,j\}}$ and $e\in\R^\dimE$ such that $f_{j,i}(e)>0$ (and hence $f_{i,j}(e)<0$ since invasion implies fixation), let
    \[
    s_1^j(\beta,e):=\inf\left\{t> 0:\exists\ell\neq i,j,\ \beta_\ell+\int_0^tf_{\ell, j}(\phi_j(s,e))ds\in\{0,1\}\text{ or }1+\int_0^tf_{i,j}(\phi_j(s,e))ds\in\{0,1\}\right\}.
    \]
    Note that $s_1(\beta,e)>0$ since $f_{i,j}(e)<0$.
    \item For all $i,j,p$ in $I$ such that $i\neq j$ and $j\neq p$ (note that $p=i$ is possible), let
    \begin{multline*}
    C_{ijp}:=\left\{(\beta,e)\in(0,1)^{I\setminus\{i,j\}}\times\R^\dimE: f_{j,i}(e)>0,\ s_1^j(\beta,e)<\infty\right.
    \\
    \left.\text{ and }\beta_p+\int_0^{s^j_1(\beta,e)}f_{p,j}(\phi_j(s,e))ds\in\{0,1\}\right\},
    \end{multline*}
    with the convention that $\beta_p=1$ if $p=i$, and
    \[
    \hat{C}_{ijp}=\left\{(\beta,e)\in C_{ijp}:\beta_p+\int_0^{s^j_1(\beta,e)}f_{p,j}(\phi_j(s,e))ds=0\right\}.
    \]
    \item For all $i,j,p$ in $I$ such that $i\neq j$ and $j\neq p$, let $\Upsilon_{ijp}:\R_+\times C_{ijp}\to \R_+\times[0,1]^{I\setminus\{j,p\}}\times\R^\dimE$ be defined for all $t\geq 0$, $(\beta,e)\in C_{ijp}$ by
    \[
    \Upsilon_{ijp}(t,\beta,e)=\left(t+s_1^j(\beta,e),\left(\beta_\ell+\int_0^{s_1^j(\beta,e)}f_{\ell, j}(\phi_j(s,e))ds\right)_{\ell\in I\setminus\{j,p\}}, \phi_j\big(s_1^j(\beta,e),e\big)\right).
    \]
\end{enumerate}

The sets $G_n$ of initial conditions such that induction stops at step $n$ defined in Section~\ref{sec:Baire-d=2} need to be modified as follows. First, we need to keep track of the set of species that are initially alive in the population. Indeed, once a species goes extinct, this set changes but contrary to the case $d=2$, we are not sure that induction never stops. Therefore, we need to proceed by induction on the number of surviving species. This is why we consider a general subset $I$ of $\{1,2,\ldots,d\}$ in the notations above.

Below, we construct for all $I\subset\{1,2,\ldots,d\}$, $i
\in I$ and all $n\geq 1$ the set $G^I_{i,n}$ of initial conditions with set $I$ of initially non-extinct species, with species $i$ as initial dominant species and such that induction might stop at step $n$. For this, we proceed by induction on $|I|$. To start the induction, we use Theorem~\ref{thm:Baire-d=2-k=1} from Section~\ref{sec:Baire-d=2}: in the proof of this theorem, we constructed for all $i\neq j\in\{1,\ldots,d\}$ sets of initial conditions $G^{\{i,j\}}_{i,n}$ with only species $i$ and $j$ initially non-extinct and species $i$ initially dominant for which induction may stop at step $n$ and which is Baire meager. In addition, we define
\begin{equation*}
    G^{\{i,j\}}_i=\bigcup_{n\geq 1}G^{\{i,j\}}_{n,i}.
\end{equation*}
For the induction step, given $I\subset\{1,\ldots, d\}$ such that $|I|\geq 3$, we assume that we have constructed $G^J_{j,n}$ for all $J\subsetneq I$, all $j\in J$ and all $n \geq 1$ and we define
may stop at step $n$ and which is Baire meager. In addition, we define
\begin{equation*}
    G^{J}_j=\bigcup_{n\geq 1}G^{J}_{n,j}.
\end{equation*}
Then the sets $G^I_{i,n}$ are constructed as follows: for all $i\neq j\in I$, let
\begin{equation}
    \label{eq:def-G_0^I}
    G^{I,i,j}_0:=\left\{(\beta,e)\in[0,1]^{I\setminus\{i,j\}}\times\R^\dimE:f_{j,i}(e)=0\text{ or }\beta\in\partial([0,1]^{I\setminus\{i,j\}})\right\},
\end{equation}
where $\partial([0,1]^{I\setminus\{i,j\}})$ is the boundary of the set $[0,1]^{I\setminus\{i,j\}}$, i.e.\ the sets of $(\beta_\ell)_{\ell\in I\setminus\{i,j\}}$ such that $\beta_\ell\in\{0,1\}$ for at least one $\ell\in I\setminus\{i,j\}$. Then, we define for all $i\in I$
\begin{equation}
    \label{eq:G^I_1}
    G^I_{1,i}:=\bigcup_{j\in I\setminus\{i\}}\Psi_{ij}^{-1}\left(\R_+\times G^{I,i,j}_0\right).
\end{equation}
With this definition, for all initial conditions $(\beta,e)\in(0,1]^I\times\R^\dimE$ with $I$ as set of initially non-extinct species of the form such that $\beta_i=1$ and $((\beta_\ell)_{\ell\in I\setminus\{i\}},e)\in G^I_{1,i}$, at time $t_1$ in the statement of Theorem~\ref{thm:BDP-max-2}, either the fitness of the invading or vanishing species is zero or more than two non-resident species (i.e.\ species different from $i$) hit 0 or 1 at the same time. In particular, the initial conditions such that species $i$ is initially resident and induction stops in Theorem~\ref{thm:BDP-max-2} at time $t_1$ belong to this set of initial conditions.

With some abuse of notation, we will call $((\beta_\ell)_{\ell\in I\setminus\{i\}},e)\in ((0,1]^{I\setminus\{i\}}\times\R^\dimE)\setminus G^I_{1,i}$ the initial conditions for which species $i$ is initially resident and induction does not stop at time $t_1$ (instead of the complete vector $(\beta,e)$).

Similarly, we define the set of initial conditions for which species $i$ is initially resident and induction may stop at time $t_2$ in Theorem~\ref{thm:BDP-max-2} as follows:
\begin{equation}
    \label{eq:G^I_2}
    G_{2,i}^I:=\left[\left(\bigcup_{j\in I\setminus\{i\}}\bigcup_{p\in I\setminus\{j\}}\Psi^{-1}_{ij}\left(\Upsilon^{-1}_{ijp}\left(\R_+\times G^{I,j,p}_0\right)\right)\right)\bigcup\left(\bigcup_{j\in I\setminus\{i\}}\Psi^{-1}_{ij}\left(\R_+\times G^{I\setminus\{j\}}_i\right)\cap \hat{A}_{i,j}\right)\right]\setminus G^I_{1,i}.
\end{equation}
This set is the union of two sets: the first one is the set of initial conditions for which induction does not stop at time $t_1$ and there is a change of resident species from $i$ to $j$ at time $t_1$, but at time $t_2$, either one of the fitnesses vanish or more than two non-resident species hit 0 or 1 at the same time; the second one is the set of initial conditions for which induction does not stop at time $t_1$ and species $j$ goes extinct at time $t_1$ and the population state at time $t_1$ with species $j$ removed belongs to the set $G^{I\setminus\{j\}}_i$ of initial conditions for which induction stops with initial surviving species $I\setminus\{j\}$ and initial resident species $i$ (see~\eqref{eq:def-G^I_i} below for the definition of $G^{I\setminus\{j\}}$).

Similarly, for all $n\geq 2$ and $i\in I$, we define
\begin{multline*}
    G_{n+1,i}^I:=\left[\left(\bigcup_{\substack{(i_0,\ldots,i_{n+1})\in I^{n+2},\ i_0=i, \\ \text{and } i_{\ell+1}\neq i_\ell\ \forall 0\leq \ell\leq n}}\Psi^{-1}_{i_0i_1}\left(\Upsilon^{-1}_{i_0i_1i_2}\left(\ldots\Upsilon^{-1}_{i_{n-1}i_n i_{n+1}}\left(\R_+\times G^{I,i_n,i_{n+1}}_0\right)\ldots\right)\right)\right)\bigcup\right. 
    \\ \left.
    \left(\bigcup_{\substack{(i_0,\ldots,i_n)\in I^{n+1},\ i_0=i, \\ \text{and }i_{\ell+1}\neq i_\ell\ \forall 0\leq \ell\leq n-1}}\Psi^{-1}_{i_0i_1}\left(\Upsilon^{-1}_{i_0i_1i_2}\left(\ldots\left(\Upsilon^{-1}_{i_{n-2}i_{n-1}i_{n}}\left(\R_+\times G_i^{I\setminus\{i_n\}}\right)\cap\hat{C}_{i_{n-2}i_{n-1}i_{n}}\right)\ldots\right)\right)\right)\right]\setminus (G^I_{1,i}\cup\ldots\cup G^I_{n,i}). 
\end{multline*}
Finally, we define for all $i\in I$
\begin{equation}
    \label{eq:def-G^I_i}
    G^I_i=\bigcup_{n\geq 1}G^I_{n,i}.
\end{equation}
This ends the inductive construction of the sets $G^I_{i,n}$.

\medskip

In view of the previous discussion, Proposition~\ref{prop:F_k} can be extended to the above sets as follows. We leave the proof to the reader.

\begin{prop}
    \label{prop:G_k}
    For all $(\beta,e)\in D$ which does not belong to the set
    \[
    \bigcup_{i\in I}\left\{(\beta,e)\in(0,1]^I\times\R^\dimE:((\beta_\ell)_{\ell\in I\setminus\{i\}},e)\in G^I_i\text{ and }\beta_i=1\right\},
    \]
    induction never stops in Theorem~\ref{thm:BDP-max-2}.
\end{prop}

From Assumption~\ref{hyp:not-many-zeros} of Theorem~\ref{thm:Baire-d>2}, the sets $G^{I,i,j}_0$ are closed with empty interior, hence are nowhere dense. Hence, as in Section~\ref{sec:Baire-d=2}, the key to prove that the sets $G^I_i$ are Baire meager is the following lemma. 


\begin{lem}\label{lemma:homeo-d>2}
\begin{description}
    \item[\textmd{(i)}] For all $i\neq j$ in $I$, $\Psi_{ij}$ is a homeomorphism on the set 
    \[
    A_{ij}\setminus \Psi^{-1}_{ij}\left(\R_+\times G^{I,i,j}_0 
    \right).
    \]
    \item[\textmd{(ii)}] For all $i,j,p$ in $I$ such that $i\neq j$ and $j\neq p$, $\Upsilon_{ijp}$ is a homeomorphisms on the set 
    \[
    \R_+\times C_{ijp}\setminus \Upsilon^{-1}_{ijp}\left(\R_+\times G^{I,j,p}_0 
    \right).
    \]
    \item[\textmd{(iii)}] {For all $I\subset\{1,\ldots,d\}$ with $|I|\geq 2$ and $i\neq j\in I$, the set $\Psi_{ij}^{-1}(\R_+\times G_0^{I,i,j})$ is nowhere dense.}
    \item[\textmd{(iv)}] {For all $I\subset\{1,\ldots,d\}$ with $|I|\geq 2$, $i\neq j\in I$ and $p\neq j\in I$, the set $\Upsilon_{ijp}^{-1}(\R_+\times G_0^{I,j,p})$ is nowhere dense.} 
\end{description}
\end{lem}

This result follows from similar arguments as Lemma~\ref{lemma:homeo}. We do not detail the proof, but we only give the main ideas. The first two points can be proved exactly as Lemma~\ref{lemma:homeo}~(i) and~(ii). {To prove the third point, we first show that $\Psi_{ij}^{-1}(\R_+\times \{(\beta,e)\in (0,1)^{I\setminus \{i,j\}}\times \R^\dimE:f_{j,i}(e)=0\})$ is nowhere dense can be proved similarly as Lemma~\ref{lemma:homeo}~(iii). The fact that $\Psi_{ij}^{-1}(\R_+\times \partial([0,1]^{I\setminus\{i,j\}})\times \R^\dimE)$ is nowhere dense can be proved following similar ideas as Lemma~\ref{lemma:homeo}~(iii): if $\ell\in I\setminus\{i,j\}$ is such that the $\beta_\ell$ component of $\Psi_{ij}(\beta,e)$ is equal to 1 (resp. 0), we can consider as in Lemma~\ref{lemma:homeo}~(iii) an open set in the neighborhood of $(\beta,e)$ with smaller (resp. larger) values at the $\ell$-th component of $\beta$. Similarly, the fact that $\Upsilon_{ijp}^{-1}(\R_+\times G_0^{I,j,p})$ is nowhere dense can be proved combining ideas from Lemma~\ref{lemma:homeo}~(iii) and~(iv) as follows: we first prove that $\Upsilon_{ijp}^{-1}(\R_+\times \{(\beta,e)\in (0,1)^{I\setminus \{j,p\}}\times \R^\dimE:f_{p,j}(e)=0\})$ is nowhere dense using either the method of Lemma~\ref{lemma:homeo}~(iii) if $p\neq i$ or the method of Lemma~\ref{lemma:homeo}~(iv) if $p=i$. The fact that $\Upsilon_{ijp}^{-1}(\R_+\times \partial([0,1]^{I\setminus\{i,j\}})\times \R^\dimE)$ is nowhere dense can be proved as above.}

The conclusion of the proof of Theorem~\ref{thm:Baire-d>2} then proceeds as follows: in view of~\eqref{eq:G^I_1} and by the previous lemma, $G^I_{1,i}$ is a finite union of nowhere dense sets, hence is also nowhere dense. Assuming by induction that $G^J_j$ is a Baire meager set for all $J\subsetneq I$ and all $j\in J$ (note that this is true for $|J|=2$ by Theorem~\ref{thm:Baire-d=2-k=1}), we deduce from~\eqref{eq:G^I_2} and the previous lemma that $G^I_{2,i}$ is a finite union of Baire meager sets, hence is also Baire meager. Similarly, $G^I_{n,i}$ is Baire meager for all $n\geq 2$ and hence $G^I_i$ is also Baire meager.

Hence we have proved by induction that $G^I_i$ is Baire meager for all $I\subset\{1,2,\ldots,d\}$ and thus Theorem~\ref{thm:Baire-d>2} follows from Proposition~\ref{prop:G_k}.

\subsection{Conditions ensuring that $S=+\infty$ when induction never stops: proof of Proposition~\ref{prop:S-infty}}
\label{sec:S-infty}

We gave in the last subsections conditions ensuring that induction never stops, but this not necessarily implies that $S=+\infty$. We
now study under which conditions this last property holds true. We start with the simplest case $d=2$. The next proposition
corresponds to Proposition~\ref{prop:S-infty} when $d=2$. The generalization to higher dimension is given at the end of the section.

\begin{prop}
\label{prop:S=infty-d=2}
Assume that $d=2$, $\dimE\geq 1$, that the assumptions of Theorem~\ref{thm:BDP-max-2} are satisfied, that induction never stops for a given initial condition $(\beta(0),e(0))$ and that Assumption~\ref{hyp:no-heteroclinic-loop} holds true. 
Then $S=+\infty$. In particular, the convergence in Theorem~\ref{thm:BDP-max-2} holds true in probability in $L^\infty_\text{loc}([0,+\infty),\mathbb{R}^{\dimE+2})$.
\end{prop}


\begin{proof}
    Recall that, in the case $d=2$, Assumption~\ref{hyp:no-heteroclinic-loop} is equivalent to assuming that the functions $f_{2,1}$ and $f_{1,2}$ have no common zeroes.

    Without loss of generality, assume $(\beta(0),e(0))\in D_1$. By contradiction, if induction never stops and $S<+\infty$, $t_k<+\infty$ for all $k\geq 0$ and hence $\beta_1(t_k)=\beta_2(t_k)=1$ for all $k\geq 1$, i.e.\ each time $t_k$ is a time of change of resident population. Since $d=2$, necessarily, the two species 1 and 2 are successively resident, which implies that $f_{1,2}(e(t_{2p}))>0$, $f_{2,1}(e(t_{2p}))<0$, $f_{1,2}(e(t_{2p+1}))<0$ and $f_{2,1}(e(t_{2p+1}))>0$ for all $p\geq 1$.

    Now, $e(t)$ follows the flow associated to the vector field $F_1$ on $[t_{2p},t_{2p+1})$ and $F_2$ on $[t_{2p+1},t_{2p+2})$ for all $p\geq 0$. Since we assumed $F_1$ and $F_2$ bounded this means that the function $e(t)$ is Lipschitz, hence $e^\star:=\lim_{t\to S}e(t)$ is well-defined as we assumed $S<+\infty$. By continuity of $f_{1,2}$ and $f_{2,1}$, we deduce that $f_{1,2}(e^\star)=f_{2,1}(e^\star)=0$, which contradicts our assumption.
\end{proof}



The case $d\geq 3$ is more intricate because we cannot expect that $S=+\infty$ generically. We first give an example for $d=3$ with a heteroclinic cycle for which induction never stops and $S<\infty$. 

\paragraph{A simple example where $d=3$ and $S<\infty$.}


We assume $d=3$, that the rates are defined by
\begin{gather}
    b_i(e) = 2, \qquad d_i(e) = 0, \qquad c_{i,i}(e)=1 \qquad  \forall i\in\{1,2,3\} \label{eq.simple.ex.Sfinite.1}\\
    c_{2,1}(e)=c_{3,2}(e)=c_{1,3}(e)=\frac{1}{2},\label{eq.simple.ex.Sfinite.2}\\
    c_{1,2}(e)=c_{2,3}(e)=c_{3,1}(e)=\frac{5}{4}\label{eq.simple.ex.Sfinite.3}
\end{gather}
and that $F=0$, i.e.\ the environment satisfies $e(t)=e(0)$.
 
\begin{figure}
\centering
a\includegraphics[width=0.5\textwidth]{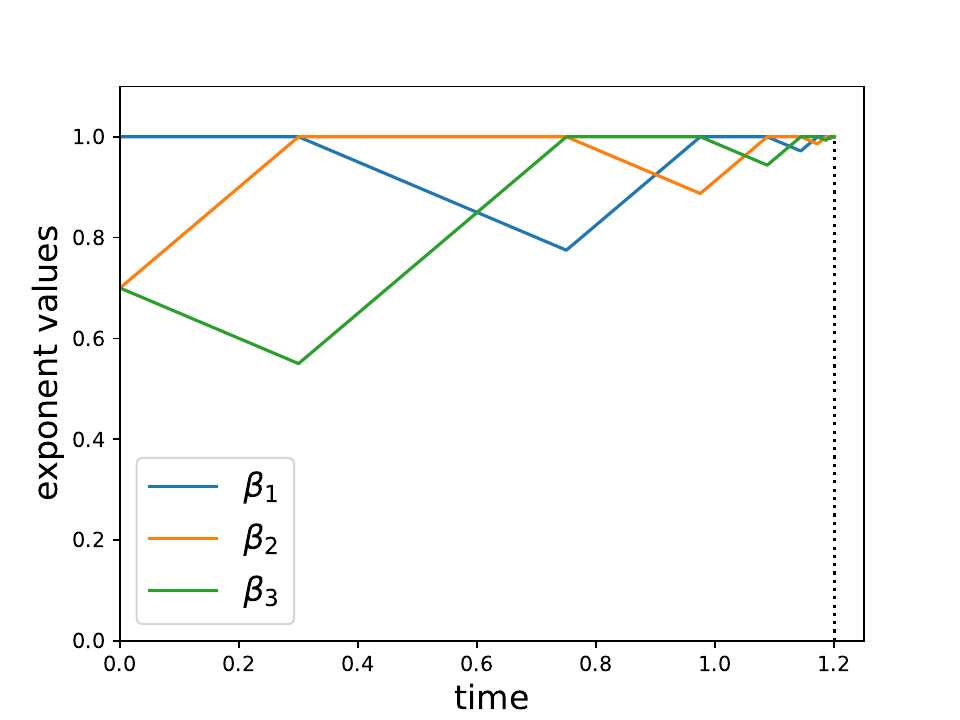}
\caption{\label{fig:simple.ex.Sfinite}Example of the asymptotical exponent dynamics for the model \eqref{eq.simple.ex.Sfinite.1}-\eqref{eq.simple.ex.Sfinite.3}, starting from the initial condition $\beta_0(1,1-0.3,1-0.3)$.}
\end{figure}

The fitnesses are then given by
\begin{gather*}
    f_{2,1}(e)=f_{3,2}(e)=f_{1,3}(e)=1, \qquad \ f_{1,2}(e)=f_{2,3}(e)=f_{3,1}(e)=-\frac{1}{2},
\end{gather*}
so that species 3 dominates species 2 (in the sense that species 3 invades species 2 and species 2 does not invade species 3), which dominates species 1, which dominates species 3. Moreover, notice that the dynamics of the species does not depend on the environmental variable $e(t)$.
Starting from the initial condition $\beta(0)=(1,1-\delta,1-\delta)$, with $\delta\in\left]0,\frac{2}{3}\right[$, the reader can check that $t_1=\delta$, the minimal exponent at the time $t_k$ (corresponding to the exponent of the invading species for the next step) is $\min_{i\in\{1,2,3\}}\beta_i(t_k)=1-\frac{3\delta}{2^k}$ and the time intervals, for $k\geq 1$, are $t_{k+1}-t_k=\frac{3\delta}{2^k}$. We then obtain $S=\sup_k t_k = 4\delta<+\infty$ (see Figure~\ref{fig:simple.ex.Sfinite}).

It is clear that this simple example is not specific and a similar result can be proved with a varying environment converging to a limit $e^\star$ such that 
\begin{gather}
    f_{2,1}(e^\star)>0,\ f_{1,2}(e^\star)<0,\ f_{3,2}(e^\star)>0,\ f_{2,3}(e^\star)<0, f_{1,3}(e^\star)>0,\ f_{3,1}(e^\star)<0, \label{eq:ex-heteroclinic} \\
    f_{2,1}(e^\star)>|f_{1,2}(e^\star)|,\ f_{3,2}(e^\star)>|f_{2,3}(e^\star)|\ \text{and}\ f_{1,3}(e^\star)>|f_{3,1}(e^\star)|. \notag
\end{gather}

\medskip

This example shows that for values of the model parameters satisfying the conditions~\eqref{eq:ex-heteroclinic}, we cannot expect in
general to have $S=+\infty$. This explains why Assumption~\ref{hyp:no-heteroclinic-loop} is needed in Proposition~\ref{prop:S-infty}.
We conclude this section by completing the proof of Proposition~\ref{prop:S-infty}.



\begin{proof}[Proof of Proposition~\ref{prop:S-infty}]
  The case $d=2$ is Proposition~\ref{prop:S=infty-d=2}. So let us assume $d\geq 3$.
  
  Assume induction never stops and, by contradiction, $S<+\infty$. This means that $t_k$ is well-defined for all $k\geq 1$ and hence
  there exists a sequence $(j_k)_k$ in $\{1,2,\ldots,d\}$ such that $\beta_{j_k}(t)=1$ for all $t\in[t_{k-1},t_k]$. Therefore, there
  exists a subset $I\neq\emptyset$ of $\{1,2\ldots,d\}$ such that $j_k=i$ for infinitely many $k$ if and only if $i\in I$. In
  addition, for all $i\in I$, there exists $\ell(i)\in I\setminus\{i\}$ such that $j_k=i$ and $j_{k+1}=\ell(i)$ for infinitely many
  $k$. Note that this implies that, for these $k$,
    \[
    f_{\ell(i),i}(e(t_{k+1}))>0\text{ and }f_{i,\ell(i)}(e(t_{k+1}))<0.
    \]
    Let $i\in I$ be fixed. Since the sequence $(\ell^p(i))_{p\geq 0}$, where $\ell^p$ is the $p$-th iterate of $\ell$, takes only finitely many values, there exists $i_1\in I$ and
    $n\geq 2$ such that $i_1=\ell^{n+1}(i_1)$ and $\ell^p(i_1)\neq i_1$ for all $p\leq n$. We denote by $i_p:=\ell^p(i_1)$ for
    $p\leq n$. As in the proof of Proposition~\ref{prop:S=infty-d=2}, the sequence $(e(t_k))_k$ converges to some $e^\star$, so we
    deduce that for all $1\leq \ell\leq n$
    \[
    f_{i_{\ell+1},i_\ell}(e^\star)\geq 0\text{ and }f_{i_{\ell},i_{\ell+1}}(e^\star)\leq 0.
    \]
    This contradicts Assumption~\ref{hyp:no-heteroclinic-loop}, hence Proposition~\ref{prop:S-infty} is proved.
\end{proof}

\section*{Acknowledgements}
This work was supported by Inria Associate Team aStoNiche. N.C.\ is funded by the European Union (ERC, SINGER, 101054787 and ERC-2024-COG MUSEUM-101170884). Views and opinions expressed are however those of the author(s) only and do not necessarily reflect those of the European Union or the European Research Council. Neither the European Union nor the granting authority can be held responsible for them. N.C.\ and C.F.\ are supported by the Chair ``Modélisation Mathématique et Biodiversité" of Veolia Environnement-Ecole Polytechnique-Museum National d’Histoire Naturelle-Fondation X. L.V. has been  supported by ANID through FONDECYT Iniciaci\'on, project number 11240158-2024 \textit{Adaptive behavior in stochastic population dynamics and non-linear Markov processes in ecoevolutionary modeling}. N.Z. is supported by ANID through FONDECYT Postdoctorado, project number 3260842. PM acknowledges support from Exploration Grant 13220168, from The Agencia Nacional de Investigación y Desarrollo (ANID), support from Centro de Modelamiento Matem\'atico (CMM), Grant FB210005, BASAL funds for centers of excellence from ANID-Chile, and support from the International Center for Theoretical Physics (ICTP) through their Associates Programme funded by the Simons Foundation through grant number 284558FY19.



\bibliographystyle{abbrv}
\bibliography{biblio_niche_construction}

\end{document}